\documentclass[showpacs,preprintnumbers,amsmath,amssymb,12pt]{article}
\usepackage{graphicx,amsfonts,amssymb,amsthm,amsmath,graphics,epsfig,pstricks,fancyhdr,fancybox}
\usepackage{dcolumn}
\usepackage{bm}

\newtheorem{thm}{Teorema}[section]

\newtheorem{prop}[thm]{Proposition}
\newtheorem{defn}[thm]{Definition}

\newcommand{\ha}{\widehat{a}\,}
\newcommand{\hb}{\widehat{b}\,}
\newcommand{\sdiff}{\emph{sdif}}
\newcommand{\sde}{\emph{SDE}}
\newcommand{\ode}{\emph{ODE}}
\newcommand{\pde}{\emph{PDE}}
\newcommand{\qho}{\emph{QHO}}
\newcommand{\fpe}{\emph{FPE}}
\newcommand{\sm}{stochastic mechanics}

\newcommand{\refeq}[1]{~(\ref{#1})}
\newcommand{\myref}[1]{~\ref{#1}}
\newcommand{\mycite}[1]{~\cite{#1}}

\newcommand{\R}{\mathbb{R}}

\newcommand{\RE}{\bm{R}}
\newcommand{\rv}{\textit{rv}}

\newcommand{\ac}{\textit{ac}}

\newcommand{\pdf}{\textit{pdf}}
\newcommand{\cdf}{\textit{cdf}}
\newcommand{\qf}{\textit{qf}}

\newcommand{\Pqo}{\bm{P}\hbox{-\emph{a.s.}}}

\newcommand{\PR}[1]{\bm{P}\left\{{#1}\right\}}
\newcommand{\EXP}[1]{\bm{E}\left[{#1}\right]}
\newcommand{\VAR}[1]{\bm{V}\left[{#1}\right]}
\newcommand{\cov}[2]{\bm{cov}\left[{#1},{#2}\right]}
\newcommand{\MED}[1]{\bm{M}\left[{#1}\right]}
\newcommand{\Lop}[1]{\mathcal{L}\left[{#1}\right]}

\newcommand{\norm}{\mathfrak{N}}
\newcommand{\lnorm}{\mathfrak{lnN}}

\newcommand{\gam}{\mathfrak{G}}

\newcommand{\ou}{\emph{OU}}

\newcommand{\fvel}{\overset{\rightarrow}{v}}
\newcommand{\bvel}{\overset{\leftarrow}{v}}

\newcommand{\indim}{\noindent{\bf Proof:}\hspace{0.2cm}}
\newcommand{\findim}{\hfill$\blacksquare$\vspace{0.3cm}\noindent}

\def\ito{It\={o}}

\begin{document}
\thispagestyle{empty}

\title{\Huge \textbf{Gompertz and logistic \\ stochastic dynamics: \\ Advances in an ongoing quest}
\author{
%Andrea \textsc{Andrisani}\footnote{deflema@yahoo.it}\\
%Department of \textsl{Mathematics}, University of Bari \\
%\vspace{7pt} via E. Orabona 4, 70125 Bari, Italy \\
Nicola \textsc{Cufaro Petroni}\footnote{cufaro@ba.infn.it}  \\
Department of \textsl{Mathematics} and \textsl{TIRES}, University of Bari\\
\textsl{INFN} Sezione di Bari\\ \vspace{7pt}
via E. Orabona 4, 70125 Bari, Italy\\
Salvatore \textsc{De Martino}\footnote{demartino@sa.infn.it}\\
Full prof. of General Physics (Ret), University of Salerno\\
\vspace{7pt}
personal address: via Luigi Settembrini 61, 80139 Napoli, Italy\\
Silvio \textsc{De Siena}\footnote{silvio.desiena@gmail.com}\\
Full prof. of Theoretical Physics (Ret), University of Salerno\\
%\vspace{5pt}
personal address: via Bastioni 15, 84122 Salerno, Italy } }

\date{}

\maketitle

\begin{abstract}
\noindent In this report we summarize a few methods for solving the
stochastic differential equations (\sde) and the corresponding
Fokker-Planck equations describing the Gompertz and logistic random
dynamics. It is shown that the solutions of the Gompertz \sde\ are
completely known, while for the logistic \sde's we can provide the
solution as an explicit process, but we find much harder to write
down its distribution in closed form. Many details of possible ways
out of this maze are listed in the paper and its appendices. We also
briefly discuss the prospects of performing a suitable averaging, or
a deterministic limit. The possibility is also suggested of
associating these equations to the stochastic mechanics of a quantum
harmonic oscillator adopted as a tool serviceable also in the field
of stochastic control: in particular we propose to investigate the
equations associated to the quantum stationary states
\end{abstract}

%\vspace{5pt}

%\noindent \emph{PACS}:

\vspace{5pt}


%\newpage

%\begin{center}

%\end{center}

\newpage

\tableofcontents

%\newpage

\section{Introduction}

%\section{Gompertz and Logistic random evolutions}

In the present paper we will mainly deal with two classes of
one-dimensional \sde's (\emph{stochastic differential equations},
see Appendix\myref{theito}) with a non linear drift that have been
recently discussed in a number of papers (see for
instance\mycite{sa,skiad,pasq} and references quoted therein) and
that are the noisy version of their deterministic counterparts (see
Appendix\myref{determ}): the Gompertz and the logistic equations.
The general form of a \textbf{\emph{Gompertz} \sde}\ for the process
$\xi(s)$ is
\begin{equation}\label{gomp}
    d\xi(s)=\big[a_1\xi(s)-a_2\,\xi(s)\ln \left(b_1\xi(s)\right)\big]\,ds+b_1\xi(s)\,d\beta(s)
\end{equation}
where $\beta(s)$ is a Wiener process with diffusion coefficient
$2\delta$, namely $\EXP{\beta^2(s)}=2\delta s$. Remark that a
possible difference between the two coefficients $b_1$ -- that in
the logarithm argument, and the other
 in front of the Wiener process -- can be easily reabsorbed by
redefining $a_1$ and $a_2$. The physical dimensions (by supposing
for instance that $\xi$ and $\beta$ are lengths $L$, while $s$ is a
time $T$) are then
\begin{equation*}
    [a_1]=\frac{1}{T}\qquad\quad[a_2]=\frac{1}{T}\qquad\quad[b_1]=\frac{1}{L}\qquad\quad[\delta]=\frac{L^2}{T}
\end{equation*}
Within a similar notation the \textbf{\emph{logistic} \sde}\ is
\begin{equation}\label{logist}
    d\xi(s)=\xi(s)\left[a_1-a_2\,\xi(s)\right]ds+b_1\xi(s)d\beta(s)
\end{equation}
with the following physical dimensions
\begin{equation*}
    [a_1]=\frac{1}{T}\qquad\quad[a_2]=\frac{1}{LT}\qquad\quad[b_1]=\frac{1}{L}\qquad\quad[\delta]=\frac{L^2}{T}
\end{equation*}
while the generalized $\bm\theta$\textbf{\emph{-logistic} \sde}\
($\theta>0$) is
\begin{equation}\label{thetalog}
    d\xi(s)=\xi(s)\left[a_1-a_2\,\xi^{\theta}(s)\right]ds+b_1\xi(s)d\beta(s)
\end{equation}
with the physical dimensions
\begin{equation*}
    [a_1]=\frac{1}{T}\qquad\quad[a_2]=\frac{1}{L^\theta T}\qquad\quad[b_1]=\frac{1}{L}\qquad\quad[\delta]=\frac{L^2}{T}
\end{equation*}
The logistic\refeq{logist} is apparently recovered for $\theta=1$.
Therefore -- within the notations of the Appendix\myref{theito} --
the time independent coefficients of the Gompertz \sde\ are
\begin{equation}\label{gompcoeff}
    a(x)=a_1x-a_2\,x\ln(b_1x)\qquad\quad b(x)=b_1x
\end{equation}
those of the logistic \sde\ are
\begin{equation}\label{logistcoeff}
    a(x)=a_1x-a_2\,x^2\qquad\quad b(x)=b_1x
\end{equation}
while for the $\theta$-logistic we finally have
\begin{equation}\label{thetalogcoeff}
    a(x)=a_1x-a_2\,x^{1+\theta}\qquad\quad
    b(x)=b_1x\qquad\quad\theta>0
\end{equation}

In order to simplify their look, however, it is expedient to recast
these equation in a dimensionless form so that only the essential
parameters will remain in evidence. As for the equation\refeq{gomp},
a transformation to the dimensionless quantities
\begin{equation*}
    t=a_1s\qquad
    X(t)=b_1\xi\left(\frac{t}{a_1}\right)\qquad
    W(t)=b_1\beta\left(\frac{t}{a_1}\right)\qquad
    D=\frac{b_1^2\delta}{a_1}\qquad\alpha=\frac{a_2}{a_1}
\end{equation*}
would give rise to the \textbf{\emph{dimensionless Gompertz} \sde}
\begin{equation}\label{gompdl}
    dX(t)=X(t)\big[1-\alpha\ln X(t)\big]dt+X(t)\,dW(t)
\end{equation}
with coefficients
\begin{equation}\label{gompcoeffdl}
    a(x)=x(1-\alpha\ln x)\qquad\quad
    b(x)=x\qquad\quad\EXP{W^2(t)}=2Dt
\end{equation}
while with the similar transformations
\begin{equation*}
    t=a_1s\qquad\quad
    X(t)=\frac{a_2}{a_1}\,\xi\left(\frac{t}{a_1}\right)\qquad\quad
    W(t)=b_1\beta\left(\frac{t}{a_1}\right)\qquad\quad
    D=\frac{b_1^2\delta}{a_1}
\end{equation*}
the equation\refeq{logist} becomes the \textbf{\emph{dimensionless
logistic} \sde}
\begin{equation}\label{logistdl}
    dX(t)=X(t)\big[1-X(t)\big]dt+X(t)\,dW(t)
\end{equation}
with coefficients
\begin{equation}\label{logistcoeffdl}
    a(x)=x(1-x)\qquad\quad b(x)=x\qquad\quad\EXP{W^2(t)}=2Dt
\end{equation}
Finally the transformations
\begin{equation*}
    t=a_1s\qquad\quad
    X(t)=\left(\frac{a_2}{a_1}\right)^{1/\theta}\xi\left(\frac{t}{a_1}\right)\qquad\quad
    W(t)=b_1\beta\left(\frac{t}{a_1}\right)\qquad\quad
    D=\frac{b_1^2\delta}{a_1}
\end{equation*}
will give rise to the \textbf{\emph{dimensionless
$\bm\theta$-logistic} \sde}
\begin{equation}\label{thetalogdl}
    dX(t)=X(t)\big[1-X(t)^\theta\big]dt+X(t)\,dW(t)
\end{equation}
with coefficients
\begin{equation}\label{thetalogcoeffdl}
    a(x)=x(1-x^\theta)\qquad\quad b(x)=x\qquad\quad\EXP{W^2(t)}=2Dt
\end{equation}
We will adopt these dimensionless notations all along the presente
paper: in the Section\myref{statsol1} we will give a virtually
complete solution of the Gompertz \sde\refeq{gompdl}, while in the
subsequent Section\myref{po} these solutions will be extended to the
parametric Gompertz equation with a time-dependent drift
coefficient. As for the logistic and $\theta$-logistic
\sde's\refeq{logistdl} and\refeq{thetalogdl} the
Section\myref{linearize} will list several partial results, in
particular the explicit solutions\refeq{logisttranssol}
and\refeq{thetalogtranssol} in the guise of processes whose
distributions however -- albeit derivable from the existing
literature\mycite{yor} -- can not be easily presented in a
manageable closed form. The same problem is addressed again in the
Section\myref{intgeowien} in the reduced form of the distribution of
the integrals of geometric Wiener processes, or even of their finite
sums at different times, but here too the results are only
preliminary while a discussion of the exact results is postponed to
a forthcoming paper

The Gompertz and logistic random dynamics will also be looked at
from the standpoint of the  \fpe's (\emph{Fokker-Planck equations},
see Appendix\myref{thefpe}) for their \pdf\ (\emph{probability
density functions}) and this will give in the future the opportunity
of exploring a further perspective. In a few previous
papers\mycite{pla,jpha} we indeed analyzed the solutions of the
\fpe's associated by the \sm\ to the quantum wave functions (see
Appendix\myref{stochmech} for the particular case of the stationary
states of a \qho, \emph{quantum harmonic oscillator}), and we looked
into the possibility of controlling the stochastic evolution by
means of suitable potentials. We plan therefore to extend in the
near future this analysis to the Gompertz and logistic random
dynamics. In particular we will focus our attention on the
relationship between these equations for the stationary states of a
\qho\ with frequency $\omega$, and the Gompertz and logistic
equations: this will be listed among other suggestions for future
research in the conclusive Section\myref{concl}. In a last
Appendix\myref{quantmed} about quantiles and medians, definitions
and results are finally collected to serve in scattered discussions
about a suitable coarse-graining of our stochastic equations: this
last step is proposed here in order to recover the deterministic
equations of the Appendix\myref{determ}, and therefore to show the
global consistency of these models

\section{Gompertz stochastic equations}\label{statsol1}

\subsection{Smoluchowsky \sde: stationary \pdf}\label{smol1}

We begin implementing first the transformations presented in the
Appendix\myref{constcoeff} ane leading to a Smoluchowsky \sde: since
our coefficients\refeq{gompcoeffdl} are time-independent, the
transformation\refeq{b1} which is now
\begin{equation*}
   y=h(x)=\ln x \qquad x=g(y)=e^y \qquad\quad Y(t)=\ln X(t)\qquad X(t)=e^{Y(t)}
\end{equation*}
leads to $\hb (y,t)=1$, while\refeq{a1} gives the drift coefficient
\begin{equation}\label{gompstat}
      \ha (y)=1-D-\alpha y
\end{equation}
Therefore, provided that $\alpha >0$, the Gompertz
\sde\refeq{gompdl} becomes a Smoluchowsky \sde\ that essentially
turns out to be an Ornstein-Uhlenbeck (\ou) \sde\ with an additional
constant drift
\begin{equation}\label{gompsmol}
    dY(t)=\big(1-D-\alpha Y(t)\big)\,dt+dW(t)
\end{equation}
This \ou\ \sde\ for $Y(t)$ can be completely solved and its Gaussian
transition \pdf\ is well known: as a consequence we can also easily
find the log-normal transition \pdf\ for the process $X(t)$, but we
will postpone to the next section a discussion of these details
taking a look for the time being only at the stationary solution of
the transformed $Y(t)$ process: with a dimensionless potential
$\chi(y)$, from\refeq{potential} and\refeq{gompstat} we first have
\begin{equation}\label{potentialdl}
    -D\,\chi'(y)=\ha(y)=1-D-\alpha y\qquad\quad\chi(y)=\frac{\phi(y)}{kT}
\end{equation}
and therefore
\begin{equation*}
    \chi(y)=\frac{\alpha}{2D}\,y^2-\frac{1-D}{D}\,y+c
\end{equation*}
so that, if $\alpha>0$ and if the integration constant $c$ is
suitably chosen, the stationary Boltzmann distribution comes out to
be a Gaussian law
$\norm\left(\frac{1-D}{\alpha},\frac{D}{\alpha}\right)$ as it is for
every \ou\ process; the original process $X(t)=e^{Y(t)}$ then is a
\emph{geometric} \ou\ and its stationary solutions are log-normal

%\subsection{Space-independent coefficients}
\subsection{Linearized \sde}

We consider next the transformation presented in the
Appendix\myref{spindep} and leading to process independent
coefficients, and we check first whether the compatibility
condition\refeq{compat} holds for our Gompertz \sde: an answer in
the affirmative follows from a direct calculation since
\begin{equation*}
    b(x)\left[D\,b''(x)-\frac{d}{dx}\left(\frac{a(x)}{b(x)}\right)\right]=\alpha
\end{equation*}
As a consequence, with $c=\alpha$, we first have from\refeq{bt} that
$\hb (t)=e^{\alpha t}$, then from\refeq{ht} it is
\begin{equation*}
    h(x,t)=e^{\alpha t}\int\frac{1}{x}\,dx=e^{\alpha t}\ln x
\end{equation*}
and finally from\refeq{at} we find that
\begin{equation*}
    \ha (t)=e^{\alpha t}\left(\alpha \ln x+\frac{x-\alpha\, x\ln x}{x}-D\right)=(1-D)\,e^{\alpha t}
\end{equation*}
The transformed Gompertz \sde\ for the process
\begin{equation*}
    Z(t)=e^{\alpha t}\ln X(t)
\end{equation*}
is then
\begin{equation}\label{gompt}
    dZ(t)=(1-D)\,e^{\alpha t}dt+e^{\alpha t}dW(t)
\end{equation}
whose solution, with $Z(0)=Z_0=\ln X_0$, from\refeq{sdetsol} is
\begin{equation*}
    Z(t)=Z_0+(1-D)\frac{e^{\alpha t}-1}{\alpha}+\int_0^te^{\alpha
    u}dW(u)
\end{equation*}
so that the solution of the \sde\refeq{gompdl} is
\begin{equation}\label{gompproc}
    X(t)=e^{e^{-\alpha t}Z(t)}=X_0^{e^{-\alpha t}}e^{(1-D)(1-e^{-\alpha
    t})/\alpha}\,\,e^{\,\int_0^t e^{-\alpha(t-u)}dW(u)}
\end{equation}
It is interesting to remark moreover that, with a degenerate initial
condition $X_0=x_0,\:\Pqo$ and by switching off the Wiener noise
($D=0$) this solution exactly coincides with the
solution\refeq{gompdet} of the deterministic Gompertz \emph{ODE}
discussed in the Appendix\myref{detgomp}. It is easy to see on the
other hand that for the process $Y(t)=e^{-\alpha t}Z(t)=\ln X(t)$ we
also have
\begin{equation*}
    dZ(t)=\alpha e^{\alpha t}Y(t)dt+e^{\alpha t}dY(t)
\end{equation*}
so that by comparing it with\refeq{gompt} we get
\begin{equation*}
    dY(t)=\left(1-D-\alpha Y(t)\right)dt+dW(t)
\end{equation*}
apparently coincident with the Smoluchowsky equation\refeq{gompsmol}
already discussed in the Section\myref{smol1}. Since moreover it is
easy to see that
\begin{equation*}
    \int_s^t\ha(u)du=(1-D)e^{\alpha s}\,\frac{e^{\alpha (t-s)}-1}{\alpha }\qquad\quad \int_s^t\hb\,^2(u)du=e^{2\alpha s}\,\frac{e^{2\alpha (t-s)}-1}{2\alpha }
\end{equation*}
from\refeq{sdetsol} and \refeq{sdettrans} we have for the solution
of\refeq{gompt} with $Z(s)=ze^{\alpha s},\;\;\Pqo$ at $t=s$
\begin{eqnarray*}
    Z(t)&=&ze^{\alpha s}+(1-D)e^{\alpha s}\,\frac{e^{\alpha (t-s)}-1}{\alpha }+\int_s^te^{\alpha u}dW(u)\\
    &\sim&\norm\left(ze^{\alpha s}+(1-D)e^{\alpha s}\,\frac{e^{\alpha (t-s)}-1}{\alpha }\,,\,De^{2\alpha s}\,\frac{e^{2\alpha (t-s)}-1}{\alpha }\right)
\end{eqnarray*}
and for the solution of the Gompertz \sde\refeq{gompdl} with $\ln
X(s)=e^{\alpha s}\ln y=ze^{\alpha s},\;\;\Pqo$ at $t=s$
\begin{eqnarray*}
   \ln X(t)&=&e^{-\alpha t}Z(t)\\
    &=&e^{-\alpha (t-s)}\ln y+\frac{1-D}{\alpha}(1-e^{-\alpha (t-s)})+e^{-\alpha t}\int_s^te^{\alpha u}dW(u)\\
    &\sim&\norm\left(e^{-\alpha (t-s)}\ln y+\frac{1-D}{\alpha}(1-e^{-\alpha (t-s)})\,,\,D\,\frac{1-e^{-2\alpha (t-s)}}{\alpha}\right)
\end{eqnarray*}
Therefore the \textbf{\emph{transition} \pdf}\ $f(x,t|y,s)$ of the
Gompertz process $X(t)$ is the log-normal law associated to the
previous Gaussian distribution which also asymptotically goes to the
stationary log-normal distribution
\begin{equation*}
    \lnorm\left(\frac{1-D}{\alpha}\,,\,\frac{D}{\alpha}\right)
\end{equation*}
already found in the section\myref{smol1}, with expectation
\begin{equation*}
    \EXP{X(+\infty)}=e^{\frac{2-D}{2\alpha}}
\end{equation*}
corresponding to the asymptotic value $e^{1/\alpha}$ of the
deterministic Gompertz equation when $D=0$ (see
Appendix\myref{detgomp})

\subsection{Parametric equations}\label{po}

In a generalization of the previous investigations, let us consider
next the so called \textbf{\emph{parametric} \ou\ \sde}\ (with a
time dependent drift velocity) for a process $\eta(s)$
\begin{equation*}
    d\eta(s)=-\omega a(s)\eta(s)\,ds+b\,d\beta(s)
\end{equation*}
where $\beta(s)$ again is a Wiener process with diffusion
coefficient $2\delta$ so that $\EXP{\beta^2(s)}=2\delta s$. Here
$\omega$ is a frequency, $b$ a dimensionless constant and $a(s)$ a
dimensionless function of the time $s$. Since
$\sigma=\sqrt{^\delta/_\omega}$ has the same physical dimensions of
$\eta$ and $\beta$ (a length, for instance), we can now switch to a
dimensionless formulation by taking $t=\omega s$ and
\begin{equation*}
   Y(t)=\frac{1}{\sigma}\eta\left(\frac{t}{\omega}\right)\qquad
    W(t)=\frac{b}{\sigma}\beta\left(\frac{t}{\omega}\right)\qquad
    D=\frac{\delta b^2}{\omega\sigma^2}=b^2\qquad \alpha(t)=a\left(\frac{t}{\omega}\right)
\end{equation*}
a transformation leading to the dimensionless parametric \ou\ \sde
\begin{equation}\label{paramou}
    dY(t)=-\alpha(t)Y(t)dt+dW(t)\qquad\quad Y(0)=Y_0
\end{equation}
where $W(t)$ is a Wiener process with diffusion coefficient $2D$,
namely $\EXP{W^2(t)}=2Dt$. For further purposes $\alpha(t)$ could
also become a stochastic process, but for the time being we will
take it just as a suitable deterministic function. If we now define
the new process $X(t)$ as
\begin{equation*}
    X(t)=e^{Y(t)}\qquad\quad Y(t)=\ln X(t)
\end{equation*}
taking $a(y,t)=-\alpha(t)y$ and $b(y,t)=1$, with $x=h(y)=e^y$ and
$y=g(x)=\ln x$, from\refeq{itoy1} and\refeq{itoy2} -- swapping the
symbols $x,y$  -- we obtain the coefficients for the transformed
equation
\begin{eqnarray*}
  \ha(x,t) &=& \big[e^y\big(1-\alpha(t)y\big)\big]_{y=\ln x}=x\big(1-\alpha(t)\ln x\big) \\
  \hb(x,t) &=& \big[e^y\big]_{y=\ln x}=x
\end{eqnarray*}
finally leading to the \textbf{\emph{parametric Gompertz} \sde}
\begin{equation}\label{paramgomp}
   dX(t)= X(t)\big(1-\alpha(t)\ln X(t)\big)\,dt+X(t)\,dW(t)
\end{equation}
Therefore this new equation\refeq{paramgomp}, that apparently
generalizes\refeq{gompdl}, can be solved by looking first at a
solution of\refeq{paramou}: to this end we remark that, since
$a(y,t)=-\alpha(t)y$ and $b(y,t)=1$, using again\refeq{itoy1}
and\refeq{itoy2}, the transformation
\begin{equation*}
    Z(t)=Y(t)e^{\,\int_0^t\!\alpha(u)du}=h(Y(t),t)\qquad\quad h(y,t)=ye^{\,\int_0^t\!\alpha(u)du}
\end{equation*}
leads to the coefficients
\begin{equation*}
    \ha(z,t)=0\qquad\qquad\hb(z,t)=e^{\,\int\!\alpha(t)dt}
\end{equation*}
namely to the \sde
\begin{equation*}
    dZ(t)=e^{\,\int_0^t\!\alpha(u)du}dW(t)\qquad\quad Z(0)=Y_0
\end{equation*}
whose solution is
\begin{equation*}
    Z(t)=Y_0+\int_0^te^{\,\int_0^r\!\alpha(u)du}dW(r)
\end{equation*}
so that we finally have
\begin{equation}\label{paramousol}
    Y(t)=Y_0e^{-\int_0^t\alpha(u)du}+\int_0^te^{-\int_r^t\alpha(u)du}dW(r)
\end{equation}
that with a condition at a time $s$ becomes
\begin{equation}\label{paramousols}
    Y(t)=Y_0e^{-\int_s^t\alpha(u)du}+\int_s^te^{-\int_r^t\alpha(u)du}dW(r)
\end{equation}
As long as the \rv\ $Y_0\sim\norm(y_0,\sigma_0^2)$ is Gaussian, the
process\refeq{paramousols} is Gaussian too, and if $Y_0$ is
degenerate the \pdf\ of $Y(t)$ also is the \emph{transition} \pdf.
On the other hand when\refeq{paramousols} is Gaussian their laws are
completely determined by the expectations and the covariance
functions that can be explicitly calculated. In fact we first have
\begin{equation}\label{xmean}
    \EXP{Y(t)}=y_0\,e^{-\int_s^t\alpha(u)du}
\end{equation}
then, by taking $\widetilde{Y}_0=Y_0-y_0$ and
\begin{equation*}
    \widetilde{Y}(t)=Y(t)-\EXP{Y(t)}=\widetilde{Y}_0\,e^{-\int_s^t\alpha(u)du}+\int_s^te^{-\int_r^t\alpha(u)du}dW(r)
\end{equation*}
from the independence of $X_0$ and $W(t)$ and from the usual
properties\refeq{wiendiff} of the increments $dW(t)$, we get for
$t<t'$
\begin{eqnarray*}
    \cov{Y(t)}{Y(t')}&=&\EXP{\widetilde{Y}(t)\widetilde{Y}(t')}\\
    &=&\sigma_0^2\,e^{-\int_s^t\alpha(u)du-\int_s^{t'}\alpha(u)du}\\
    &&\qquad\quad+\EXP{\int_s^te^{-\int_r^t\alpha(u)du}dW(r)\int_s^{t'}e^{-\int_{r'}^{t'}\alpha(u)du}dW(r')}\\
    &=&\sigma_0^2\,e^{-\int_s^t\alpha(u)du-\int_s^{t'}\alpha(u)du}\\
    &&\qquad\quad+\EXP{\int_s^te^{-\int_r^t\alpha(u)du}dW(r)\int_s^te^{-\int_{r'}^{t'}\alpha(u)du}dW(r')}\\
    &=&\sigma_0^2\,e^{-\int_s^t\alpha(u)du-\int_s^{t'}\alpha(u)du}+2D\int_s^te^{-\int_r^t\alpha(u)du-\int_r^{t'}\alpha(u)du}dr
\end{eqnarray*}
and hence in any case
\begin{equation}\label{paramcovar}
    \cov{Y(t)}{Y(t')}=\sigma_0^2\,e^{-\int_s^t\alpha(u)du-\int_s^{t'}\alpha(u)du}+2D\int_s^{t\wedge t'}e^{-\int_r^t\alpha(u)du-\int_r^{t'}\alpha(u)du}dr
\end{equation}
By the way remark that the results\refeq{xmean}
and\refeq{paramcovar} hold even if $Y(t)$ is not Gaussian, but when
$Y_0\sim\norm(y_0,\sigma_0^2)$ they also completely determines the
distribution of the process\refeq{paramousols} and in particular
\begin{equation}\label{parampdf}
    Y(t)\sim\norm\left(y_0\,e^{-\int_s^t\alpha(u)du}\,,\,\sigma_0^2\,e^{-\int_s^t\alpha(u)du}+2D\int_s^te^{-2\int_r^t\alpha(u)du}dr\right)
\end{equation}
With $\sigma_0=0$ (degenerate initial condition) we moreover have
\begin{equation}\label{paramtrans}
    Y(t)\sim\norm\left(y_0\,e^{-\int_s^t\alpha(u)du}\,,\,2D\int_s^te^{-2\int_r^t\alpha(u)du}dr\right)
\end{equation}
which plays the role of the transition \pdf\ for the processes
solution of\refeq{paramou}.

Going back now to the parametric Gompetrz \sde\refeq{paramgomp}, we
find that when the parametric \ou\ process $Y(t)$ is Gaussian, then
the process $X(t)$ is log-normal and we can explicitly give all the
details of its distribution, in particular from\refeq{paramtrans}
and with $X(s)=y=e^{y_0}$ its transition \pdf\ $f(x,t|y,s)$ is
\begin{equation}\label{paramgomptrans}
    X(t)\sim\lnorm\left(e^{-\int_s^t\alpha(u)du}\ln y\,,\,2D\int_s^te^{-2\int_r^t\alpha(u)du}dr\right)
\end{equation}
It is possible then to calculate also its expectation and variance
according to\refeq{lnorm}
\begin{equation*}
  \EXP{X(t)} = e^{\EXP{Y(t)}+\VAR{Y(t)}/2}\qquad\quad
  \VAR{X(t)} = e^{2\EXP{Y(t)}+\VAR{Y(t)}}\left(e^{\VAR{Y(t)}}-1\right)
\end{equation*}
that in the case\refeq{paramgomptrans} of degenerate initial
conditions become
\begin{eqnarray}
  \EXP{X(t)} &=& y^{\,e^{-\int_s^t\alpha(u)du}}e^{D\int_s^te^{-2\int_r^t\alpha(u)du}dr}\label{lognormmean} \\
  \VAR{X(t)} &=&
  y^{\,2e^{-\int_s^t\alpha(u)du}}e^{2D\int_s^te^{-2\int_r^t\alpha(u)du}dr}\left(e^{2D\int_s^te^{-2\int_r^t\alpha(u)du}dr}-1\right)\label{lognormvar}
\end{eqnarray}
As for the median, when $Y(t)$ is Gaussian, from\refeq{med} we
simply get
\begin{equation}\label{lognormmedian}
        \MED{X(t)}=e^{\EXP{Y(t)}}=y^{\,e^{-\int_s^t\alpha(u)du}}
\end{equation}
that can be used along with\refeq{lognormmean} and\refeq{lognormvar}
to analyze the oscillations of the system. However, as hinted also
in the Appendix\myref{medians}, we must keep into account that these
results about means and medians are not completely general since
they hold only in so far as the parametric \ou\ process $Y(t)$ is
Gaussian: this is not always the case because it requires that the
initial condition itself should be Gaussian. We could alternatively
take advantage of the more general relation\refeq{medianexp}, but
also in this case there is a snag because we should calculate the
median $\MED{Y(t)}$ which is not always a straightforward job for
arbitrary initial conditions. The adoption of medians as a way to
retrieve the deterministic evolution has been extensively discussed
in a few previous papers\mycite{sa}

\section{Logistic stochastic equations}\label{linearize}

\subsection{Smoluchowsky \sde: stationary \pdf}

The transformation\refeq{b1} for the
coefficients\refeq{logistcoeffdl}, namely
\begin{equation*}
   y=h(x)=\ln x \qquad x=g(y)=e^y \qquad Y(t)=\ln X(t)\qquad X(t)=e^{Y(t)}
\end{equation*}
applied to the logistic \sde\refeq{logistdl} leads to $\hb (y,t)=1$,
and from\refeq{a1} to the drift coefficient
\begin{equation*}
      \ha (y)=1-D-e^y
\end{equation*}
namely to the Smoluchowsky \sde
\begin{equation*}
    dY(t)=\left(1-D-e^{Y(t)}\right)\,dt+dW(t)
\end{equation*}
and from\refeq{potential} to a dimensionless potential
\begin{equation*}
    \chi(y)=\frac{\phi(y)}{kT}=\frac{e^y}{D}-\frac{1-D}{D}\,y+c
\end{equation*}
that -- provided now that $1>D$ -- gives rise to the following
stationary log-gamma Boltzmann distribution (see\mycite{grad} 3.328
for the normalization integral)
\begin{equation*}
    \frac{e^{-\frac{e^y}{D}+\frac{1-D}{D}\,y}}{D^{\frac{1-D}{D}}\,\Gamma\left(\frac{1-D}{D}\right)}\qquad\quad1>D
\end{equation*}
We will not elaborate further about this stationary distribution for
the process $Y(t)$, and we will rather confine ourselves to remark
that by transforming back to the original process $X(t)=e^{Y(t)}$ we
find that its stationary density is
\begin{equation*}
    \frac{1}{D\,\Gamma\left(\frac{1-D}{D}\right)}\,\left(\frac{x}{D}\right)^{\frac{1-D}{D}-1}\,e^{-\frac{x}{D}}\qquad\quad
    1>D\qquad x>0
\end{equation*}
namely that the stationary distribution of $X(t)$ is the gamma law
$\gam\left(\frac{1-D}{D},\frac{1}{D}\right)$. There is no easy way
instead at this stage to find the form of the transition \pdf\ for
both the processes $Y(T)$ and $X(t)$. Remark finally that the
condition $1>D$, required to have a normalizable stationary
solution, amounts to the explicit condition
\begin{equation*}
    a_1>b_1^2\delta
\end{equation*}
among the coefficients of the original logistic \sde\refeq{logist}
laden with its physical dimensions, and hence it represent an
equilibrium condition between the dynamical and the diffusive
components of the process

\subsection{Linearized \sde}\label{linsde}

Considering first to the transformation to process-independent new
coefficients discussed in the Appendix\myref{spindep}, since the
coefficients\refeq{logistcoeffdl} are time-independent, we must
preliminarily check the compatibility condition\refeq{compat}, but
we find
\begin{equation*}
    b(x)\left[D\,b''(x)-\frac{d}{dx}\left(\frac{a(x)}{b(x)}\right)\right]=x
\end{equation*}
that is not constant, and hence the said compatibility
condition\refeq{compat} is not satisfied by the logistic
coefficients

As next step we then explore the possibility of linearizing the
\sde\refeq{logistdl} in the sense discussed in the
Appendix\myref{linapp}: in order to check the
condition\refeq{lincond} we find from the
coefficients\refeq{logistcoeffdl} that
\begin{equation*}
    q(x)=1-D-x\qquad\quad
    b(x)q'(x)=-x\qquad\quad
    \frac{1}{q'(x)}\frac{d}{dx}\left[b(x)q'(x)\right]=1
\end{equation*}
so that the compatibility condition\refeq{lincond} is satisfied and
from\refeq{tratolin1} and\refeq{lincond2} we also find
\begin{equation*}
    \hb_1=-1\qquad p(x)=\int\frac{dx}{b(x)}=\ln x\qquad\quad h(x)=c\,e^{-p(x)}=c\,e^{-\ln x}=\frac{c}{x}
\end{equation*}
The reciprocal transformation relations then are
\begin{equation*}
    Y(t)=\frac{c}{X(t)}\qquad\qquad X(t)=\frac{c}{Y(t)}
\end{equation*}
and if we choose $c=1$ as integration constant we get the
transformation
\begin{equation*}
    y=h(x)=\frac{1}{x}\qquad\quad
    x=g(y)=\frac{1}{y}\qquad\quad
    Y(t)=\frac{1}{X(t)}\qquad\quad X(t)=\frac{1}{Y(t)}
\end{equation*}
so that with
\begin{equation*}
    h'(x)=-\frac{1}{x^2}\qquad\qquad h''(x)=\frac{2}{x^3}
\end{equation*}
from\refeq{itoy1} and\refeq{itoy2} we have
\begin{eqnarray*}
  \ha(y) &=& h'(g(y))a(g(y))+D\,h''(g(y))b^2(g(y))=(2D-1)y+1 \\
  \hb(y) &=& h'(g(y))b(g(y))=-y
\end{eqnarray*}
namely from\refeq{linearizedcoeff}
\begin{equation*}
    \ha_0=1\qquad\quad\ha_1=2D-1\qquad\qquad\hb_0=0\qquad\quad\hb_1=-1
\end{equation*}
and hence the new \sde\ is
\begin{equation}\label{logistdlt}
    dY(t)=  \big[(2D-1)Y(t)+1\big]\,dt-Y(t)\,dW(t)
\end{equation}
As a consequence, by taking as in\refeq{linearhsolbar1}
\begin{equation*}
    \overline{Z}(t)=(D-1)t-W(t)\sim\norm\big((D-1)t\,,\,2Dt\big)
\end{equation*}
the general solution\refeq{linearsolnot} of the linearized
\sde\refeq{logistdlt} for $Y(0)=Y_0$ is
\begin{equation*}
    Y(t)=e^{\overline{Z}(t)}\left(Y_0+\int_0^te^{-\overline{Z}(u)}du\right)
\end{equation*}
while the solution $X(t)=\frac{1}{Y(t)}$ of the logistic
\sde\refeq{logistdl} for $X(0)=X_0=\frac{1}{Y_0}$ is
\begin{equation}\label{logistsol}
    X(t)=\frac{X_0e^{-\overline{Z}(t)}}{1+X_0\int_0^te^{-\overline{Z}(u)}du}
\end{equation}
Remark that with a degenerate initial condition $X_0=x_0,\:\Pqo$ and
by switching off the Wiener noise ($D=0$) we get
$\overline{Z}(t)=-t$, and the solution\refeq{logistsol} exactly
coincides with the solution\refeq{logistdet} of the deterministic
logistic \emph{ODE} discussed in the Appendix\myref{detlogist}. By
taking on the other hand $X_0=y,\;\Pqo$ at a time $0\leq s\leq t$ we
have the solution
\begin{equation}\label{logisttranssol}
    X(t)=\frac{y\,e^{-\overline{Z}(t-s)}}{1+y\int_s^te^{-\overline{Z}(u)}du}
\end{equation}
whose \pdf\ $f(x,t|y,s)$ will be the \textbf{\emph{transition}
\pdf}\ of our logistic process. If moreover we define the derivable
process
\begin{equation}\label{Aproc}
    A(t)=X_0\int_0^te^{-\overline{Z}(u)}du\qquad\qquad\dot A(t)=X_0\,e^{-\overline{Z}(t)}
\end{equation}
the solution\refeq{logistsol} takes the equivalent forms (see
also\mycite{grad} $3.434.2$)
\begin{eqnarray}\label{logistsol2}
    X(t)&=&\frac{\dot A(t)}{1+A(t)}=\frac{d}{dt}\ln\big(1+A(t)\big)=\frac{d}{dt}\int_0^\infty e^{-u}\frac{1-e^{-uA(t)}}{u}\,du\nonumber\\
    &=&-\int_0^\infty\frac{e^{-u}}{u}\,\frac{d}{dt}\,e^{-uA(t)}du=\dot A(t)\int_0^\infty e^{-u(1+A(t))}du
\end{eqnarray}
hinting additionally to a possible direct connection between the
moments of $X(t)$ and the generating function of $A(t)$, namely
$\EXP{e^{-uA(t)}}$

At first sight these results seem to be coherent with those
elaborated in other papers\mycite{skiad}, where also a few
procedures leading to the calculation of expectations and variances
of $X(t)$ are discussed. For the time being however we will neglect
a detailed analysis of these claims, noting instead that an
explicit, exact form of the \pdf\ of\refeq{logisttranssol} could be
retrieved only by taking advantage of a few rather intricate results
available in the literature\mycite{yor}. While indeed
$\overline{Z}(t)$ is a Gaussian process (it is just a re-scaled
Wiener process plus a uniform drift) so that $e^{-\overline{Z}(t)}$
is a geometric Gaussian process with a log-normal law, it is totally
another matter to find the law of the integral process
\begin{equation*}
    \int_s^te^{-\overline{Z}(u)}du
\end{equation*}
Extensive research\mycite{yor,dufr,schrod,carr} has been devoted to
this problem, but the available answers are far from being easy to
handle (see also the subsequent discussion in the
Section\myref{intgeowien}). We will therefore present in the
following, for the time being, only a few \emph{elementary}
approaches with their associated partial results, looking forward
instead to scrutinize the question in further detail in a
forthcoming paper within the framework of a more general setting

\subsubsection{Semi-explicit transition \pdf}\label{semiexpl}

We will present first a semi-explicit form of the transition \pdf\
by recalling the notations and the results referred to in the
Appendix\myref{transpdf}. Since the
coefficients\refeq{logistcoeffdl} are time-independent we begin by
defining the functions
\begin{align*}
   & y=h(x)=\int\frac{dx}{b(x)}=\ln x\qquad\quad
   x=g(y)=e^y\\
   &
   \ha(y)=\frac{a(g(y))}{b(g(y))}-D\,b'(g(y))=1-D-e^y
\end{align*}
so that we also find
\begin{eqnarray*}
  \beta(y) &=& -\frac{\ha\,^2(y)}{4D}-\frac{\ha\,'(y)}{2}=-\frac{1}{4D}\left(1-D-e^y\right)^2+\frac{e^y}{2} \\
   &=&-\frac{1}{4D}\left(1-e^y\right)^2+\frac{1}{2}-\frac{D}{4}
\end{eqnarray*}
As a consequence we have the following expressions
\begin{eqnarray*}
  \overline{h}(r) &=& r\,\ln x+(1-r)\ln y\qquad\quad e^{\overline{h}(r)}=x^r y^{1-r}\\
  \overline{W}_{st}(r) &=& W\big(s+(t-s) r\big)-\big( r W(t)+(1-  r)
  W(s)\big)\\
  \beta\left(\overline{W}_{st}(r)+\overline{h}(r)\right)
         &=& \frac{2-D}{4}-\frac{1}{4D}\left(1-x^r y^{1-r} e^{\overline{W}_{st}(r)}\right)^2
\end{eqnarray*}
and hence we get
\begin{eqnarray*}
  Z(s,t) &=&
  \int_0^1\beta\left(\overline{W}_{st}(r)+\overline{h}(r)\right)dr \;=\; \frac{2-D}{4}-\frac{1}{4D}\int_0^1\left(1-x^r y^{1-r} e^{\overline{W}_{st}(r)}\right)^2dr\\
  &=& \frac{2-D}{4}-\frac{1}{4D}+\frac{y}{2D}\int_0^1\left(\frac{x}{y}\right)^r e^{\overline{W}_{st}(r)}dr -\frac{y^2}{4D}\int_0^1\left(\frac{x}{y}\right)^{2r} e^{2\overline{W}_{st}(r)}dr  \\
  &=& -\frac{(1-D)^2}{4D}+\frac{y}{2D}\int_0^1\left(\frac{x}{y}\right)^r e^{\overline{W}_{st}(r)}dr -\frac{y^2}{4D}\int_0^1\left(\frac{x}{y}\right)^{2r} e^{2\overline{W}_{st}(r)}dr
\end{eqnarray*}
We thus find
\begin{eqnarray*}
  e^{(t-s)Z(s,t)} &=& e^{\frac{2-D}{4}(t-s)}\exp\left\{-\frac{t-s}{4D}\int_0^1\left(1-x^r y^{1-r}
  e^{\overline{W}_{st}(r)}\right)^2dr\right\}\\
  &=& e^{-\frac{(1-D)^2}{4D}(t-s)}\\
  &&\quad\exp\left\{-\frac{t-s}{4D}\left(
            y^2\int_0^1\left(\frac{x}{y}\right)^{2r}e^{2\overline{W}_{st}(r)}dr
                                                 -2y\int_0^1\left(\frac{x}{y}\right)^r e^{\overline{W}_{st}(r)}dr\right)\right\}
\end{eqnarray*}
so that finally for the expectation factor in our transition
\pdf\refeq{pdf} we have
\begin{align*}
 &\qquad\qquad\qquad\qquad\quad \EXP{e^{(t-s)Z(s,t)}} = e^{-\frac{(1-D)^2}{4D}(t-s)}\mu(x,t;y,s)\\
 & \mu(x,t;y,s)=\EXP{\exp\left\{-\frac{t-s}{4D}\left(
            y^2\int_0^1\!\!\left(\frac{x}{y}\right)^{2r}\!\!e^{2\overline{W}_{st}(r)}dr
                                                 -2y\int_0^1\!\!\left(\frac{x}{y}\right)^r\!\! e^{\overline{W}_{st}(r)}dr\right)\right\}}
\end{align*}
As for the other factors in\refeq{pdf} we first remark that
\begin{figure}
\begin{center}
\includegraphics*[width=13cm]{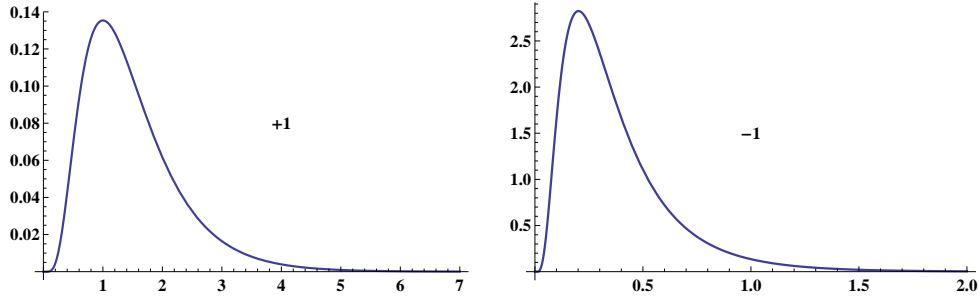}
\caption{Behavior of the function\refeq{firstterm} for both the
possible signs (beware: the plot scales in the two pictures are
rather different)}\label{pdffig}
\end{center}
\end{figure}
\begin{eqnarray*}
  \frac{1}{2D}\int_y^x\frac{a(z)}{b^2(z)}\,dz &=& \frac{1}{2D}\int_y^x\frac{1-z}{z}\,dz=\frac{1}{2D}\left(\ln\frac{x}{y}-x+y\right) \\
  -\frac{1}{4D(t-s)}\left(\int_y^x\frac{dz}{b(z)}\right)^2&=&-\frac{1}{4D(t-s)}\left(\int_y^x\frac{dz}{z}\right)^2=-\frac{1}{4D(t-s)}\ln^2\frac{x}{y}
\end{eqnarray*}
and then we find out for the transition \pdf
\begin{eqnarray}
  f(x,t|y,s) &=& \frac{e^{-\frac{(1-D)^2}{4D}(t-s)}}{x\sqrt{4\pi D(t-s)}}\sqrt{\frac{y}{x}}\,
            e^{\frac{1}{2D}\left(\ln\frac{x}{y}-x+y\right)}e^{-\frac{1}{4D(t-s)}\ln^2\frac{x}{y}}\,\mu(x,t;y,s)\nonumber\\
   &=&\frac{e^{-\frac{(1-D)^2}{4D}(t-s)-\frac{1}{2}\ln\frac{x}{y}+\frac{1}{2D}\left(\ln\frac{x}{y}-x+y\right)-\frac{1}{4D(t-s)}\left(\ln\frac{x}{y}\right)^2}}{x\sqrt{4\pi
                                     D(t-s)}}\,\mu(x,t;y,s)\nonumber\\
   &=&\frac{e^{-\frac{x-y}{2D}-\frac{1}{4D(t-s)}\left((1-D)(t-s)-\ln\frac{x}{y}\right)^2}}{x\sqrt{4\pi
   D(t-s)}}\label{logistpdf}\,\mu(x,t;y,s)
\end{eqnarray}
Here too, however, despite the fact that the transition \pdf\ is
given in closed form, the calculation of the expectation
in\refeq{logistpdf} depends on the knowledge of the law of the
integral of a geometric Gaussian process similar to that of the
solution\refeq{logisttranssol}. Therefore -- at variance with the
case of the Gompertz \sde\ -- the expression\refeq{logistpdf} seems
to represent the farthest point we can reach at present along this
path in our quest for an explicit formula of the transition \pdf\
for the logistic \sde. The main hurdle apparently is the computation
of the expectation term $\mu(x,t;y,s)$ which contains integrals of
geometric Wiener processes with non elementary
distributions\mycite{yor}. All that we can easily assess for the
time being is the behavior of the explicit term in front of $\mu$
that is of the type
\begin{equation}\label{firstterm}
    \frac{e^{-x-(\pm1-\ln x)^2}}{x}
\end{equation}
This function turns out to be regular in the origin $x=0$ for both
the possible signs and, coherently with the behavior of the
stationary distribution, it displays a gamma-like shape for $x>0$ as
can be seen from the Figure\myref{pdffig}

\subsubsection{Fokker-Planck equation}

We will next turn our attention to the possible solutions of the
corresponding \fpe\ along the lines presented in the
Appendix\myref{thefpe}: for a process X(t) solution of the logistic
\sde\refeq{logistdl} the \fpe\refeq{fpeq} is
\begin{eqnarray}
  \partial_tf_X(x,t)&=&-\partial_x\left[\fvel(x)f_X(x,t)\right]+\partial^2_x\left[B(x)f_X(x,t)\right] \nonumber\\
   &=& -\partial_x\left[x(1-x)f_X(x,t)\right]+D\,\partial^2_x\left[x^2f_X(x,t)\right] \label{logistfpx}\\
   &=&
   Dx^2\,\partial^2_xf_X(x,t)+x(4D-1+x)\partial_xf_X(x,t)+(2D-1+2x)f_X(x,t)\nonumber
\end{eqnarray}
where
\begin{equation*}
    \fvel(x)=a(x)=x(1-x)\qquad\qquad B(x)=Db^2(x)=Dx^2
\end{equation*}
while for the transformed process $Y(t)=\frac{1}{X(t)}$ solution of
the \sde\refeq{logistdlt} we have
\begin{eqnarray}
  \partial_tf_Y(y,t) &=&-\partial_y\left[\fvel(y)f_Y(y,t)\right]+\partial^2_y\left[B(y)f_Y(y,t)\right] \nonumber\\
   &=& -\partial_y\left[((2D-1)y+1)f_Y(y,t)\right]+D\,\partial^2_y\left[y^2f_Y(y,t)\right] \label{logistfpy}\\
   &=&
   Dy^2\,\partial^2_yf_Y(y,t)+[(2D+1)y-1]\partial_yf_Y(x,t)+f_Y(y,t)\nonumber
\end{eqnarray}
with
\begin{equation*}
    \fvel(y)=\ha(y)=(2D-1)y+1\qquad\qquad B(y)=D\,\hb^2(y)=Dy^2
\end{equation*}
Remark by the way that the \pdf\ of $Y(t)=\frac{1}{X(t)}$ can always
be derived from that of $X(t)$ as
\begin{equation*}
    f_Y(y,t)=\frac{1}{y^2}\,f_X\left(\frac{1}{y}\,,t\right)
\end{equation*}
so that the corresponding equations could be deduced one from the
other by means of this transformation

%\subsubsection{Eigenfunction expansions}

We can then look for the solutions with eigenfunction expansions
starting with\refeq{logistfpx}: since we already know that the \pdf\
of the gamma law $\gam\left(\frac{1-D}{D},\frac{1}{D}\right)$
\begin{equation}\label{logistat}
    \widetilde{f}_X(x)=\frac{\left(\frac{1}{D}\right)^\frac{1-D}{D}}{\Gamma\left(\frac{1-D}{D}\right)}\,x^{\frac{1-D}{D}-1}e^{-\frac{x}{D}}
    =\frac{1}{D}\,\frac{\left(\frac{x}{D}\right)^{\frac{1-D}{D}-1}}{\Gamma\left(\frac{1-D}{D}\right)}\,e^{-\frac{x}{D}}\qquad\quad
    1>D
\end{equation}
is a stationary solution of\refeq{logistfpx} (this can be also
checked by direct calculation), from the Section\myref{eigenexp}
stems that by taking
\begin{equation*}
    f_X(x,t)=\sqrt{\widetilde{f}_X(x)}\,g_X(x,t)
\end{equation*}
we get for $g_X$ the new equation
\begin{equation}\label{slx}
    \partial_tg_X(x,t)=\Lop{g_X}(x,t)
\end{equation}
where $\mathcal{L}$ is an operator of the Sturm-Liouville form
\begin{equation}\label{slop}
    \Lop{\varphi}(x)=\frac{d}{dx}\left[p(x)\frac{d\varphi(x)}{dx}\right]-q(x)\varphi(x)
\end{equation}
that is  is self-adjoint for functions satisfying suitable boundary
conditions in $x=0$ and $x=+\infty$. It can be shown (and
cross-checked by direct calculation) that for our
equation\refeq{logistfpx} we have in particular
\begin{eqnarray*}
    p(x)&=&B(x)=Dx^2\\
    q(x)&=&\frac{\left[B'(x)-\fvel(x)\right]^2}{4B(x)}-\frac{\left[B'(x)-\fvel(x)\right]'}{2}=\frac{(x-1)^2-2D}{4D}
\end{eqnarray*}
so that we finally get
\begin{equation*}
    \partial_tg_X(x,t)=Dx^2\,\partial_x^2g_X(x,t)+2Dx\,\partial_xg_X(x,t)+\frac{2D-(x-1)^2}{4D}\,g_X(x,t)
\end{equation*}
We then separate the variables by taking
\begin{equation*}
    g_X(x,t)=e^{-\lambda t}G_X(x)
\end{equation*}
obtaining the eigenvalue equation
\begin{equation*}
    \Lop{G_X}(x)+\lambda G_X(x)=0
\end{equation*}
that can be explicitly written as
\begin{equation}\label{eigeneqx}
  x^2\,G_X''(x)+2x\,G_X'(x)+\left[\frac{2D-(x-1)^2}{4D^2}+\frac{\lambda}{D}\right]G_X(x)=0
\end{equation}
Now this is a totally Fuchsian equation with two singularities in
$x=0$ and $x=+\infty$ and consequently can be treated with the usual
methods: first of all it is possible to check by direct calculation
that $\lambda_0=0$ is an eigenvalue for the eigenfunction
$G_0(x)=\sqrt{\widetilde{f}_X(x)}$. Then to simplify the notation we
change the variable according to
\begin{equation*}
    z=\frac{x}{D}\qquad\qquad
    G_X(x)=G_X\left(Dz\right)=\psi(z)
\end{equation*}
and we get
\begin{equation}\label{eigeneqx2}
    z^2\psi''(z)+2z\psi'(z)+\left[\frac{1+2\lambda}{2D}-\left(\frac{Dz-1}{2D}\right)^2\right]\psi(z)=0
\end{equation}
Now we take
\begin{eqnarray*}
  \psi(z) &=& \frac{e^{z/2}}{z}\,u(z) \\
  \psi'(z) &=& \frac{e^{z/2}}{z}\left[u'(z)+\frac{z-2}{2z}\,u(z)\right] \\
  \psi''(z) &=& \frac{e^{z/2}}{z}\left[u''(z)+\frac{z-2}{z}\,u'(z)+\frac{4+(z-2)^2}{4z^2}\,u(z)\right]
\end{eqnarray*}
and we find
\begin{equation*}
    u''(z)+u'(z)+\left[\frac{1}{2Dz}+\frac{(1+2\lambda)2D-1}{4D^2z^2}\right]u(z)=0
\end{equation*}
that can be put in the form of a confluent hypergeometric equation
(see\mycite{grad} formula $9.202.1$)
\begin{equation}\label{confluent}
    u''(z)+u'(z)+\left[\frac{\frac{1}{2D}}{z}+\frac{\frac{1}{4}-\frac{(1-D)^2-4\lambda D}{4D^2}}{z^2}\right]u(z)=0
\end{equation}
where moreover the term
\begin{equation*}
    \mu^2=\frac{(1-D)^2-4\lambda
    D}{4D^2}=\left(\frac{1-D}{2D}\right)^2-\frac{\lambda}{D}
\end{equation*}
is required to be positive, which would happen only if
\begin{equation*}
    \lambda<D\left(\frac{1-D}{2D}\right)^2
\end{equation*}
When this happens two linearly independent solutions are
(\mycite{grad} formula $9.202.2-3$)
\begin{eqnarray*}
  u_1(z) &=& z^{\frac{1}{2}+\mu}\,e^{-z}\Phi\left(\frac{D-1}{2D}+\mu,1+2\mu;z\right) \\
  u_2(z) &=& z^{\frac{1}{2}-\mu}\,e^{-z}\Phi\left(\frac{D-1}{2D}-\mu,1-2\mu;z\right)
\end{eqnarray*}
where $\Phi(\alpha,\gamma;z)$ is the confluent hypergeometric
function (see\mycite{grad} formula $9.210.1$). It is well known that
the eigenvalues are found by requiring that $\Phi(\alpha,\gamma;z)$
degenerates in a Laguerre polynomial (see\mycite{grad} formula
$8.970.1$) and that this happens when $\alpha=-n$ is a negative
integer (see\mycite{grad} formula $8.972.1$). As a consequence our
eigenvalues are selected by the requirement
\begin{equation*}
    \frac{D-1}{2D}\pm\mu=\frac{D-1}{2D}\pm\sqrt{\left(\frac{1-D}{2D}\right)^2-\frac{\lambda}{D}}=-n
\end{equation*}
namely
\begin{equation}\label{eigenvalx}
    \lambda_n=D\left(\frac{1-D}{D}\,n-n^2\right)
\end{equation}
However, while $\lambda_0=0$ is confirmed as an eigenvalue, we find
that just a finite number of eigenvalues are possibly positive, and
that they turn negative as soon as
\begin{equation*}
    n>\frac{1-D}{D}
\end{equation*}
For small $D$ this limit can be quite large, but that
notwithstanding it remains a finite number, and hence it seems
apparent that we can not have the infinite sequence of (increasing)
positive eigenvalues that we would have supposed to have: this is a
puzzling point, and these results can be also checked with a
shortcut through \textsc{Mathematica} asking for the solutions
of\refeq{eigeneqx}. Remark that the positivity of the eigenvalues is
directly linked to the ergodicity of the system, because it would
entail that all the eigenfunctions other than the stationary
solution are wiped out in time exponentially fast: failing to have
positive eigenvalues would instead present the case of exploding
terms in the eigenfunction expansion

If on the other hand we try to look as an alternative for the
eigenfunction expansion of the solution of the Fokker-Planck
equation\refeq{logistfpy} for the transformed process
$Y(t)=\frac{1}{X(t)}$ we would meet again seeming insurmountable
problems: we indeed immediately find by direct calculation that the
stationary solution is now, quite understandably, the inverse gamma
law
$\hbox{\emph{\textbf{inv}}-}\gam\left(\frac{1-D}{D},\frac{1}{D}\right)$
with \pdf
\begin{equation}\label{logistaty}
    \widetilde{f}_Y(y)=\frac{\left(\frac{1}{D}\right)^\frac{1-D}{D}}{\Gamma\left(\frac{1-D}{D}\right)}\,y^{-\frac{1-D}{D}-1}e^{-\frac{1}{Dy}}
    =D\,\frac{\left(\frac{1}{Dy}\right)^{\frac{1}{D}}}{\Gamma\left(\frac{1-D}{D}\right)}\,e^{-\frac{1}{Dy}}\qquad\quad
    1>D
\end{equation}
then by taking
\begin{equation*}
    f_Y(y,t)=\sqrt{\widetilde{f}_Y(y)}\,g_Y(y,t)
\end{equation*}
we get the equation
\begin{equation}\label{sly}
    \partial_tg_Y(y,t)=\Lop{g_Y}(y,t)
\end{equation}
where $\Lop{\,\cdot\,}$ is now a Sturm-Liouville
operator\refeq{slop} with
\begin{eqnarray*}
    p(y)&=&B(y)=Dy^2\\
    q(y)&=&\frac{\left[B'(y)-\fvel(y)\right]^2}{4B(y)}-\frac{\left[B'(y)-\fvel(y)\right]'}{2}=\frac{(y-1)^2-2Dy^2}{4Dy^2}
\end{eqnarray*}
so that we finally get
\begin{equation*}
    \partial_tg_Y(y,t)=Dy^2\,\partial_y^2g_Y(y,t)+2Dy\,\partial_yg_Y(y,t)+\frac{2Dy^2-(y-1)^2}{4Dy^2}\,g_Y(y,t)
\end{equation*}
We next separate the variables by taking
\begin{equation*}
    g_Y(y,t)=e^{-\lambda t}G_Y(y)
\end{equation*}
to have the eigenvalue equation
\begin{equation*}
    \Lop{G_Y}(y)+\lambda G_Y(y)=0
\end{equation*}
which can be explicitly written as
\begin{equation}\label{eigeneqy}
  y^2\,G_Y''(y)+2y\,G_Y'(y)+\left[\frac{2Dy^2-(y-1)^2}{4D^2y^2}+\frac{\lambda}{D}\right]G_Y(y)=0
\end{equation}
or equivalently as
\begin{equation*}
    G_Y''(y)+\frac{2}{y}\,G_Y'(y)+\left[\frac{2Dy^2-(y-1)^2}{4D^2y^4}+\frac{\lambda}{Dy^2}\right]G_Y(y)=0
\end{equation*}
Now it is apparent that this equation has a Fuchsian singularity at
$y=+\infty$ because its coefficients asymptotically vanish quickly
enough, but also has non Fuchsian singularity in $y=0$ where its
second coefficient displays a 4$^{th}$ order pole. As a consequence
there is no standard procedure to solve it

%\subsubsection{Laplace transforms}

The two Fokker-Planck equations\refeq{logistfpx}
and\refeq{logistfpy} could finally also be recast in the form of
\ode's by means of a Laplace transform in $t$
\begin{align*}
   &\qquad\qquad\phi_X(x)=\int_0^{+\infty}\!\!e^{-pt}f_X(x,t)\,dt\qquad\qquad\phi_Y(y)=\int_0^{+\infty}\!\!e^{-pt}f_Y(y,t)\,dt
   \\
   &\int_0^{+\infty}\!\!e^{-pt}\partial_tf_X(x,t)\,dt=p\,\phi_X(x)-f_0(x)\qquad\int_0^{+\infty}\!\!e^{-pt}\partial_tf_Y(y,t)\,dt=p\,\phi_Y(y)-g_0(y)
\end{align*}
where, to keep the notation to a minimum, we hid the explicit
dependence on $p$: on the other hand the variable $p$ becomes an
external parameter in the subsequent \ode's which respectively are
\begin{align}
 & Dx^2\phi_X''(x)+[x^2+(4D-1)x]\phi_X'(x)+(2x+2D-1-p)\phi_X(x) = -f_0(x)\label{ode1} \\
 & Dy^2\phi_Y''(y)+\left[(2D+1)y-1\right]\phi_Y'(y)+(1-p)\phi_Y(y) =
 -g_0(y)\label{ode2}
\end{align}
where $f_0(x)$ and $g_0(y)$ are the initial \pdf's. The associated,
homogeneous equations
\begin{align}
 & Dx^2\phi_X''(x)+[x^2+(4D-1)x]\phi_X'(x)+(2x+2D-1-p)\phi_X(x) = 0\label{hode1} \\
 & Dy^2\phi_Y''(y)+\left[(2D+1)y-1\right]\phi_Y'(y)+(1-p)\phi_Y(y) =
 0\label{hode2}
\end{align}
are exactly solved by \textsc{Mathematica} in terms of confluent
hypergeometric functions. In particular, taking
\begin{equation*}
    \beta=\frac{1-D}{D}\qquad\qquad\eta(p)=\sqrt{\beta^2+\frac{p}{D}}
\end{equation*}
the general solution of\refeq{hode1} for instance is
\begin{eqnarray*}
  \phi_X(x) &=& e^{-\frac{x}{D}}x^{\beta-1+\eta(p)}\left[C_1\Psi\left(-\beta+\eta(p),1+2\eta(p),\frac{x}{D}\right)\right. \\
   &&\qquad\quad+\left.\frac{C_2}{2\eta(p)B\left(1+\beta-\eta(p),2\eta(p)\right)}\,\Phi\left(-\beta+\eta(p),1+2\eta(p),\frac{x}{D}\right)\right]
\end{eqnarray*}
where $\Phi$ and $\Psi$ are confluent hypergeometric functions
(see\mycite{grad} 9.201.1, 9.210.2). However, even if we can manage
to find the complete solution of the non homogeneous equation, the
problem of inverting such Laplace transforms would still stay with
us making these results rather incomplete, at least for the time
being

\subsection{$\theta$-logistic \sde}

We can now generalize the previous results to the $\theta$-logistic
\sde\refeq{thetalogdl} and we start by looking for a stationary
Boltzmann distribution. The transformation\refeq{b1} for the
coefficients\refeq{thetalogcoeffdl}, namely
\begin{equation*}
   y=h(x)=\ln x \qquad x=g(y)=e^y \qquad Y(t)=\ln X(t)\qquad X(t)=e^{Y(t)}
\end{equation*}
applied to the $\theta$-logistic \sde\refeq{thetalogdl} leads to
$\hb (y,t)=1$, and from\refeq{a1} to the drift coefficient
\begin{equation*}
      \ha (y)=1-D-e^{\theta y}
\end{equation*}
namely to the Smoluchowsky \sde
\begin{equation*}
    dY(t)=\left(1-D-e^{\theta Y(t)}\right)\,dt+dW(t)
\end{equation*}
and from\refeq{potential} to a dimensionless potential
\begin{equation*}
    \chi(y)=\frac{\phi(y)}{kT}=\frac{e^{\theta y}}{\theta D}-\frac{1-D}{D}\,y+c
\end{equation*}
that -- provided now that $1>D$ -- gives rise to the following
stationary generalized log-gamma Boltzmann distribution
(see\mycite{grad} 3.328 for the normalization integral)
\begin{equation*}
    \frac{e^{-\frac{e^{\theta y}}{\theta D}+\frac{1-D}{D}\,y}}{\frac{1}{\theta}(\theta D)^{\frac{1-D}{\theta D}}\,\Gamma\left(\frac{1-D}{\theta D}\right)}\qquad\quad1>D
\end{equation*}
This can finally be transformed back to the original process
$X(t)=e^{Y(t)}$ giving as stationary density
\begin{equation*}
    \frac{\theta\, x^{\frac{1-D}{D}-1}e^{-\frac{x^\theta}{\theta D}}}{(\theta D)^{\frac{1-D}{\theta D}}\Gamma\left(\frac{1-D}{\theta D}\right)}
    \qquad\quad
    1>D
\end{equation*}
which is the \pdf\ of the generalized gamma law
$\gam_\theta\left(\frac{1-D}{D},\frac{1}{(\theta
D)^{1/\theta}}\right)$

We then go on to linearize the \sde\refeq{thetalogdl}: from the
coefficients\refeq{thetalogcoeffdl} we find
\begin{equation*}
    q(x)=1-D-x^\theta\qquad\quad
    b(x)q'(x)=-\theta x^\theta\qquad\quad
    \frac{1}{q'(x)}\frac{d}{dx}\left[b(x)q'(x)\right]=\theta
\end{equation*}
so that the compatibility condition\refeq{lincond} is satisfied and
from\refeq{tratolin1},\refeq{lincond2} we have
\begin{equation*}
    \hb_1=-\theta\qquad p(x)=\int\frac{dx}{b(x)}=\ln x\qquad\quad h(x)=c\,e^{\hb_1 p(x)}=c\,e^{-\theta\ln x}=\frac{c}{x^\theta}
\end{equation*}
If we then choose $c=1$ as integration constant the reciprocal
transformation relations are
\begin{equation*}
    y=h(x)=\frac{1}{x^\theta}\qquad\quad
    x=g(y)=\frac{1}{y^{1/\theta}}\qquad\quad
    Y(t)=\frac{1}{X(t)^\theta}\qquad\quad X(t)=\frac{1}{Y(t)^{1/\theta}}
\end{equation*}
so that with
\begin{equation*}
    h'(x)=-\frac{\theta}{x^{1+\theta}}\qquad\qquad h''(x)=\frac{\theta(1+\theta)}{x^{2+\theta}}
\end{equation*}
from\refeq{itoy1} and\refeq{itoy2} we have
\begin{eqnarray*}
  \ha(y) &=& h'(g(y))a(g(y))+D\,h''(g(y))b^2(g(y))=\big((1+\theta)D-1\big)\,\theta y+\theta \\
  \hb(y) &=& h'(g(y))b(g(y))=-\theta y
\end{eqnarray*}
namely from\refeq{linearizedcoeff}
\begin{equation*}
    \ha_0=\theta\qquad\quad\ha_1=\big((1+\theta)D-1\big)\,\theta\qquad\qquad\hb_0=0\qquad\quad\hb_1=-\theta
\end{equation*}
and hence the new \sde\ is
\begin{equation}\label{thetalogdlt}
    dY(t)=  \big[\big((1+\theta)D-1\big)Y(t)+1\big]\,\theta dt-\theta Y(t)\,dW(t)
\end{equation}
Taking now as in\refeq{linearhsolbar1}
\begin{equation*}
    \overline{Z}(t)=\theta(D-1)t-\theta W(t)\sim\norm\big((D-1)\theta t\,,\,2D\theta^2t\big)
\end{equation*}
the general solution\refeq{linearsolnot} of the linearized
\sde\refeq{thetalogdlt} for $Y(0)=Y_0$ is
\begin{equation*}
    Y(t)=e^{\overline{Z}(t)}\left(Y_0+\theta\int_0^te^{-\overline{Z}(u)}du\right)
\end{equation*}
while the solution $X(t)$ of the $\theta$-logistic
\sde\refeq{thetalogdl} for $Y_0=X_0^{-\theta}$ is
\begin{equation}\label{thetalogsol}
    X(t)=\left(\frac{X_0^\theta\, e^{-\overline{Z}(t)}}{1+\theta
    X_0^\theta\int_0^te^{-\overline{Z}(u)}du}\right)^{1/\theta}
\end{equation}
Remark that here too, with a degenerate initial condition
$X_0=x_0,\:\Pqo$ and by switching off the Wiener noise ($D=0$) we
get $\overline{Z}(t)=-\theta t$, and the solution\refeq{thetalogsol}
exactly coincides with the solution\refeq{thetalogdet} of the
deterministic $\theta$-logistic \emph{ODE} discussed in the
Appendix\myref{detlogist}. Then again, taking $X_0=y,\;\Pqo$ at a
time $0\leq s\leq t$ we have the solution
\begin{equation}\label{thetalogtranssol}
    X(t)=\left(\frac{y^\theta\, e^{-\overline{Z}(t-s)}}{1+\theta
    y^\theta\int_s^te^{-\overline{Z}(u)}du}\right)^{1/\theta}
\end{equation}
whose \pdf\ $f(x,t|y,s)$ will be the \textbf{\emph{transition}
\pdf}\ of our $\theta$-logistic process. If moreover we define the
derivable process
\begin{equation}\label{thetaAproc}
    A(t)=X_0^\theta\int_0^te^{-\overline{Z}(u)}du\qquad\qquad\dot A(t)=X_0^\theta e^{-\overline{Z}(t)}
\end{equation}
the solution\refeq{thetalogsol} takes the equivalent forms (see
also\mycite{grad} $3.434.2$)
\begin{eqnarray}\label{thetalogsol2}
    X(t)&=&\!\!\left(\frac{\dot A(t)}{1+\theta
    A(t)}\right)^{\frac{1}{\theta}}=\left({\frac{1}{\theta}\,\frac{d}{dt}\ln\big[1+\theta A(t)\big]}\right)^{\frac{1}{\theta}}
    =\left(\frac{d}{dt}\int_0^\infty e^{-u}\frac{1-e^{-\theta uA(t)}}{\theta u}\,du\right)^{\frac{1}{\theta}}\nonumber\\
    &=&\!\!\left(-\int_0^\infty\frac{e^{-u}}{\theta u}\,\frac{d}{dt}\,e^{-\theta uA(t)}du\right)^{\frac{1}{\theta}}
    =\left(\dot A(t)\int_0^\infty e^{-u(1+\theta A(t))}du\right)^{\frac{1}{\theta}}
\end{eqnarray}

Retracing finally the procedure of  the Appendix\myref{transpdf}
leading to the semi-explicit transition \pdf\refeq{logistpdf}, from
the coefficients\refeq{thetalogcoeffdl} we begin by defining the
functions
\begin{equation*}
   y=h(x)=\ln x\qquad\quad
   x=g(y)=e^y\qquad\quad
   \ha(y)=1-D-e^{\theta y}
\end{equation*}
so that we also find
\begin{equation*}
  \beta(y) = -\frac{1}{4D}\left(1-e^{\theta y}\right)^2+\frac{1+(\theta-1)e^{\theta y}}{2}-\frac{D}{4}
\end{equation*}
Keeping then for $\overline{h}(r)$ and $\overline{W}_{st}(r)$ the
same definitions of the Section\myref{semiexpl}, we have now
\begin{align*}
 & \beta\left(\overline{W}_{st}(r)+\overline{h}(r)\right) \\
 & \qquad\qquad = \frac{2-D+2(\theta-1)x^{r\theta} y^{(1-r)\theta} e^{\theta\overline{W}_{st}(r)}}{4}-\frac{1}{4D}\left(1-x^{r\theta} y^{(1-r)\theta} e^{\theta\overline{W}_{st}(r)}\right)^2
\end{align*}
and hence we get
\begin{align*}
 & Z(s,t) \\
 & = -\frac{(1-D)^2}{4D}+\frac{1+(\theta-1)D}{2D}\,y^\theta\int_0^1\left(\frac{x}{y}\right)^{r\theta} e^{\theta\overline{W}_{st}(r)}dr
   -\frac{y^{2\theta}}{4D}\int_0^1\left(\frac{x}{y}\right)^{2r\theta} e^{2\theta\overline{W}_{st}(r)}dr
\end{align*}
We thus find for the expectation factor in the transition
\pdf\refeq{pdf}
\begin{align*}
 &\qquad\qquad\EXP{e^{(t-s)Z(s,t)}} = e^{-\frac{(1-D)^2}{4D}(t-s)}\mu_\theta(x,t;y,s)\\
 & \mu_\theta(x,t;y,s)=\bm{E}\Bigg[\exp\left\{-\frac{t-s}{4D}\Bigg(
            y^{2\theta}\int_0^1\!\!\Bigg(\frac{x}{y}\right)^{2r\theta}\!\!e^{2\theta\overline{W}_{st}(r)}dr\\
 &\qquad\qquad\qquad\qquad\qquad\qquad-2\big(1+(\theta-1)D\big)y^\theta\int_0^1\!\!\left(\frac{x}{y}\right)^{r\theta}\!\! e^{\theta\overline{W}_{st}(r)}dr\Bigg)\Bigg\}\Bigg]
\end{align*}
As for the other factors in\refeq{pdf} we now have
\begin{eqnarray*}
  \frac{1}{2D}\int_y^x\frac{a(z)}{b^2(z)}\,dz &=& \frac{1}{2D}\left(\ln\frac{x}{y}-\frac{x^\theta-y^\theta}{\theta}\right) \\
  -\frac{1}{4D(t-s)}\left(\int_y^x\frac{dz}{b(z)}\right)^2&=&-\frac{1}{4D(t-s)}\ln^2\frac{x}{y}
\end{eqnarray*}
and therefore with a little algebra we find out for the transition
\pdf
\begin{equation}
  f(x,t|y,s) =\frac{e^{-\frac{x^\theta-y^\theta}{2D\theta}-\frac{1}{4D(t-s)}\left((1-D)(t-s)-\ln\frac{x}{y}\right)^2}}{x\sqrt{4\pi
   D(t-s)}}\label{thetalogpdf}\,\mu_\theta(x,t;y,s)
\end{equation}

\subsection{Integrals of a geometric Wiener
process}\label{intgeowien}

In the Section\myref{linsde} we have discussed the solutions of the
\sde\ ruled by a logistic dynamics and we have found them contingent
on processes basically of the type
\begin{equation}\label{gwp}
    X(t)=\int_0^te^{W(s)}ds
\end{equation}
where $W(t)\sim\norm(0,2Dt)$ is a Wiener process with diffusion
coefficient $2D$. The processes\refeq{gwp} are also known as
\emph{Exponential Functionals of Brownian Motion} and have been
extensively studied in the financial context in a rather
mathematical setting\mycite{yor}. By postponing a more accurate
analysis to a forthcoming paper, in the present section we will
instead scrutinize the distribution of $X(t)$ with more elementary
tools leading of course just to partial results.

\subsubsection{Moments}

We will begin by looking at the moments of\refeq{gwp}
\begin{eqnarray}\label{mom}
    M_n(t)=\EXP{X^n(t)}&=&\EXP{\int_0^tds_1\ldots\int_0^tds_ne^{W(s_1)+\ldots+W(s_n)}}\nonumber \\
    &=&\int_0^tds_1\ldots\int_0^tds_n\EXP{e^{\sum_{k=1}^nW(s_k)}}\label{mom}
\end{eqnarray}
Since it is
\begin{equation*}
    \cov{W(s)}{W(t)}=2D\,(s\wedge t)
\end{equation*}
we have for every choice of $t_1,\ldots,t_n$ that
\begin{equation*}
    \big(W(t_1),\ldots,W(t_n)\big)\sim\norm\left(\bm 0,2D\,\mathbb{A}\right)
\end{equation*}
where $\bm 0=(0,\ldots,0)$ and
\begin{equation*}
    \mathbb{A}=\left(
                 \begin{array}{ccccc}
                   t_1 & t_1\wedge t_2 & t_1\wedge t_3 & \ldots & t_1\wedge t_n \\
                   t_2\wedge t_1 & t_2 & t_2\wedge t_3 &  & t_2\wedge t_n \\
                   t_3\wedge t_1 & t_3\wedge t_2 & t_3 &  & t_3\wedge t_n \\
                   \vdots &  &  & \ddots &  \vdots \\
                   t_n\wedge t_1 & t_n\wedge t_2 & t_n\wedge t_3 & \ldots & t_n \\
                 \end{array}
                \right)
\end{equation*}
As a consequence
\begin{equation*}
    \sum_{k=1}^nW(s_k)=(1,\ldots,1)\left(\begin{array}{c}
    W(s_1) \\ \vdots \\ W(s_n)
    \end{array}\right)\sim\norm\left(0\,,\,2D\sum_{j,k=1}^n(s_j\wedge s_k)\right)
\end{equation*}
and $e^{\sum_{k=1}^nW(s_k)}$ is the corresponding log-normal \rv, so
that
\begin{equation}\label{symmexp}
   \EXP{e^{\sum_{k=1}^nW(s_k)}}=e^{D\sum_{j,k=1}^n(s_j\wedge s_k)}
\end{equation}
which is apparently invariant under every permutation of the
variables, while the moments are
\begin{equation*}
    M_n(t)=\int_0^tds_1\ldots\int_0^tds_n\,e^{D\sum_{j,k=1}^n(s_j\wedge s_k)}=\int_{[0,t]^n}ds_1\ldots ds_n\,e^{D\sum_{j,k=1}^n(s_j\wedge s_k)}
\end{equation*}
In this integral the variables $s_k$ are not ordered, but we can go
around this problem in the following way: consider the subset of
$[0,t]^n$
\begin{equation*}
    B=\left\{(s_1,\ldots,s_n)\in[0,t]^n\,:\,0\le s_1\le\ldots\le s_n\le t\right\}
\end{equation*}
where the variables are ordered according to their indices, and let
$\mathcal{P}$ be the family of the $n!$ permutations $\Pi$ of
$s_1,\ldots,s_n$. The $n!$ subsets $\Pi(B)$ obtained by permutations
of the variables in $B$ are then such that
\begin{equation*}
    \bigcup_{\Pi\in\mathcal{P}}\Pi(B)=[0,t]^n
\end{equation*}
while on the other hand -- because of the symmetry of\refeq{symmexp}
under permutations -- all the integrals of\refeq{symmexp} on every
$\Pi(B)$ take the same value. We have then that
\begin{equation*}
    M_n(t)=n!\int_B ds_1\ldots ds_n\,e^{D\sum_{j,k=1}^n(s_j\wedge s_k)}
    =n!\int_0^t\!ds_n\ldots\!\int_0^{s_3}\!\!ds_2\int_0^{s_2}\!\!ds_1\,e^{D\sum_{j,k=1}^n(s_j\wedge s_k)}
\end{equation*}
and since for $s_1\le s_2\le\ldots\le s_n$ it is
\begin{eqnarray*}
  \sum_{j,k=1}^n(s_j\wedge s_k) &=& \sum_{k=1}^ns_k+2\sum_{j<k}(s_j\wedge s_k)=\sum_{k=1}^ns_k+2\sum_{j=1}^{n-1}\sum_{k=j+1}^n(s_j\wedge s_k)\\
   &=&\sum_{k=1}^ns_k+2\sum_{j=1}^{n-1}\sum_{k=j+1}^n
   s_j=\sum_{k=1}^ns_k+2\sum_{j=1}^{n-1}(n-j)s_j \\
   &=&\sum_{k=1}^ns_k+2\sum_{j=1}^n(n-j)s_j=\sum_{k=1}^n[2(n-k)+1]s_k
\end{eqnarray*}
we finally get
\begin{equation*}
    M_n(t)
    =n!\int_0^t\!ds_n\ldots\!\int_0^{s_3}\!\!ds_2\int_0^{s_2}\!\!ds_1\,e^{D\sum_{k=1}^n[2(n-k)+1]s_k}
\end{equation*}
and with the changes of variables $v_k=Ds_k$
\begin{eqnarray}\label{moment}
    M_n(t)&=&\frac{n!}{D^n}\int_0^t\!dv_n\ldots\int_0^{v_3}\!\!dv_2\int_0^{v_2}\!\!dv_1\,e^{\sum_{k=1}^n[2(n-k)+1]v_k}\nonumber\\
    &=&\frac{n!}{D^n}\int_0^t\!dv_ne^{v_n}\int_0^{v_n}\!\!dv_{n-1}e^{3v_{n-1}}\ldots\int_0^{v_3}\!\!dv_2\,e^{(2n-3)v_2}\int_0^{v_2}\!\!dv_1\,e^{(2n-1)v_1}
\end{eqnarray}
By looking now at the explicit calculations for the first few values
of $n$
\begin{eqnarray*}
  M_1(t) &=& \frac{e^{Dt}-1}{D} \\
  M_2(t) &=& \frac{e^{4Dt}-4e^{Dt}+3}{6D^2} \\
  M_3(t) &=& \frac{e^{9Dt}-6e^{4Dt}+15e^{Dt}-10}{60D^3} \\
  M_4(t) &=&
  \frac{e^{16Dt}-8e^{9Dt}+28e^{4Dt}-56e^{Dt}+35}{840D^4}\;\ldots
\end{eqnarray*}
we can conjecture the following general form for the moments
of\refeq{gwp}
\begin{equation}\label{mom1}
    M_n(t)=(-1)^n\frac{n!}{D^n}\sum_{k=0}^n(-1)^k\frac{2-\delta_{k0}}{(n-k)!\,(n+k)!}\,e^{k^2Dt}
\end{equation}
but we have yet no proof by recurrence and induction, so
that\refeq{moment} remains for the time being our last validated
result

\subsubsection{Characteristic function}\label{chfexpansion}

Leaving aside for now every convergence question\footnote{This is
not at all a small detail, as it will be clear at the end of the
present section. On the other hand it is known that this problem
already exists for the lognormal distributions: while all the
moments exist and are finite the generating function does not exist,
and the characteristic function can not be represented as a
convergent series\mycite{lnorm}. This is related indeed to the fact
that the lognormal distribution is not uniquely determined by its
moments, and it would not be surprising then to find that this
behavior extends also to the integrals of lognormal processes},
starting from\refeq{mom1} we could surmise that the characteristic
function of $X(t)$ takes the form
\begin{eqnarray*}
  \varphi(u,t) &=& \sum_{n=0}^\infty\frac{(iu)^n}{n!}M_n(t)
      =\sum_{n=0}^\infty\left(\frac{-iu}{D}\right)^n\sum_{k=0}^n(-1)^k\frac{2-\delta_{k0}}{(n-k)!\,(n+k)!}\,e^{k^2Dt}\\
   &=&\sum_{k=0}^\infty(-1)^k(2-\delta_{k0})e^{k^2Dt}\sum_{n=k}^\infty\left(\frac{-iu}{D}\right)^n\frac{1}{(n-k)!\,(n+k)!}
\end{eqnarray*}
and by changing the index $n$ into $\ell=n-k$
\begin{equation*}
    \varphi(u,t) =\sum_{k=0}^\infty\left(\frac{iu}{D}\right)^k(2-\delta_{k0})e^{k^2Dt}\sum_{\ell=0}^\infty\left(\frac{-iu}{D}\right)^\ell\frac{1}{\ell!\,(\ell+2k)!}
\end{equation*}
It is known on the other hand that (see\mycite{grad} 8.402)
\begin{equation*}
    \sum_{\ell=0}^\infty\frac{z^\ell}{\ell!\,(\ell+2k)!}=(-z)^{-k}J_{2k}(2i\sqrt{z})
\end{equation*}
where $J_n(x)$ are the Bessel functions of the first kind, and hence
\begin{equation*}
     \varphi(u,t)
     =\sum_{k=0}^\infty(2-\delta_{k0})e^{k^2Dt}J_{2k}\left(2i\sqrt{\frac{-iu}{D}}\right)
\end{equation*}
Since in general $\varphi(-u)=\overline{\varphi(u)}$, we can
restrict ourselves to $u>0$ and in this case we have
\begin{equation*}
    2i\sqrt{\frac{-iu}{D}}=\pm2e^{i\frac{\pi}{4}}\sqrt{\frac{u}{D}}=\pm(1+i)\sqrt{\frac{2u}{D}}\qquad\qquad
    u>0
\end{equation*}
On the other hand we also know that $J_{2k}(-z)=J_{2k}(z)$ so that
we finally have
\begin{eqnarray}
     \varphi(u,t)
     &=&\sum_{k=0}^\infty(2-\delta_{k0})e^{k^2Dt}J_{2k}\left((1+i)\sqrt{\frac{2u}{D}}\right)\nonumber
     \\
     &=&2\sum_{k=1}^\infty
     e^{k^2Dt}J_{2k}\left((1+i)\sqrt{\frac{2u}{D}}\right)+J_0\left((1+i)\sqrt{\frac{2u}{D}}\right)\qquad u>0\label{chf}
\end{eqnarray}
This result looks however only formal because the presence of terms
$e^{k^2Dt}$ in the sums lends no hope for a convergence whatsoever.
This was moreover a foregone conclusion since the (absolute)
moments\refeq{mom1} utterly fail the convergence test for the Taylor
expansion of the characteristic function: the moments\refeq{mom1}
coincide indeed with the absolute moments because out \rv's are
always positive (as for every exponential function), and hence the
convergence of the Taylor expansion would require
\begin{equation*}
    \overline{\lim_n}\,\frac{\sqrt[n]{M_n(t)}}{n}<+\infty
\end{equation*}
while a few numerical trials show that the limit diverges for every
choice of $t$ and $D$

\subsubsection{Finite sums of a geometric Wiener process}

Since the process $X(t)$ in\refeq{gwp} is the integral of a
geometric Wiener process, we could first of all investigate the laws
of sums of a geometric Wiener process at different times. Let us
begin with the simplest case
\begin{equation}\label{z}
    Z=e^{W(s)}+e^{W(t)}\qquad\quad s\leq t
\end{equation}
by remarking first that the two log-normal \rv's $e^{W(s)}$ and
$e^{W(t)}$ are not independent. The \rv\ $Z$ can however be put in
the form of the product of two independent \rv's
\begin{equation*}
    Z=e^{W(s)}\left(1+e^{W(t)-W(s)}\right)=XY
\end{equation*}
where we know that $X>0$ is log-normal, while $Y>1$ is a $1$-shifted
log-normal:
\begin{equation*}
    X=e^{W(s)}\sim\lnorm(0,2Ds)\qquad\quad Y-1=e^{W(t)-W(s)}\sim\lnorm(0,2D(t-s))
\end{equation*}
namely the \pdf's respectively are
\begin{equation}\label{xypdf}
    f_X(x)=\frac{e^{-\frac{\ln^2 x}{4Ds}}}{x\sqrt{4\pi Ds}}\,\vartheta(x)\qquad\quad f_Y(y)=\frac{e^{-\frac{\ln^2(y-1)}{4D(t-s)}}}{(y-1)\sqrt{4\pi D(t-s)}}\,\vartheta(y-1)
\end{equation}
$\vartheta(x)$ being the Heaviside function. To find the \pdf\
$f_Z(z)$ of $Z$ we could then remark that $\ln Z=\ln X + \ln Y$ is
the sum of two independent \rv's where in particular $\ln
X=W(s)\sim\norm(0,2Ds)$. We could hence first calculate the \pdf\ of
$\ln Z$ as the convolution of the \pdf's of $W(s)$ and $\ln Y$, and
then transform it back to the \pdf\ of $Z$. It is important to
remark however that $\ln Y$, as the logarithm of a $1$-shifted
log-normal, by no means is a normal \rv\ as can be apparently argued
from the simple remark that, being $Y>1$, we always get $\ln Y>0$.
The \pdf\ of $\ln Y$ can of course be explicitly calculated with the
usual procedure, but it turns out to have an involuted form which
makes hard to calculate the required convolution, and still harder
to find back $f_Z(z)$. Alternatively we could try to directly
calculate the \pdf\ of the product of two non-negative, independent
\rv's according to the following result
\begin{prop}
If $Z=XY$ is the product of two \ac\ \rv's with joint \pdf\
$f(x,y)$, then its \pdf\ is
\begin{equation}\label{zpdf}
    f_Z(z)=\int_0^{\infty}\frac{dx}{x}\left[f\left(x,\frac{z}{x}\right)+f\left(-x,-\frac{z}{x}\right)\right]
\end{equation}
and when in particular we take $X\ge0,\;Y\ge0$ it becomes
\begin{equation}\label{zpdf0}
    f_Z(z)=\vartheta(z)\int_0^{\infty}\frac{dx}{x}f\left(x,\frac{z}{x}\right)
\end{equation}
Finally, if we also suppose that $X,Y$ are independent with
marginals $f_X(x)$ and $f_Y(y)$ the \pdf\ is
\begin{equation}\label{zpdf0ind}
    f_Z(z)=\vartheta(z)\int_0^{\infty}\frac{dx}{x}f_X\left(x\right)f_Y\left(\frac{z}{x}\right)
\end{equation}
\end{prop}
\begin{figure}
\begin{center}
\includegraphics*[width=15cm]{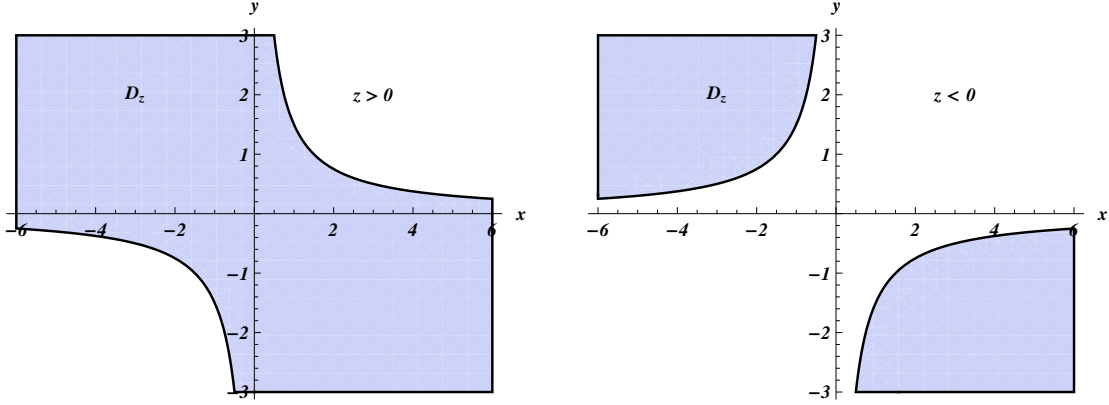}
\caption{The shadowed regions show the integration domain
$D_z$}\label{domain}
\end{center}
\end{figure}
 \indim
Starting from the \cdf\ of $Z$ we have
\begin{equation*}
    F_Z(z)=\PR{Z\le z}=\PR{XY\le z}=\int\!\!\!\int_{D_z}f(x,y)\,dx\,dy
\end{equation*}
where $D_z=\{(x,y)\in\RE^2\,:\,xy\le z\}$. Looking at the
Figure\myref{domain} we see then that
\begin{eqnarray*}
  F_Z(z) &=& \int_{-\infty}^0dx\int_{z/x}^\infty dy\,f(x,y)+\int_0^{\infty}dx\int_{-\infty}^{z/x} dy\,f(x,y) \\
   &=&  \int_0^{\infty}dx\int_{-z/x}^\infty dy\,f(-x,y)+\int_0^{\infty}dx\int_{-\infty}^{z/x}
   dy\,f(x,y) \\
   &=&  \int_0^{\infty}dx\int_{-\infty}^{z/x} dy\,f(-x,-y)+\int_0^{\infty}dx\int_{-\infty}^{z/x}
   dy\,f(x,y)
\end{eqnarray*}
and by introducing a new variable $u=xy$
\begin{eqnarray*}
  F_Z(z) &=& \int_0^{\infty}dx\int_{-\infty}^z
  \frac{du}{x}\,f\left(-x,-\frac{u}{x}\right)+\int_0^{\infty}dx\int_{-\infty}^z
   \frac{du}{x}\,f\left(x,\frac{u}{x}\right) \\
    &=&\int_{-\infty}^z du\int_0^{\infty}\frac{dx}{x}\left[f\left(x,\frac{u}{x}\right)+f\left(-x,-\frac{u}{x}\right)\right]
\end{eqnarray*}
so that at once we get the \pdf\refeq{zpdf}, while the other two
formulas\refeq{zpdf0} and\refeq{zpdf0ind} immediately follow
 \findim

\begin{figure}
\begin{center}
\includegraphics*[width=15cm]{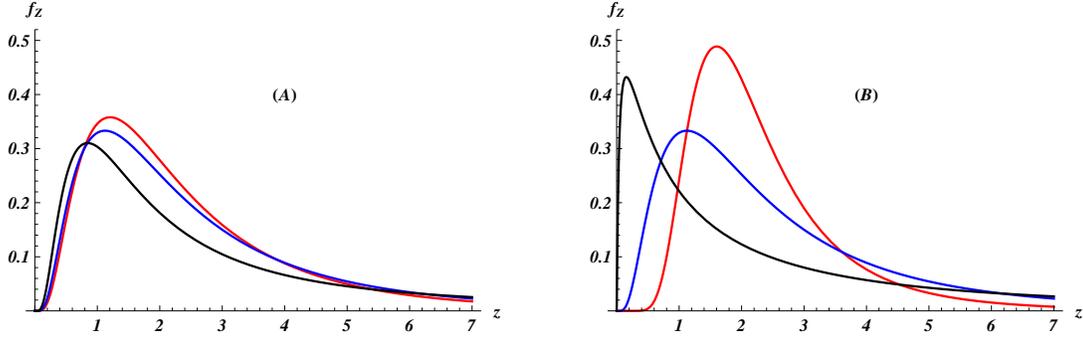}
\caption{Numerical instances of the \pdf\ $f_Z(z)$ in\refeq{zpdf}:
$\bm{(\!A)}$ $4Ds=1.0$, while $4D(t-s)=0.1$ (\emph{red}) $1.0$
(\emph{blue}) $8.0$ (\emph{black}); $\bm{(\!B)}$ $4D(t-s)=1$, while
$4Ds=0.2$ (\emph{red}) $1.0$ (\emph{blue}) $5.0$
(\emph{black})}\label{fz}
\end{center}
\end{figure}

\noindent From\refeq{zpdf0ind} and\refeq{xypdf} we then have the
following \pdf\ for $Z$
\begin{equation}\label{zpdf}
    f_Z(z)=\vartheta(z)\int_0^z\frac{e^{-\frac{\ln^2 x}{4Ds}}}{x\sqrt{4\pi
    Ds}}\,\frac{e^{-\frac{\left[\ln(z-x)-\ln x\right]^2}{4D(t-s)}}}{(z-x)\sqrt{4\pi D(t-s)}}\,dx
\end{equation}
which again is not an easy calculation to perform, even if it looks
tantalizingly near to an explicit answer. Numerical integration
shows that\refeq{zpdf} is correctly normalized, and numerical plots
in Figure\myref{fz}
 display a very reasonable behavior confirming
that our calculation is so far acceptable: that notwithstanding, the
unavailability of a complete result for such a simple case as the
\rv\ $Z$ in\refeq{z} also uphold the view that finding the law of
$X(t)$ in\refeq{gwp} is a problem hard to crack

\subsubsection{Differential equations}

When an \ac\ process $X(t)$ is solution of a \sde\
\begin{equation}\label{xsde}
    dX(t)=a(X(t),t)\,dt+b(X(t),t)\,dW(t)
\end{equation}
then $X(t)$ is Markovian, with almost every trajectory everywhere
continuous, and its \pdf\ $f_X(x,t)$ is solution of a Fokker-Planck
equation
\begin{equation}\label{xfpe}
    \partial_tf_X(x,t)=-\partial_x\left[A(x,t)f_X(x,t)\right]+\frac{1}{2}\,\partial_x^2\left[B(x,t)f_X(x,t)\right]
\end{equation}
where
\begin{equation*}
    a(x,t)=A(x,t)\qquad\qquad B(x,t)=2Db^2(x,t)
\end{equation*}
In this case it also satisfies a \emph{continuity equation} which
represents a requirement of probability conservation: we can indeed
immediately recast\refeq{xfpe} into the form
\begin{equation}\label{xconteq}
    \partial_tf_X(x,t)+\partial_x\left[v(x,t)f_X(x,t)\right]=0
\end{equation}
provided that
\begin{equation}\label{xvel}
    v(x,t)=A(x,t)-\frac{\partial_x\left[B(x,t)f_X(x,t)\right]}{2f_X(x,t)}
\end{equation}
It is apparent however that in the present context the continuity
equation\refeq{xconteq} is not a new equation really different from
the Fokker-Planck equation\refeq{xfpe}, and this is made clear in
particular by the fact that the velocity field\refeq{xvel} is
contingent on the solution $f_X(x,t)$ of\refeq{xfpe}: in other words
here $v(x,t)$ does not represent an external, given field but
depends on the solution $f_X(x,t)$ so that (if $A(x,t)$ and $B(x,t)$
are given) we can directly calculate $v$ from $f_X$, and conversely
$f_X$ from $v$

Not every stochastic process, however, is Markovian, and in
particular $X(t)$ defined in\refeq{gwp} is not. In this case neither
the process trajectories will satisfy a \sde\ of the
type\refeq{xsde}, nor its \pdf\ will be a solution of a \pde\ of the
type\refeq{xfpe}. This is quite understandable, and in fact every
quest for some other kind (for instance) of \pde, even if possible,
is doomed to futility: since $X(t)$ is not Markovian, in order to
find the law of the process it would be far from enough to know its
one time \pdf\ $f_X(x,t)$ together with its transition \pdf\
$f_X(x,t|y,s)$. We would need instead the knowledge of every joint
\pdf\ $f_X(x_1,t_;\ldots;x_n,t_n)$ that in any case could not be
extracted from a single \pde.

On the other hand the conservation of the probability should be
guaranteed in any case, and hence the \pdf\ $f_X(x,t)$ is supposed
to satisfy some kind of continuity equation of the
type\refeq{xconteq}, but for the fact that now this continuity
equation can no longer be derived from a corresponding \fpe. We must
at once remark, however, that\refeq{xconteq} in no way can surrogate
the role of a \fpe: first of all its possible solutions will not
constitute the basis to build the process laws; and furthermore --
as we already have remarked -- $v(x,t)$ is not a given function
independent from the solution. In general the connection between
$v(x,t)$ and $f_X(x,t)$ will not be as simple as\refeq{xvel}, but in
any case every possible solution will be associated to its own
velocity field.

As a matter of fact, however, when the process is not Markovian we
should give up our old habit of thinking to the \emph{different
processes} selected by different initial conditions for a single
transition \pdf\ as to a \emph{unified process}: now every global
law (represented by the said  joint \pdf's
$f_X(x_1,t_;\ldots;x_n,t_n)$) defines a different process and we do
not see in general a way to detect homogeneous classes among them.
This means, among others, that for every continuity equation with a
given $v(x,t)$ there will be just one possible solution of interest
for us, and that a family of processes could be located only through
their mating with the velocity fields

Two remarks are in order here: first, we could revert to our
initial, more narrow aim of finding just the law of the \rv's $X(t)$
in\refeq{gwp} and not the global law of the geometric Wiener process
that these \rv's represent for $t>0$. We have found however that
even restricting the scope of our enquiry to this carefully
circumscribed problem will not make easy to pick up a meaningful
solution. Second, we could try to circumvent the non Markovianity of
$X(t)$ in a way reminiscent of an Ornstein-Uhlenbeck procedure: the
process $X(t)$ is apparently derivable with
$\dot{X}(t)=e^{W(t)}=Y(t)$ so that its stochastic differential is
\begin{equation*}
    dX(t)=e^{W(t)}dt=Y(t)\,dt
\end{equation*}
This however does not constitute a \sde, and hence the \pdf\ of
$X(t)$ (which arguably is non Markovian, as in the
Ornstein-Uhlenbeck case) does not satisfy a Fokker-Planck equation.
The process $Y(t)=e^{W(t)}$ on the other hand -- from the \ito\
formula -- is a solution of the \sde\ for a geometric Wiener process
\begin{equation*}
    dY(t)=DY(t)\,dt+Y(t)\,dW(t)
\end{equation*}
and hence it is Markovian so that its \pdf\ $f_Y(y,t)$ obeys the
corresponding Fokker-Planck equation
\begin{eqnarray}
  \partial_tf_Y(y,t) &=& -D\partial_x\left[yf_Y(y,t)\right]+D\partial^2_y\left[y^2f_Y(y,t)\right] \nonumber\\
   &=& Dy^2\partial^2_Yf_Y(y,t)+3D\partial_yf_Y(y,t)+Df_Y(y,t)
   \label{yfpe}
\end{eqnarray}
The pair $X(t),Y(t)$ will thus satisfy the system
\begin{equation*}
    \left\{
      \begin{array}{lcl}
        dX(t) &=& Y(t)\,dt \\
        dY(t) &=& DY(t)\,dt+Y(t)\,dW(t)
      \end{array}
    \right.
\end{equation*}
If then we define the vector process
\begin{equation*}
    \bm Z(t)=\left(
               \begin{array}{c}
                 X(t) \\
                 Y(t) \\
               \end{array}
             \right)
\end{equation*}
it will be a solution of the two-components, vector \sde
\begin{equation}\label{vsde}
    d\bm Z(t)=\bm a(\bm Z(t))\,dt+\mathbb{C}(\bm Z(t))\,d\bm W(t)
\end{equation}
where
\begin{equation*}
    \bm a(\bm z)=\bm a(x,y)=\left(
               \begin{array}{c}
                 y \\
                 Dy \\
               \end{array}
             \right)
    \qquad\quad
    \mathbb{C}(\bm z)=\mathbb{C}(x,y)=\left(
                                        \begin{array}{cc}
                                          0 & 0 \\
                                          0 & y \\
                                        \end{array}
                                      \right)
\end{equation*}
while the vector Wiener process can be taken as
\begin{equation*}
    \bm W(t)=\left(
               \begin{array}{c}
                 W_X(t) \\
                 W(t) \\
               \end{array}
             \right)
\end{equation*}
$W_X(t)$ being any auxiliary Wiener process apparently playing no
role in the discussion of the \sde\refeq{vsde}. The vector process
$\bm Z(t)$ is then Markovian and to its vector \sde\refeq{vsde} it
is possible to associate a (2+1)-dimensional Fokker-Planck equation
for the \pdf\ $f(\bm z,t)=f(x,y,t)$: solving this equation would
lead in principle to the complete law of the process $\bm Z(t)$ and
hence, by marginalization, to the much sought-after law of its
component $X(t)$
\begin{prop}
The Fokker-Planck equation of the process  $\bm Z(t)$ is
\begin{equation}\label{fpe}
    \partial_tf(x,y,t)=D\partial_y^2\left[y^2f(x,y,t)\right]-y\partial_xf(x,y,t)-D\partial_y\left[yf(x,y,t)\right]
\end{equation}
while the marginal \pdf\ $f_X(x,t)$ of the process $X(t)$ obeys the
continuity equation
\begin{equation}\label{conteq}
    \partial_tf_X(x,t)+\partial_x\left[v(x,t)f_X(x,t)\right]=0
\end{equation}
where the velocity field is
\begin{equation}\label{vfield}
    v(x,t)=\EXP{\left.\dot{X}(t)\,\right|\,X(t)=x}=\EXP{e^{W(t)}\left|\int_0^te^{W(s)}ds=x\right.}
\end{equation}
\end{prop}
 \indim
The vector process $\bm Z(t)$ is Markovian and to its vector
\sde\refeq{vsde} it is possible to associate a (2+1)-dimensional
Fokker-Planck equation for the \pdf\ $f(\bm z,t)=f(x,y,t)$ with the
coefficients
\begin{equation*}
    \bm A(\bm z)=\bm a(\bm z)=\left(
               \begin{array}{c}
                 y \\
                 Dy \\
               \end{array}
             \right)
    \qquad\quad
    \mathbb{B}(\bm z)=2D\,\mathbb{C}(\bm z)\mathbb{C}^T(\bm z)=2D\left(
                                        \begin{array}{cc}
                                          0 & 0 \\
                                          0 & y^2 \\
                                        \end{array}
                                      \right)
\end{equation*}
giving rise to
\begin{eqnarray*}
    \partial_tf(x,y,t)\!\!&=&\!\!D\partial_y^2\left[y^2f(x,y,t)\right]-\partial_x\left[yf(x,y,t)\right]-\partial_y\left[Dyf(x,y,t)\right]
    \\
    \!\!&=&\!\!Dy^2\partial_y^2f(x,y,t)-y\left[\partial_xf(x,y,t)-3D\partial_yf(x,y,t)\right]+Df(x,y,t)
\end{eqnarray*}
namely\refeq{fpe}. This \pde\ essentially is confined to the
quadrant $x>0, y>0$ because the processes $X(t)$ and $Y(t)$ are
positive and never vanish. When the \pdf\ $f(x,y,t)$ has been found,
we can calculate the marginals $f_X(x,t)$ and $f_Y(y,t)$ as
\begin{equation*}
    f_X(x,t)=\int_0^{+\infty}f(x,y,t)\,dy\qquad\qquad
    f_Y(y,t)=\int_0^{+\infty}f(x,y,t)\,dx
\end{equation*}
As a matter of fact, however, from the laws of a Wiener process
$W(t)\sim\norm(0,2Dt)$ we already know that the log-normal \pdf\ of
$Y(t)$ which is
\begin{equation}\label{lgn}
    f_Y(y,t)=\frac{e^{-\frac{\ln^2y}{4Dt}}}{y\sqrt{4\pi Dt}}
\end{equation}
and hence we can also perform a first check of the coherence of our
joint equation\refeq{fpe}: by $x$-marginalization of\refeq{fpe} we
indeed have
\begin{eqnarray}
    \partial_tf_Y(y,t)&=&Dy^2\partial_y^2f_Y(y,t)+3Dy\partial_yf_Y(y,t)+Df_Y(y,t)-y\int_0^{+\infty}\partial_xf(x,y,t)\,dx\nonumber
    \\
    &=&Dy^2\partial_y^2f_Y(y,t)+3Dy\partial_yf_Y(y,t)+Df_Y(y,t)-y\left[f(x,y,t)\right]_{x=0}^{x=+\infty}\nonumber\\
    &=&Dy^2\partial_y^2f_Y(y,t)+3Dy\partial_yf_Y(y,t)+Df_Y(y,t)+y
    f(0,y,t)\label{xmarg}
\end{eqnarray}
that coincides with\refeq{yfpe} provided that $f(0,y,t)=0$, as it is
reasonable to require. On the other hand, by taking\refeq{lgn} into
account, we also have by direct calculation that
\begin{equation*}
    \frac{\partial_tf_Y(y,t)}{f_Y(y,t)}=\frac{\ln^2y-2Dt}{4Dt^2}=\frac{Dy^2\partial_y^2f_Y(y,t)+3Dy\partial_yf_Y(y,t)}{f_Y(y,t)}+D
\end{equation*}
so that the \pdf\refeq{lgn} of $Y(t)$ actually is a solution of the
$x$-marginalized equation\refeq{xmarg}. In the same vein we can then
study the $y$-marginalized equation of\refeq{fpe}: let us first
remark that
\begin{equation*}
    \left[yf(x,y,t)\right]_{y=0}^{y=\infty}=\left[y^2f(x,y,t)\right]_{y=0}^{y=\infty}=0
\end{equation*}
because all the moments of a log-normal distribution (which is the
marginal of $f(x,y,t)$) are finite so that $y^nf(x,y,t)$ must be
infinitesimal for $y\to+\infty$, while $f(x,y,t)$ can diverge in
$y\to0$ to an order strictly lesser than 1. As a consequence by
$y$-marginalizing\refeq{fpe} with integrations by part (the finite
terms vanish) we have the continuity equation\refeq{conteq}
\begin{eqnarray*}
  \partial_tf_X(x,t) &=& D\int_0^\infty\!\! y^2\partial_y^2f(x,y,t)\,dy-\partial_x\int_0^\infty\!\! y f(x,y,t)\,dy\\
   &&\qquad\qquad\qquad\qquad+3D\int_0^\infty\!\! y\partial_yf(x,y,t)\,dy+Df_X(x,t) \\
   &=&2Df_X(x,t)-\partial_x\left[f_X(x,t)\int_0^\infty\!\! y\,
   \frac{f(x,y,t)}{f_X(x,t)}\,dy\right]\\
   &&\qquad\qquad\qquad\qquad\qquad\qquad-3Df_X(x,t)+Df_X(x,t)\\
   &=&-\partial_x[v(x,t)f_X(x,t)]
\end{eqnarray*}
where $v(x,t)$ is defined in\refeq{vfield}
 \findim

\noindent All this is hardly surprising: we would get the same
continuity equation for the vector process constituted by the pair
position/velocity of a Brownian motion in the Ornstein-Uhlenbeck
dynamical model. The continuity equation\refeq{conteq} however is
not very useful for us because the velocity field $v(x,t)$ is in
some sense dependent from the form of the solution: the defining
conditional expectation\refeq{vfield} is indeed calculated with a
law also involving the marginal $f_X(x,t)$. In other words a well
behaved \pdf\ always satisfies a continuity equation when the
velocity field is rightly defined as in\refeq{vfield}, and hence it
expresses a consistence requirement, rather than a true equation ...
unless you already know the (well behaved) velocity field, and hence
the solution. On the other hand we do not have here an explicit
expression of $v(x,t)$, but it would be interesting to remark here
that its definition\refeq{vfield} seems to hint to the need of some
kind of \emph{mean-conditioning} since after all the condition is
expressed in terms of a sum (integral) of \rv's of the same kind of
the averaged one: this is a point worth of a further enquiry

Alternatively, by postponing every marginalization, we could try
first to solve the joint Fokker-Planck equation\refeq{fpe}:
separating the variables with $f(x,y,t)=g(x,y)h(t)$ we have
\begin{equation*}
    \frac{\dot{h}(t)}{h(t)}=\frac{Dy^2\partial_y^2g(x,y)-y\left[\partial_xg(x,y)-3D\partial_yg(x,y)\right]+Dg(x,y)}{g(x,y)}=\lambda
\end{equation*}
and hence
\begin{align*}
   & \dot{h}(t) = \lambda h(t)\\
   &
   Dy^2\partial_y^2g(x,y)-y\left[\partial_xg(x,y)-3D\partial_yg(x,y)\right]+(D-\lambda)
   g(x,y)=0
\end{align*}
Take then $g(x,y)=u(x)v(y)$ to have
\begin{equation*}
    \frac{u'(x)}{u(x)}=Dy\,\frac{v''(y)}{v(y)}+3D\,\frac{v'(y)}{v(y)}+\frac{D-\lambda}{y}=-\mu
\end{equation*}
and finally with $f(x,y,t)=u(x)v(y)h(t)$
\begin{align*}
    &  \dot{h}(t)  = \lambda h(t)\\
    &  u'(x)  = -\mu u(x) \\
    &  Dy^2\,v''(y)+3Dy\,v'(y)+(D-\lambda+\mu y)v(y)=0
\end{align*}
From the first two equations we simply have
\begin{equation*}
    h(t)=a e^{\lambda t}\qquad\quad u(x)= b e^{-\mu x}
\end{equation*}
while the third can be written as
\begin{equation*}
    v''(y)+\frac{3}{y}\,v'(y)+\left(\frac{D-\lambda}{Dy^2}+\frac{\mu}{Dy}\right)v(y)=0
\end{equation*}
Comparing now this equation with the Bessel equation 8.491.12
in\mycite{grad}, we see that within the notations adopted there we
must require
\begin{equation*}
    2\alpha-2\beta\nu+1=3\quad\qquad\beta=\frac{1}{2}\quad\qquad\beta^2\gamma^2=\frac{\mu}{D}\quad\qquad\alpha(\alpha-2\beta\nu)=\frac{D-\lambda}{D}
\end{equation*}
namely
\begin{equation*}
    \alpha=1\pm\sqrt{\frac{\lambda}{D}}\quad\qquad\beta=\frac{1}{2}\quad\qquad\gamma=\pm2\sqrt{\frac{\mu}{D}}\quad\qquad\nu=\pm2\sqrt{\frac{\lambda}{D}}
\end{equation*}
so that the solutions will have the form
\begin{equation*}
    v(y)=\frac{1}{y}\,Z_{\pm2\sqrt{\frac{\lambda}{D}}}\left(\pm2\sqrt{\frac{\mu y}{D}}\,\right)
\end{equation*}
where the symbol $Z_\nu(z)$ denotes one of the Bessel functions $J,
N, H^{(1)}, H^{(1)}$, as any linear combination of them. This seems
to confirm the relation with the Bessel functions found in the
Section\myref{chfexpansion}.

Several remarks, however, are in order: first, to keep $\alpha$ and
$\nu$ real we must require $\lambda \geq0$; on the other hand if
$\lambda>0$ the factor $h(t)$ will result in a time-exploding term,
so that the most reasonable option seems to be $\lambda =0$
(stationary solution). This choice would result in the $Z_0(z)$
Bessel functions, but empirical evidence (from \textit{Mathematica})
seems to imply that either $y^{-1}Z_0(\sqrt{y})$ is not always
non-negative, or it diverges for $z\to+\infty$, and in any case it
diverges for $y\to0$ in a non integrable way so that the
normalization of these functions appears to be hopeless. For
example, the unique Bessel function giving rise to non negative,
asymptotically infinitesimal solutions is $K_0(z)$, but
$y^{-1}K_0(\sqrt{y})$ diverges for $y\to0$ at an order $1+\epsilon$
with $\epsilon>0$ arbitrarily small. This is in any case coherent
with the remark that a stationary solution is hardly conceivable for
processes based on the exponentials of a Wiener process

\section{Conclusions and outlooks}\label{concl}

The present paper summarizes both a rather conclusive discussion of
the solutions of the Gompertz \sde\refeq{gompdl}, and several
partial results related to the solutions of the logistic
\sde\refeq{logistdl}. By postponing to future enquiries the
completion of this program, the definition of a viable deterministic
coarse-graining of these equations and a connection between their
solutions and the Nelson stochastic mechanics recalled in the
Appendix\myref{stochmech}, we will end our discussion by listing
here a few among the many points that would deserve a further
elaboration
\begin{enumerate}
%    \item \textbf{Stochastic Gompertz systems with time-dependent
%    coefficients:} A stochastic Gompertz system is a process
%    $X(t)$ solution of a \sde\refeq{gompdl}. A first problem then
%    could be:
%    \emph{what happens when we take a time-dependent, instead of a
%    constant,
%    coefficient $\alpha$} in this \sde\refeq{gompdl}
    \item \textbf{Random parametric Gomperts \sde's:}\label{point2} The parametric Gompertz
    \sde\
    \refeq{paramgomp} with a time-dependent frequency $\alpha(t)$ is associated to the parametric \ou\
    \sde\refeq{paramou}: ask then {what happens when $\alpha(t)$ is also random}, for instance of the
    type $\alpha(t)=\alpha_0(1+U(t))$ where $U(t)$ is a suitable
    external process (for instance either another Wiener process independent from $W(t)$, or $W(t)$
    itself). The case of random coefficients can be compared to the
    case of systems of \sde\ where the second \sde\ defines the new
    process $U(t)$. The two equations can be either coupled, or
    uncoupled
    \item \textbf{Modified Gompertzian growth:} {Can we devise some
    modified Gompertz equation giving rise -- beside the usual growth -- either to oscillations,
    or to some decrease to some other stable level?} This could be
    done by means of two mechanisms: either an external forcing term
    inscribed in the time depending (but not necessarily random)
    coefficients; or some shrewd random term (as in the previous
    point\myref{point2}) ruling, for instance, in an unpredictable way the
    sign of the exponentials of the Gompertz functions. This enquiry
    could lead to compare these systems with the random
    {\emph{Lotka-Volterra}} systems where the oscillations are induced by a
    coupling in a system of equations describing populations: should
    we think, then, to coupled Gompertzian systems?
    \item \textbf{Two kinds of deterministic correpondence:}\label{point1} Take the \sde
    \begin{equation*}
        dX(t)=a(X(t),t)dt+b(X(t),t)dW(t)
    \end{equation*}
    we can {recover a deterministic equation} in two, non-equivalent
    ways: either we can just consider that the diffusion coefficient
    $D$ vanishes, drop the random term and consider the
    \ode\footnote{This in some sense reverses the usual procedures leading to the Langevin equation starting from a Newton equation}
    \begin{equation*}
        \dot{x}(t)=a(x(t),t)
    \end{equation*}
    or we can take the expectation of the \sde
    \begin{equation*}
        \frac{d\EXP{X(t)}}{dt}=\EXP{a(X(t),t)}
    \end{equation*}
    which is rather different from the previous one because
    $\EXP{a(X(t),t)}$ is not $a(\EXP{X(t)},t)$, but for the linear
    case. The previous ambiguity
    is not relevant indeed for linear \sde\ of the Smoluchowsky type
    \begin{equation*}
        dX(t)=[q(t)-p(t)X(t)]dt+dW(t)
    \end{equation*}
    whose deterministic counterpart in both the cases is
    \begin{equation*}
        \dot{x}(t)+p(t)x(t)=q(t),\qquad x(t)=\EXP{X(t)}
    \end{equation*}
    namely the most general, first order, linear \ode\ whose solutions are
    completely known. On the other hand also the general solution of
    the corresponding \sde\ is completely known, so that it would be
    telling to compare first these two solutions. Contrariwise the problem for
    the non linear $a(x,t)$ still stay with us
    \item \textbf{Equations for the medians $\MED{X(t)}$:} In the quest
    for the deterministic counterpart of a \sde\ it could be useful to
    look to the medians (see Appendix\myref{quantmed}) rather than to the expectations. It is not easy however
    to find the equations for the medians of a
    process satisfying even the simple linear \sde
    \begin{equation*}
        dY(t)=[b(t)-a(t)Y(t)]dt+dW(t)
    \end{equation*}
    because in general the medians of sums ar not sums of the medians
    \item \textbf{The $n=1$ state of a quantum harmonic oscillator:}
    In the Appendix\myref{stochmech} we explicitly solved the \fpe\ associated by the stochastic mechanics to the
    $n=1$ eigenstate of a quantum harmonic
    oscillator, but at the same time we were not able to solve the
    corresponding \sde. As a matter of fact we were not able to find
    a process displaying as a transition \pdf\ that derived from the
    harmonic oscillator solution. The problem is that this solution looks like a
    mixture (and even in this case: what about the \rv\ with a
    mixture law in terms of the \rv's obeying the separate mixing
    laws?), but in fact it is not a true mixture because the coefficients sum up to 1,
    but are not in the range $[0,1]$ as should be for every convex
    combination. It looks then more as an affine combination, but nothing
    seems to be known about this case. Alternatively we could consider our
    solution as  a convex combination of \pdf\ taking also negative
    values. Either way we seem to be obliged
    to take seriously the existence of (at least virtual) negative
    probabilities\mycite{mueck,feynm1,feynm2} and follow this  path to its bitter end trying to
    guess what kind of \rv's -- if any -- can be distributed in this way
    \item \textbf{Logistic systems:} Finally, even if we choose to focus our attention on
    the more manageable Gompertz systems, at least from a mathematical standpoint it
    would be relevant to complete the investigation about the solutions of the
    Logistic equations that we gave in two different explicit
    ways\refeq{logisttranssol} and\refeq{logistpdf}:
    these, however, are for the time being more formal than substantial because, as we pointed out in
    the corresponding sections, we did not yet explicitly elaborate the laws of the
    involved \rv's that will be instead scrutinized in a forthcoming paper
\end{enumerate}

\newpage

\begin{appendix}

\section{Deterministic Gompertz and logistic
equations}\label{determ}

In the present appendix we will briefly recall the deterministic
variants of the \sde's presented all along the present paper: we
refer in particular to the dimensionless
\sde's\refeq{gompdl},\refeq{logistdl} and\refeq{thetalogdl}. These
deterministic versions are attained simply by switching off the
Wienerian noise $W(t)$ of the said \sde's (namely by taking $D=0$),
and are first order \emph{ODE}'s (\emph{ordinary differential
equations}) of the type
\begin{equation*}
    \dot{x}(t)=v\big[x(t),\,t\big]
\end{equation*}
Their solutions will be used in the present paper as a benchmark to
check the consistency of the solutions of the corresponding \sde's.
The reason why in the said \sde's the noise turns out to be also
proportional to the process itself will be addressed in a
forthcoming paper devoted to the deterministic ad stochastic models
of population dynamics

\subsection{Gompertz equation}\label{detgomp}

The deterministic version of the Gompertz \sde\refeq{gompdl}
apparently is
\begin{equation*}
    \dot{x}(t)=x(t)\big(1-\alpha\ln x(t)\big)
\end{equation*}
that, with an initial condition $x(0)=x_0$, is easily integrated by
separating the variables providing the solution
\begin{equation}\label{gompdet}
    x(t)=x_0^{e^{-\alpha t}}e^{(1-e^{-\alpha t})/\alpha}
\end{equation}
that is a monotonic function starting from $x_0$ and asymptotically
relaxing to the limiting value $e^{1/\alpha}$. This asymptotic
value, however, is always larger than $1$ because we require
$\alpha>0$: we could nevertheless outflank this limitation by adding
a new parameter $\beta$ into the equation
\begin{equation*}
    \dot{x}(t)=x(t)\big(\beta-\alpha\ln x(t)\big)
\end{equation*}
which now has the solution
\begin{equation*}
    x(t)=x_0^{e^{-\alpha t}}e^{\beta(1-e^{-\alpha t})/\alpha}
\end{equation*}
with an asymptotic value $e^{\beta/\alpha}$ that can now be both
larger and smaller than $1$ according to the sign of $\beta$

\subsection{Logistic and $\theta$-logistic
equations}\label{detlogist}

The logistic \emph{ODE} associated to the \sde\refeq{logistdl} is
\begin{equation*}
    \dot{x}(t)=x(t)\big(1-x(t)\big)
\end{equation*}
with a $v(x)=x(1-x)$, and its solution with the initial condition
$x(0)=x_0$ is again retrieved by separating the variables:
\begin{equation}\label{logistdet}
    x(t)=\frac{x_0}{x_0+(1-x_0)e^{-t}}
\end{equation}
Here too we are dealing with a function monotonically going, for
$t\to+\infty$, from an arbitrary $x_0$ to $1$. The $\theta$-logistic
\emph{ODE} is finally its generalization with $\theta>0$
\begin{equation*}
    \dot{x}(t)=x(t)\big(1-x^\theta(t)\big)
\end{equation*}
whose solution is easily seen to be now
\begin{equation}\label{thetalogdet}
    x(t)=\left(\frac{x_0^\theta}{x^\theta_0+(1-x^\theta_0)e^{-\theta
    t}}\right)^{1/\theta}
\end{equation}
with the same qualitative behavior as before

\section{Solving stochastic differential equations}\label{theito}

\subsection{An epitome of \ito\ calculus}\label{intro}

First of all let us recall in a simplified, mnemonic form a few
results of the \ito\ calculus that can be found for example
in\mycite{gs}, p.\ 11-27. Leaving aside a rigorous presentation, we
will work here in the setting of the \emph{stochastic differentials}
(\sdiff) whose precise meaning must be retrieved from their role and
use in the standard definition of the \ito\ integral
\begin{equation*}
    \int_a^bX(t)\,dW(t)
\end{equation*}
where $W(t)$ is a Wiener process with diffusion coefficient $2D$,
while $X(t)$ is a well behaved processes non-anticipative w.r.t.\
$W(t)$\mycite{gs}. In this framework the most relevant innovation
w.r.t.\ the usual differential calculus is the fact that the
differential $dW(t)$ of the Wiener process no longer can be deemed
to be of the order $dt$, but we will have instead
\begin{equation}\label{wiendiff}
    \left[dW(t)\right]^2=2Ddt\qquad\EXP{dW(t)}=0\qquad\EXP{dW(t)dW(s)}=2D\delta(s-t)\,ds\,dt
\end{equation}
that are shorthand, symbolic notations for the following integral
results:
\begin{eqnarray*}
    \int_a^bX(t)\,[dW(t)]^2&=&2D\int_a^bX(t)\,dt \\
    \EXP{\int_a^bX(t)\,dW(t)}&=&0 \\
    \EXP{\int_a^bX(t)\,dW(t)\int_a^bY(s)\,dW(s)}&=&2D\int_a^b\EXP{X(t)Y(t)}\,dt
\end{eqnarray*}
Further conceivable differentials, as for instance
\begin{equation}\label{horddiff}
    dW(t)dt\qquad\quad [dW(t)]^2dt\qquad\quad dW(t)(dt)^2\qquad\quad \ldots
\end{equation}
will instead behave as higher order infinitesimals and will then be
treated as zero in the sense that the corresponding \ito\ integrals
vanish. These simple rules will be enough to formally deduce all the
relevant results needed in our discussion.

In the following our one-dimensional processes $X(t), Y(t)\ldots$
defined for $t\in[0,T]$ will be taken in $H_2$, namely in the space
of the processes such that
\begin{equation*}
    \PR{\int_0^TX^2(t)\,dt<+\infty}=1
\end{equation*}
in order to make sure that -- within a framework of suitable
hypotheses pointed out in\mycite{gs} -- the \ito\ integral is well
defined, and that in general our processes will have (again
symbolically, as a surrogate for \ito\ integral expressions)
\sdiff's of the type
\begin{equation}\label{sdiff}
    dX(t)=A(t)dt+B(t)dW(t)
\end{equation}
where again the coefficients $\sqrt{|A(t)|}$ and $B(t)$ are taken in
$H_2$: of course the Wiener process $W(t)$ is retrieved when $A=0$
and $B=1$. Remark that when a process has the \sdiff\refeq{sdiff}
its infinitesimals follow rules similar to that of the Wiener
process; in particular, as can be easily seen from\refeq{wiendiff}
and\refeq{horddiff}, we have
\begin{equation}\label{sqdiff}
    [dX(t)]^2=2DB^2(t)dt
\end{equation}
while further powers of $dX(t)$ and $dt$ would be higher order
infinitesimals, and hence will be neglected. In this context we also
say that, with suitable initial conditions, a process $X(t)$
satisfies a \emph{stochastic differential equations} (\sde)
\begin{equation}\label{sde}
    dX(t)=a\left(X(t),t\right)dt+b\left(X(t),t\right)dW(t)
\end{equation}
when for every $t\in[0,T]$ it has the previous \sdiff\ whose
coefficients $A(t)=a(X(t),t)$ and $B(t)=b(X(t),t)$, $\,\Pqo$\ are
contingent on the process $X(t)$ himself. In particular when the
functions $a(x)$ and $b(x)$ are time-independent the \sde
\begin{equation}\label{sdetind}
    dX(t)=a\left(X(t)\right)dt+b\left(X(t)\right)dW(t)
\end{equation}
may also admit stationary solutions
 \begin{prop}
\textbf{\emph{(\ito\ formula)}} If $X(t)$ has the
\emph{\sdiff}\refeq{sdiff} and $h(x,t)$ is a fairly differentiable
function, then the process $Y(t)=h(X(t),t)$ will have the
\emph{\sdiff}
\begin{eqnarray}\label{ito}
   d\,Y(t)&=&d\,h(X(t),t)\nonumber\\
          &=&\big[h_t(X(t),t)+h_x(X(t),t)\,A(t)+D\,h_{xx}(X(t),t)\,B^2(t)\big]\,dt\\
          &&\qquad\qquad\qquad\qquad\qquad\qquad\qquad\qquad+\,h_x(X(t),t)\,B(t)\,dW(t)\nonumber
\end{eqnarray}
When in particular $h(x)$ is time independent the \ito\ formula
becomes
\begin{eqnarray}\label{dhX}
   d\,Y(t)&=&d\,h(X(t))\\
          &=&\big[h'(X(t))\,A(t)+D\,h''(X(t))\,B^2(t)\big]\,dt+h'(X(t))\,B(t)\,dW(t)\nonumber
\end{eqnarray}
 \end{prop}
 \indim
The formula\refeq{dhX} can be proved by using the Taylor formula for
$h$ together with\refeq{sqdiff}, and neglecting the infinitesimals
of order larger than 1:
\begin{eqnarray*}
  dh(X(t)) &=& h(X(t+dt))-h(X(t))= h(X(t)+dX(t))-h(X(t)) \\
   &=&\sum_{k\geq0}\frac{h^{(k)}(X(t))}{k!}\,[dX(t)]^k-h(X(t))\\
   &=&h(X(t))+h'(X(t))dX(t)+\frac{1}{2}h''(X(t))\,[dX(t)]^2\\
   &&\qquad\qquad\qquad\qquad\qquad+\sum_{k>2}\frac{h^{(k)}(X(t))}{k!}\,[dX(t)]^k-h(X(t))\\
   &=&h'(X(t))\,[A(t)dt+B(t)dW(t)]+h''(X(t))\,DB^2(t)dt\\
   &=&\left[h'(X(t))A(t)+D\,h''(X(t))B^2(t)\right]dt+h'(X(t))B(t)\,dW(t)
\end{eqnarray*}
The more general\refeq{ito} can be proved -- with minor additions --
along the same lines. These results coincide with those of the usual
calculus but for the second derivative terms that show up as soon as
the diffusion coefficient $D$ does not vanish
 \findim
 \begin{prop}
Take two \emph{\sdiff}'s on the {same} Wiener process $W(t)$
\begin{equation*}
    dX(t)=A_X(t)dt+B_X(t)dW(t)\qquad\qquad dY(t)=A_Y(t)dt+B_Y(t)dW(t)
\end{equation*}
then for their product we have the rule
\begin{equation}\label{diffprod}
    d\,\big[X(t)Y(t)\big]=X(t)\,dY(t)+Y(t)\,dX(t)+2B_X(t)B_Y(t)D\,dt
\end{equation}
 \end{prop}
 \indim
To prove\refeq{diffprod} it would be enough to use\refeq{wiendiff},
by neglecting other higher order infinitesimals:
\begin{eqnarray*}
  d\,[X(t)Y(t)] &=& X(t+dt)Y(t+dt)-X(t)Y(t)\\
  &=&[X(t)+dX(t)][Y(t)+dY(t)]-X(t)Y(t) \\
   &=&X(t)dY(t)+Y(t)dX(t)+dX(t)dY(t)\\
   &=& X(t)dY(t)+Y(t)dX(t)+B_X(t)B_Y(t)[dW(t)]^2\\
   &=&X(t)dY(t)+Y(t)dX(t)+2B_X(t)B_Y(t)Ddt
\end{eqnarray*}
which again shows an extra term with respect to the usual calculus.
 \findim

\subsection{\sde\ transformations}

When a process $X(t)$ satisfies the \sde\refeq{sde} we can often
exploit the \ito\ formula\refeq{ito} to transform this \sde\ in a
more manageable form (see\mycite{gs} p 33-39): take a fairly
differentiable, monotonic (in $x$) function $y=h(x,t)$ and for every
$t$ denote $x=g(y,t)$ its spatial inverse, namely
\begin{equation*}
    h\big(g(y,t),t\big)=y\qquad\qquad g\big(h(x,t),t\big)=x
\end{equation*}
Then, if $X(t)$ satisfies the \sde\refeq{sde} the \ito\
formula\refeq{ito} implies that the transformed process
$Y(t)=h(X(t),t)$ will satisfy a new \sde
\begin{equation}\label{newsde}
    dY(t)=\ha \big(Y(t),t\big)\,dt+\hb \big(Y(t),t\big)\,dW(t)
\end{equation}
where now
\begin{eqnarray}
  \ha (y,t) &=& \big[h_t(x,t)+h_x(x,t)a(x,t)+D\,h_{xx}(x,t)b^2(x,t)\big]_{x=g(y,t)}\label{itoy1} \\
  \hb (y,t) &=& \big[h_x(x,t)b(x,t)\big]_{x=g(y,t)}\label{itoy2}
\end{eqnarray}
or else, in an equivalent form,
\begin{eqnarray}
  \ha (h(x,t),t)\!\! &=& \!\!h_t(x,t)+h_x(x,t)a(x,t)+D\,h_{xx}(x,t)b^2(x,t)\label{itox1} \\
  \hb (h(x,t),t)\!\! &=& \!\!h_x(x,t)b(x,t)\label{itox2}
\end{eqnarray}
When on the other hand $X(t)$ satisfies the \sde\refeq{sdetind} with
time-independent coefficients, the transformation $Y(t)=h(X(t))$
with monotonic $y=h(x)$ and $x=g(y)=h^{-1}(y)$ leads the \sde
\begin{equation}\label{newsdetind}
    dY(t)=\ha \left(Y(t)\right)dt+\hb \left(Y(t)\right)dW(t)
\end{equation}
where now
\begin{eqnarray}
  \ha (y) &=& \big[a(x)h'(x)+D\,b^2(x)h''(x)\big]_{x=g(y)}\label{itoy1tind} \\
  \hb (y) &=& \big[b(x)h'(x)\big]_{x=g(y)}\label{itoy2tind}
\end{eqnarray}
namely, in an equivalent form,
\begin{eqnarray}
  \ha (h(x)) &=& a(x)h'(x)+D\,b^2(x)h''(x)\label{itox1tind} \\
  \hb (h(x)) &=& b(x)h'(x)\label{itox2tind}
\end{eqnarray}
We are now interested in transforming a \sde\refeq{sde} into a new
form\refeq{newsde} that turns out to be simpler. To this end we will
analyze first what type of coefficients $\ha (y,t)$ and $\hb (y,t)$
lead to elementary solution, an then under which conditions the
original coefficients $a(y,t)$ and $b(y,t)$ are susceptible to be
transformed into that new form

\subsubsection{Elementary solvable \sde's}\label{elem}

\begin{enumerate}
    \item \noindent \textbf{Constant coefficients $\ha$ and $\hb$}:
In this case the \sde\refeq{newsde} becomes
\begin{equation}\label{sdeconst}
    dY(t)=\ha\,dt+\hb\,dW(t)
\end{equation}
and its solution, with $\Pqo$ initial condition $W(0)=0$ and
$Y(0)=Y_0$, apparently is
\begin{equation}\label{sdeconstsol}
    Y(t)=Y_0+\ha\,t+\hb\,W(t)
\end{equation}
Then the process $\overline{Y}(t)=Y(t)-Y_0\sim\norm\left(\ha
t\,,2D\hb t\right)$ is a Wiener process with a $\hb$ re-scaled
diffusion coefficient and a constant, deterministic drift $\ha$. For
a degenerate condition $Y_0=y,\Pqo$ at a time $s$ the law of the
process is then the transition \pdf
\begin{equation*}
    \norm\left(y+\ha (t-s),\,2D\hb (t-s)\right)
\end{equation*}
We will see in the following, however, that the compatibility
conditions on $a(x,t)$ and $b(x,t)$ to be both transformable into
new constant coefficients are in general excessively tight, so that
the corresponding solutions discussed above are indeed of little
practical interest
    \item \textbf{Smoluchowsky equation}: We consider next the case
of transformations leading to $\hb(y,t)=1$, namely to
\sde\refeq{newsde} of the type
\begin{equation}\label{smol}
    dY(t)=\ha(Y(t),t)\,dt+dW(t)
\end{equation}
that can always be considered as a Smoluchowsky equation
(see\mycite{nels}, ch. 10), namely as the overdamped
($\gamma\to+\infty$) limit of an \ou\ (Ornstein-Uhlenbeck) system of
dynamic (Newton) \sde's for a derivable process $Z(t)$
\begin{eqnarray*}
  dZ(t) &=& V(t)\,dt \\
  dV(t) &=& \gamma\big[\ha(Z(t),t)-V(t)\big]\,dt+\gamma\,dW(t)
\end{eqnarray*}
where $\ha(z,t)$ now represents an external field of forces, and the
approximation is understood in the sense that
\begin{equation*}
    \lim_{\gamma\to+\infty}Z(t)=Y(t)\qquad\quad\Pqo
\end{equation*}
There is no simple, general way to solve even the
equation\refeq{smol}, but when $\ha(y)$ is time independent the
\sde\refeq{smol} takes the form
\begin{equation}\label{sdeb1}
    dY(t)=\ha(Y(t))\,dt+dW(t)
\end{equation}
where $\ha(z)$ is now a \emph{stationary} external field of forces,
so that we can reasonably hope to find at least some
\emph{stationary} solution. If we introduce indeed a potential
$\phi(y)$ according to
\begin{equation}\label{potential}
    \ha (y)=-\frac{D}{kT}\,\phi'(y)
\end{equation}
(here $k$ is the Boltzmann constant and $T$ the absolute
temperature) then it easy to show that there is a stationary
solution with an invariant Boltzmann distribution
\begin{equation}\label{boltz}
    \frac{e^{-\frac{\phi(y)}{kT}}}{\int_\R e^{-\frac{\phi(z)}{kT}}dz}
\end{equation}
provided that $\phi(y)$ is such that $e^{-\frac{\phi(y)}{kT}}$ is an
integrable function (for details see Appendix\myref{boltzdistr})
    \item \textbf{Process-independent coefficients $\ha(t)$ and
    $\hb(t)$}: When on the other hand $\ha (t)$ and $\hb (t)$ of the
\sde\refeq{newsde} are both $y$-independent the equation is reduced
to a \sdiff\ with deterministic coefficients
\begin{equation}\label{sdet}
    dY(t)=\ha (t)\,dt+\hb (t)\,dW(t)
\end{equation}
and its solution is just the integral of the \sdiff, namely
\begin{equation}\label{sdetsol}
    Y(t)=Y_0+\int_0^t\ha (u)du+\int_0^t\hb (u)dW(u)
\end{equation}
with the initial condition $Y(0)=Y_0,\,\Pqo$ This means that the
process $\overline{Y}(t)=Y(t)-Y_0$ is Gaussian
\begin{equation*}
    \overline{Y}(t)\sim\norm\left(\int_0^t\ha (u)du\,,\,2D\!\int_0^t\hb \,^2(u)du\right)
\end{equation*}
and by taking a degenerate condition $Y_0=y,\,\Pqo$ at a time $s$
the law of the process is the transition \pdf\ $f(x,t|y,s)$
\begin{equation}\label{sdettrans}
    \norm\left(y+\int_s^t\ha (u)du\,,\,\,2D\!\int_s^t\hb \,^2(u)du\right)
\end{equation}
    \item \textbf{Process-linear coefficients}: We consider next the case of coefficients which are linear in $y$
 \begin{equation}\label{lincoeff}
    \ha(y,t)=\ha_0(t)+\ha_1(t)y\qquad\quad\hb(y,t)=\hb_0(t)+\hb_1(t)y
 \end{equation}
namely of the \sde's of the form
\begin{equation}\label{linear}
    dY(t)=\Big[\ha_0(t)+\ha_1(t)Y(t)\Big]\,dt+\left[\,\hb_0(t)+\hb_1(t)Y(t)\right]dW(t)
\end{equation}
and in particular its homogeneous counterpart with $\ha_0(t)=0$ and
$\hb_0(t)=0$
\begin{equation}\label{linearh}
    dY(t)=\ha_1(t)Y(t)dt+\hb_1(t)Y(t)dW(t)
\end{equation}
We look first for a solution of the \textbf{homogeneous
equation}\refeq{linearh} and -- as long as $Y(t)\neq0$ -- we define
a new process $Z(t)$ through the transformation
\begin{align*}
    &z=h(y,t)=\ln y\qquad y=g(z,t)=e^z\qquad\quad Z(t)=\ln Y(t)\qquad
    Y(t)=e^{Z(t)}\\
    &\qquad\qquad h_t(y,t)=0\qquad\quad h_y(y,t)=\frac{1}{y}\qquad\quad h_{yy}(y,t)=-\frac{1}{y^2}
\end{align*}
It is easy to see then from\refeq{itoy1} and\refeq{itoy2} that the
\sde\ for $Z(t)$ is
\begin{equation}\label{zeq}
    dZ(t)=\left[\ha_1(t)-D\,\hb\,^2_1(t)\right]dt+\hb_1(t)\,dW(t)
\end{equation}
with process independent coefficients as in\refeq{sdet}, and hence
from\refeq{sdetsol} the solution of\refeq{zeq} with $Z(0)=Z_0$ is
\begin{equation*}
    Z(t)=Z_0+\int_0^t\left[\ha_1(s)-D\,\hb\,^2_1(s)\right]ds+\int_0^t\hb_1(s)\,dW(s)
\end{equation*}
The (never vanishing) solution of the homogeneous
\sde\refeq{linearh} with $Y_0=e^{Z_0}$ then is
\begin{equation}\label{linearhsol}
    Y(t)=Y_0\,e^{\overline{Z}(t)}
\end{equation}
where we have introduced the process
\begin{equation}\label{linearhsolbar}
    \overline{Z}(t)=Z(t)-Z_0=\int_0^t\left[\ha_1(s)-D\,\hb\,^2_1(s)\right]ds+\int_0^t\hb_1(s)\,dW(s)
\end{equation}
which is again a solution of\refeq{zeq}, but with the initial
condition $Z(0)=0$. It turns out of course that $\overline{Z}(t)$ is
a Gaussian process with
\begin{equation*}
    \overline{Z}(t)\sim\norm\left(\int_0^t\left[\ha_1(s)-D\,\hb\,^2_1(s)\right]\!ds\,,\,\,2D\int_0^t\hb\,^2_1(s)ds\right)
\end{equation*}
Going back then to the \textbf{non homogeneous
equation}\refeq{linear} we first remark that from\refeq{dhX} and the
expression\refeq{zeq} of the \sdiff\ $d\overline{Z}(t)$ we have
\begin{eqnarray}
  d\left(e^{-\overline{Z}(t)}\right)\!\! &=&\!\! \left(-\left[\ha_1(t)-D\,\hb\,^2_1(t)\right]
        +D\,\hb\,^2_1(t)\right)e^{-\overline{Z}(t)}dt-\hb_1(t)e^{-\overline{Z}(t)}dW(t)\nonumber \\
   &=&\!\!\left[-\ha_1(t)+2D\,\hb\,^2_1(t)\right]e^{-\overline{Z}(t)}dt-\hb_1(t)e^{-\overline{Z}(t)}dW(t)\label{diffexp}
\end{eqnarray}
and then that from\refeq{diffprod},\refeq{linear} and\refeq{diffexp}
it is
\begin{eqnarray*}
  d\left(e^{-\overline{Z}(t)}Y(t)\right)\!\! &=&\!\!
  e^{-\overline{Z}(t)}dY(t)+Y(t)d\left(e^{-\overline{Z}(t)}\right)\\
   &&\qquad\qquad \qquad\qquad  -\left[\hb_0(t)+\hb_1(t)Y(t)\right]\hb_1(t)e^{-\overline{Z}(t)}2Ddt \\
   &=&\!\!e^{-\overline{Z}(t)}\left([\ha_0(t)+\ha_1(t)Y(t)]dt+[\hb_0(t)+\hb_1(t)Y(t)]dW(t)\right)\\
   &&\qquad\quad
             e^{-\overline{Z}(t)}Y(t)\left(\left[-\ha_1(t)+2D\hb\,^2_1(t)\right]dt-\hb_1(t)dW(t)\right)\\
   &&\qquad\qquad\qquad
   -\left[\hb_0(t)+\hb_1(t)Y(t)\right]\hb_1(t)e^{-\overline{Z}(t)}2Ddt
   \\
   &=&\!\!\left[\ha_0(t)-2D\,\hb_0(t)\hb_1(t)\right]e^{-\overline{Z}(t)}dt+\hb_0(t)e^{-\overline{Z}(t)}dW(t)
\end{eqnarray*}
By taking now into account the initial condition
$e^{-\overline{Z}(0)}Y(0)=Y(0)=Y_0$, the previous \sdiff\ can be
integrated as
\begin{equation*}
    e^{-\overline{Z}(t)}Y(t)=Y_0+\!\int_0^t\!\left[\ha_0(s)-2D\,\hb_0(s)\hb_1(s)\right]e^{-\overline{Z}(s)}ds
    +\int_0^t\hb_0(s)e^{-\overline{Z}(s)}dW(s)
\end{equation*}
so that the general solution of the non homogeneous
equation\refeq{linear} finally is
\begin{align}\label{linearsol}
    &Y(t)=Y_0e^{\overline{Z}(t)}+\int_0^t\left[\ha_0(s)-2D\,\hb_0(s)\hb_1(s)\right]e^{\overline{Z}(t)-\overline{Z}(s)}ds\\
    &\qquad\qquad\qquad\qquad\qquad\qquad\qquad\qquad+\int_0^t\hb_0(s)e^{\overline{Z}(t)-\overline{Z}(s)}dW(s)\nonumber
\end{align}
where the process $\overline{Z}(t)$ is defined
in\refeq{linearhsolbar}. We remark finally that in the case of
\textbf{time-independent coefficients}, namely when
\begin{align}
 & \quad \ha(y)=\ha_0+\ha_1y\qquad\qquad\hb(y)=\hb_0+\hb_1y \label{linearizedcoeff}\\
 &
 dY(t)=\bigl[\ha_0+\ha_1Y(t)\bigr]dt+\bigl[\,\hb_0+\hb_1Y(t)\bigr]dW(t)\label{linearized}
\end{align}
the process\refeq{linearhsolbar} becomes
\begin{equation}\label{linearhsolbar1}
    \overline{Z}(t)=(\ha_1-D\,\hb\,^2_1)\,t+\hb_1\,W(t)
\end{equation}
and the solution\refeq{linearsol} is reduced to
\begin{equation}\label{linearsolnot}
    Y(t)=Y_0e^{\overline{Z}(t)}+\!\left(\ha_0-2D\,\hb_0\hb_1\right)\int_0^t\!\!e^{\overline{Z}(t)-\overline{Z}(s)}ds
    +\!\hb_0\int_0^te^{\overline{Z}(t)-\overline{Z}(s)}dW(s)
\end{equation}
 \item \textbf{Time-independent coefficients $\ha(y)$ and
 $\hb(y)$}:\label{ticoeff}
There are no general formulas available, but -- as will be shown in
the Section\myref{constcoeff} -- it is always possible to manage the
transformation $Y(t)=h(X(t))$ in such a way that $\hb(y)=1$
(different constant values could be easily subsumed in a
redefinition of $W(t)$) so that the \sde\refeq{newsde} takes the
form of a time independent Smoluchowsky equation\refeq{sdeb1}. Even
so we already remarked however that, beside a possible stationary
solution, there is no simple, general way to solve the
equation\refeq{sdeb1}. If on the other hand we try to simplyfy our
problem by taking $\ha(y)=0$ and an arbitrary $\hb(y)$, then we are
led to the equation
\begin{equation}\label{sdea0}
    dY(t)=\hb(Y(t))\,dW(t)
\end{equation}
but again, in the general case, we have no simple solution to show
so that the problem must be dealt with on a case-by-case basis
\end{enumerate}

\subsubsection{Transformations to constant coefficients}\label{constcoeff}

We will first look for the transformations $y=h(x,t)$ leading to
constant coefficients  $\ha $ and $\hb$: from\refeq{itox1}
and\refeq{itox2} we have
\begin{eqnarray}
  \ha  &=& h_t(x,t)+h_x(x,t)a(x,t)+D\,h_{xx}(x,t)b^2(x,t)\label{itoxc1} \\
  \hb &=& h_x(x,t)b(x,t)\label{itoxc2}
\end{eqnarray}
From the equation\refeq{itoxc2}  -- as long as $b(x,t)$ does not
change its sign  as a function of $x$, and hence the resulting
$h(x,t)$ turns out to be monotonic and invertible in $x$ for every
$t$ -- we immediately get the transformation
\begin{equation}\label{2consth}
    h_x(x,t)=\frac{\hb}{b(x,t)}\qquad\quad
    y=h(x,t)=\hb\int\frac{dx}{b(x,t)}\qquad\quad x=g(y,t)
\end{equation}
and since this gives
\begin{equation*}
    h_{xx}(x,t)=-\hb\frac{b_x(x,t)}{b^2(x,t)}\qquad\quad
    h_t(x,t)=-\hb\int\frac{b_t(x,t)}{b^2(x,t)}\,dx
\end{equation*}
we also obtain from\refeq{itoxc1} the equation
\begin{equation*}
    \frac{\ha}{\hb}=-\int\frac{b_t(x,t)}{b^2(x,t)}\,dx+\frac{a(x,t)}{b(x,t)}-D\,b_x(x,t)
\end{equation*}
By $x$ derivation we then find that $a(x,t)$ and $b(x,t)$ must
satisfy the condition
\begin{equation*}
    a_x(x,t)-\frac{b_x(x,t)}{b(x,t)}\,a(x,t)=D\,b(x,t)b_{xx}(x,t)+\frac{b_t(x,t)}{b(x,t)}
\end{equation*}
that can be solved as a first order \ode\ (ordinary differential
equation) for $a(x,t)$ giving the explicit relation ($\alpha(t)$ is
an arbitrary integration function)
\begin{equation}\label{2constt}
    a(x,t)=\alpha(t)+b(x,t)\left[Db_x(x,t)+\int\frac{b_t(x,t)}{b^2(x,t)}\,dx\right]
\end{equation}
that for time independent $a(x)$ and $b(x)$ becomes ($\alpha$ here
is now an arbitrary integration constant)
\begin{equation}\label{2onst}
    a(x)=\alpha+Db(x)b'(x)
\end{equation}
These are however very tight conditions that can be verified only in
a few particular cases and hence we will rather go on discussing
whether at least one of the two coefficients -- either $\ha $, or
$\hb $ -- can be reasonably made constant

We begin then by requiring only $\hb (y,t)$ to be constant and,
without a loss of generality, we can take $\hb (y,t)=1$: every other
constant value can indeed be easily accounted for with a
redefinition of the diffusion coefficient of $W(t)$. The
transformation under discussion will then lead to the Smoluchowsky
\sde\refeq{smol} previously discussed in the Section\myref{elem}.
From\refeq{itoxc2} with $\hb=1$ we now get the previous
transformation\refeq{2consth} in the form
\begin{equation}\label{bconsth}
       y=h(x,t)=\int\frac{dx}{b(x,t)}\qquad\quad x=g(y,t)
\end{equation}
and then, instead of requiring the condition\refeq{itoxc1}, we must
recall\refeq{itoy1} to get the explicit expression for the new drift
coefficient
\begin{equation*}
    \ha
    (y,t)=-\left[\int\frac{b_t(x,t)}{b^2(x,t)}\,dx-\frac{a(x,t)}{b(x,t)}+D\,b_x(x,t)\right]_{x=g(y,t)}
\end{equation*}
giving rise to the time dependent Smoluchowsky equation\refeq{smol}.
When in particular $a(x)$ and $b(x)$ are time-independent the
transformation becomes
\begin{equation}\label{b1}
    h(x)=\int\frac{dx}{b(x)}
\end{equation}
and the new drift is reduced to
\begin{equation}\label{a1}
    \ha (y)=\left[\frac{a(x)}{b(x)}-D\,b'(x)\right]_{x=g(y)}
\end{equation}
so that the transformed \sde\ takes the form\refeq{sdeb1} discussed
in the point\myref{ticoeff} of the previous Section\myref{elem}. As
remarked therein, in this time-independent setting the drift
function\refeq{a1} represents a {stationary} dynamics for the
transformed process $Y(t)$, and when it is deduced from a potential
$\phi(y)$ it can lead to an \emph{invariant Boltzmann
distribution}\refeq{boltz} as stated in the Section\myref{elem} (for
details see Appendix\myref{boltzdistr})

If we instead require only $\ha (y,t)=\ha$ to be a constant, then
from\refeq{itox1} $h(x,t)$ must satisfy the \pde\ (partial
differential equation)
\begin{equation*}
    h_t(x,t)+a(x,t)h_x(x,t)+D\,b^2(x,t)h_{xx}(x,t)=\ha
\end{equation*}
which, if $a(x)$ and $b(x)$ are time-independent, is reduced to the
\ode
\begin{equation*}
    a(x)h'(x)+D\,b^2(x)h''(x)=\ha
\end{equation*}
that can be solved to give the explicit expression of the
transformation $h(x)$ and of the coefficient $\hb(x)$. For
simplicity, however, we will restrict ourselves to the case $\ha=0$:
for suitable integration constants $c_1$ and $c_2$ the general
solution is
\begin{equation*}
    h(x)=c_1+c_2\int e^{-\int\frac{a(x)}{Db^2(x)}\,dx}dx
\end{equation*}
and from\refeq{itoy2} the new diffusion coefficient will be
\begin{equation*}
    \hb(y)=c_2\left[b(x)e^{-\int\frac{a(x)}{Db^2(x)}\,dx}\right]_{x=g(y)}
\end{equation*}
that will enter into the new \sde\refeq{sdea0}. As already remarked,
however, we again have no general solutions available, so that these
\sde's could be tackled only case by case

\subsubsection{Transformations to process-independent coefficients}\label{spindep}

Given the simple results\refeq{sdetsol} and\refeq{sdettrans} it
would be interesting, as a next step, to find under what conditions
we can transform an arbitrary \sde\refeq{sde} into a new equation
with space-independent coefficients. To this end let us first remark
that the two equations\refeq{itox1} and\refeq{itox2} become in this
case
\begin{eqnarray}
  \ha (t) &=& h_t(x,t)+a(x,t)h_x(x,t)+D\,b^2(x,t)h_{xx}(x,t)\label{itot1} \\
  \hb (t) &=& h_x(x,t)b(x,t)\label{itot2}
\end{eqnarray}
From\refeq{itot2} we have
\begin{align}
   & h_x(x,t)=\frac{\hb (t)}{b(x,t)}\qquad\quad h_{xx}(x,t)=-\frac{\hb (t)}{b^2(x,t)}b_x(x,t)\label{itot3}\\
   & \qquad\qquad h_{xt}(x,t)=\frac{\dot{\hb}(t)b(x,t)-\hb (t)b_t(x,t)}{b^2(x,t)}\label{itot4}
\end{align}
while by deriving\refeq{itot1} with respect to $x$ we find
\begin{equation}\label{itoxt}
    h_{xt}(x,t)+\partial_x\left[a(x,t)h_x(x,t)+D\,b^2(x,t)h_{xx}(x,t)\right]=0
\end{equation}
Then by substituting\refeq{itot3} and\refeq{itot4}
into\refeq{itoxt}, after some manipulation we obtain
\begin{equation}\label{condt1}
    \frac{\dot{\hb}(t)}{\hb (t)}=b(x,t)\left[\frac{b_t(x,t)}{b^2(x,t)}-\partial_x\left(\frac{a(x,t)}{b(x,t)}\right)+D\,b_{xx}(x,t)\right]
\end{equation}
By deriving\refeq{condt1} with respect to $x$ again we finally get
\begin{equation}\label{condt}
    \partial_x\left\{b(x,t)\left[\frac{b_t(x,t)}{b^2(x,t)}-\partial_x\left(\frac{a(x,t)}{b(x,t)}\right)+D\,b_{xx}(x,t)\right]\right\}=0
\end{equation}
a condition -- involving only the initial coefficients $a(x,t)$ and
$b(x,t)$,  -- that must be satisfied in order to ensure the
existence of a transformation $h(x,t)$ bringing to an equation with
space-independent coefficients. In fact, when\refeq{condt} holds,
the r.h.s.\ of\refeq{condt1} must depend only on $t$, and hence this
equation enables us to find $\hb (t)$. Then by plugging this $\hb
(t)$ into\refeq{itot2} we can find the transformation $h(x,t)$, and
finally $\ha (t)$ comes out of\refeq{itot1} whose r.h.s.\ too will
turn out to be dependent only on $t$

When in particular the coefficients $a(x)$ and $b(x)$ are
time-independent the equation\refeq{condt1} becomes
\begin{equation}\label{condt2}
    \frac{\dot{\hb}(t)}{\hb (t)}=b(x)\left[D\,b''(x)-\frac{d}{dx}\left(\frac{a(x)}{b(x)}\right)\right]
\end{equation}
and since now the l.h.s.\ depends only on $t$, and the r.h.s\ only
on $x$, the two members of this equation must be both equal to a
constant $c$. Hence instead of\refeq{condt} the compatibility
condition coming out of\refeq{condt2} is just
\begin{equation}\label{compat}
    b(x)\left[D\,b''(x)-\frac{d}{dx}\left(\frac{a(x)}{b(x)}\right)\right]=c
\end{equation}
while from
\begin{equation}\label{bt}
    \frac{\hb \,'(t)}{\hb (t)}=c
\end{equation}
we first get $\hb (t)=e^{ct}$, then from\refeq{itot2} we have
\begin{equation}\label{ht}
    h(x,t)=e^{ct}\int\frac{dx}{b(x)}
\end{equation}
and finally from\refeq{itot1}
\begin{equation}\label{at}
    \ha (t)=e^{ct}\left[c\int\frac{dx}{b(x)}+\frac{a(x)}{b(x)}-D\,b'(x)\right]
\end{equation}
where the term in square brackets is in fact constant provided that
the compatibility condition\refeq{compat} is satisfied

\subsubsection{Transformations to process-linear
coefficients}\label{linapp}

Since we have shown that every \sde\refeq{linear} with space-linear
coefficients\refeq{lincoeff} can be explicitly solved, it will be
important to find under what conditions we can transform an
arbitrary \sde\refeq{sde} in the new form\refeq{linear}. We will
analyze in detail, however, only the case of
\textbf{time-independent coefficients} $a(x),\,b(x)$
and\refeq{linearizedcoeff} when the new \sde\ is supposed to take
the form\refeq{linearized}. In this case the transformation $h(x)$
is time-independent too, and the conditions\refeq{itox1}
and\refeq{itox2} become
\begin{eqnarray}
  \ha_0+\ha_1h(x) &=& h'(x)a(x)+D\,h''(x)b^2(x)\label{linreq1} \\
  \hb_0+\hb_1h(x) &=& h'(x)b(x)\label{linreq2}
\end{eqnarray}
If we choose first to have $\hb_1\neq0$, from the
equation\refeq{linreq2} we at once have
\begin{equation}\label{tratolin1}
    h(x)=c\,e^{\,\hb_1p(x)}-\frac{\hb_0}{\hb_1}\qquad\qquad
    p(x)=\int\frac{1}{b(x)}\,dx\qquad p'(x)=\frac{1}{b(x)}
\end{equation}
and hence
\begin{equation*}
    h'(x)=c\,\hb_1\frac{e^{\,\hb_1p(x)}}{b(x)}\qquad\quad
    h''(x)=c\,\hb_1\frac{e^{\,\hb_1p(x)}}{b^2(x)}\left[\hb_1-b'(x)\right]
\end{equation*}
By plugging all that into\refeq{linreq1} after some algebra we get
the equation
\begin{equation*}
    \frac{\ha_0\hb_1-\ha_1\hb_0}{c\,\hb_1}=e^{\,\hb_1p(x)}\left[\hb_1q(x)+D\,\hb\,^2_1-\ha_1\right]\qquad\qquad
    q(x)=\frac{a(x)}{b(x)}-D\,b'(x)
\end{equation*}
In order then to have a condition free from constants, by a first
derivation we get
\begin{eqnarray*}
    0&=&\frac{d}{dx}\left(e^{\,\hb_1p(x)}\left[\hb_1q(x)+D\,\hb\,^2_1-\ha_1\right]\right)\\
    &=&\hb_1e^{\,\hb_1p(x)}\left(p'(x)\left[\hb_1q(x)+D\,\hb\,^2_1-\ha_1\right]+q'(x)\right)
\end{eqnarray*}
namely from\refeq{tratolin1}
\begin{equation*}
    \ha_1-D\,\hb\,^2_1=\hb_1q(x)+b(x)q'(x)
\end{equation*}
With a second derivation we then have
\begin{equation*}
    0=\frac{d}{dx}\left[\hb_1q(x)+b(x)q'(x)\right]=\hb_1q'(x)+\frac{d}{dx}\left[b(x)q'(x)\right]
\end{equation*}
that can be recast as
\begin{equation}\label{lincond2}
    \hb_1=-\frac{1}{q'(x)}\frac{d}{dx}\left[b(x)q'(x)\right]
\end{equation}
and finally with a third derivation
\begin{equation}\label{lincond}
    \frac{d}{dx}\left(\frac{1}{q'(x)}\frac{d}{dx}\left[b(x)q'(x)\right]\right)=0
\end{equation}
which is now the wanted condition involving only the coefficients
$a(x)$ and $b(x)$ of the initial, non linear \sde\refeq{sde}.
When\refeq{lincond} is satisfied, then we can take\refeq{lincond2}
as the value of the parameter $\hb_1$, and $h(x)=ce^{\,\hb_1p(x)}$
for some suitable value of the constant $c$ as the transformation
able to reduce our equation to its linear form

If instead we require $\hb_1=0$, from the equation\refeq{linreq2}
and within the same notations we simply get
\begin{equation*}
    h(x)=\hb_0p(x)+c\qquad\quad h'(x)=\frac{\hb_0}{b(x)} \qquad\quad
    h''(x)=-\hb_0\frac{b\,'(x)}{b^2(x)}
\end{equation*}
and hence after some manipulation from\refeq{linreq1} we have
\begin{equation*}
    q(x)=\ha_1p(x)+c\ha_1+\frac{\ha_0}{\hb_0}
\end{equation*}
then by derivation we obtain $q'(x)=\ha_1p'(x)$ namely
$b(x)q'(x)=\ha_1$, and by further derivation we finally get the
condition
\begin{equation}\label{lincond3}
    \frac{d}{dx}[b(x)q'(x)]=0
\end{equation}
It is apparent that this condition\refeq{lincond3} also implies the
condition\refeq{lincond} which however remains the most general
requirement whose compliance is needed in order to be able to
transform the \sde\refeq{sde} into one with space-linear,
time-independent coefficients

\section{Solving Fokker-Planck equations}\label{thefpe}

Given the \sde\refeq{sde} the \pdf\ of its solutions can be obtained
by solving the (forward) Fokker-Planck equation
\begin{equation}\label{fpt}
    \partial_tf(x,t)=-\partial_x\left[a(x,t)f(x,t)\right]+D\,\partial^2_x\left[b^2(x,t)f(x,t)\right]
\end{equation}
that for time-independent coefficients (namely for an \sde\ of the
form\refeq{sdetind}) becomes
\begin{equation}\label{fp}
    \partial_tf(x,t)=-\partial_x\left[a(x)f(x,t)\right]+D\,\partial^2_x\left[b^2(x)f(x,t)\right]
\end{equation}
Standard solution methods are well known, as that of the
eigenfunction expansion discussed later. We will however first give
a look to a few results giving the transition \pdf\ without
solving\refeq{fpt}, but using a few results about the expectations
that however are not easily calculated explicitly

\subsection{Semi-explicit transition \pdf's}\label{transpdf}

We will first recall a formal procedure (see\mycite{gs}, \S 13, p.\
91) used to find the transition \pdf\ of the \fpe\refeq{fpt}, namely
of a process solution of a the \sde\refeq{sde}. Define -- in a
notation only partially coherent with our previous one -- the
functions
\begin{equation*}
    y=h(x,t)=\int\frac{dx}{b(x,t)}\qquad\quad x=g(y,t) \qquad\quad
    h(g(y,t),t)=y
\end{equation*}
then the function
\begin{eqnarray*}
  \ha(y,t) &=& h_t(g(y,t),t)+h_x(g(y,t),t)\,a(g(y,t),t)+Dh_{xx}(g(y,t),t)\,b^2(g(y,t),t) \\
   &=&\left[-\int\frac{b_t(x,t)}{b^2(x,t)}\,dx+\frac{a(x,t)}{b(x,t)}-Db_x(x,t)\right]_{x=g(y,t)}
\end{eqnarray*}
and finally
\begin{equation*}
    \alpha(y,t)=\frac{1}{2D}\int\ha(y,t)\,dy
    \qquad\quad \beta(y,t)=-\frac{\ha^2(y,t)}{4D}-\frac{\ha_y(y,t)}{2}-\frac{1}{2D}\int\ha_t(y,t)dy
\end{equation*}
Consider now for $0\leq s\leq t$ the two-times process
\begin{equation*}
    Z(s,t)=  \int_0^1\beta\left(\overline{W}_{st}( u)+\overline{h}_{st}( u)\,,\,s+(t-s) u\right)\,d u
\end{equation*}
where for $0\le u\le1$ we defined
\begin{eqnarray*}
  \overline{h}_{st}( u) &=&  u\, h(x,t)+(1- u)h(y,s) \\
  \overline{W}_{st}( u) &=& W\big(s+(t-s) u\big)-\big( u W(t)+(1-  u)
  W(s)\big)
\end{eqnarray*}
Remark that $\overline{W}_{st}( u)$ is now a \emph{rectilinear
Brownian bridge} with $\overline{W}_{st}(0)=\overline{W}_{st}(1)=0$
(see Appendix\myref{bridge}). Then (we just recall the results
without proofs) the transition \pdf\ solution of\refeq{fpt} is
\begin{equation}\label{pdft}
   f(x,t|y,s)=\frac{\EXP{e^{(t-s)Z(s,t)}}}{b(x,t)\sqrt{4\pi
   D(t-s)}}\,e^{-\frac{[h(x,t)-h(y,s)]^2}{4D(t-s)}+\alpha(h(x,t),t)-\alpha(h(y,s),s)}
\end{equation}
When moreover the coefficients are time-independent the \sde\ takes
the form\refeq{sdetind} and the result\refeq{pdft} can be
simplified: we have indeed now
\begin{equation*}
   y=h(x)=\int\frac{dx}{b(x)}\qquad x=g(y) \qquad
    h(g(y))=y\qquad \ha(y) =
   \left[\frac{a(x)}{b(x)}-D\,b'(x)\right]_{x=g(y)}
\end{equation*}
and hence
\begin{equation*}
  \alpha(y) = \frac{1}{2D}\int\ha(y)dy\qquad\quad
  \beta(y) = -\frac{\ha\,^2(y)}{4D}-\frac{\ha\,'(y)}{2}
\end{equation*}
In particular we find
\begin{equation*}
    \alpha\left(h(x)\right)=\frac{1}{2D}\int\left[\frac{a(x)}{b(x)}-D\,b'(x)\right]\frac{dx}{b(x)}=
           \frac{1}{2D}\int\frac{a(x)}{b^2(x)}\,dx-\frac{1}{2}\ln b(x)
\end{equation*}
As a consequence, by redefining now
\begin{equation*}
    Z(s,t)=
    \int_0^1\beta\left(\overline{W}_{st}( u)+\overline{h}( u)\right)\,d u\qquad\quad
    \overline{h}( u) =   u h(x)+(1- u)h(y)
\end{equation*}
the transition \pdf\ solution of\refeq{fp} will become
\begin{eqnarray}\label{pdf}
   f(x,t|y,s)&=&\frac{\EXP{e^{(t-s)Z(s,t)}}}{b(x)\sqrt{4\pi
   D(t-s)}}\sqrt{\frac{b(y)}{b(x)}}\\
   &&\qquad\qquad\exp\left\{\frac{1}{2D}\int_y^x\frac{a(z)}{b^2(z)}\,dz-\frac{1}{4D(t-s)}\left(\int_y^x\frac{dz}{b(z)}\right)^2\right\}\nonumber
\end{eqnarray}
Remark however that the results\refeq{pdft} and\refeq{pdf}
constitute only a semi-explicit form of the \pdf\ because they are
apparently contingent on an expectation not easy to calculate as it
is pointed out in the Section\myref{linsde}

\subsection{Random bridges}\label{bridge}

\subsubsection{Non-random interpolation}\label{intro}

Brownian bridge \sde's are stochastic versions of \ode's for
trajectories interpolating two, or more, fixed points (see for
example~\cite{kars,revuz} pp.\ 358-360). In general the non-random
interpolating trajectories, coinciding with the expectation of the
corresponding random bridges, are supposed to be \emph{linear
functions} of the time $t$, but we will argue here that there is no
really compelling reason for this choice.

Let us start with a trajectory $x(t)$ connecting two possible values
$a$ and $b$ at the endpoints of a compact time interval $[0,T]$,
namely
\begin{equation}\label{nrbridge}
    x(t)=a\,g\!\left(\frac{t}{T}\right)+b\,h\!\left(\frac{t}{T}\right)\quad\qquad
    0\leq t\leq T
\end{equation}
where we will suppose that $g(s)$ and $h(s)$ are derivable at least
once in  $[0,T]$, and
\begin{equation}\label{boundaryc01}
    \left\{
       \begin{array}{l}
         g(0)=1 \\
         g(1)=0
       \end{array}
     \right.
     \qquad\qquad
     \left\{
       \begin{array}{l}
         h(0)=0 \\
         h(1)=1
       \end{array}
     \right.
\end{equation}
so that we trivially get that
\begin{equation}\label{boundaryc}
    x(0)=a\qquad\qquad x(T)=b
\end{equation}
As a matter of fact, every possible function $x(t)$ complying with
the extremal conditions\refeq{boundaryc} can be cast in the
form\refeq{nrbridge}: for a given $x(t)$ just choose an arbitrary
$g(s)$ satisfying\refeq{boundaryc01}, and then it will be enough to
take
\begin{equation*}
    h(s)=\frac{x(Ts)-a\,g(s)}{b}
\end{equation*}
in order to obtain\refeq{nrbridge}. Remark by the way that nothing
forbids an explicit dependence of $g(s)$ and $h(s)$ from $a$ and
$b$, in so far as the conditions\refeq{boundaryc01} are preserved.
Therefore the expression\refeq{nrbridge} can be considered as
general enough for our purposes

We will look then for a first order \ode\ such that the
trajectory\refeq{nrbridge} will be its (unique) solution for the
initial condition $x(0)=a$: it is straightforward to understand that
the form of this equation, albeit independent from the initial
condition $a$, will however explicitly depend on the final condition
$b$ aimed at by our trajectory. A first order \ode\ allows indeed a
free choice just for one initial condition, while in general no
independent final condition can be arbitrarily added if we want to
have a chance to find solutions. As a consequence an \ode\ admitting
both the extremal conditions\refeq{boundaryc} must depend on one of
them: in other words there is no unique equation fitting both the
conditions\refeq{boundaryc} for arbitrary values of $a$ and $b$. In
order to eliminate the initial condition $x(0)=a$ from our \ode\ let
us remark that from\refeq{nrbridge} we get
\begin{equation}\label{nrder}
    \dot{x}(t)=\frac{a}{T}\,\dot{g}\!\left(\frac{t}{T}\right)+\frac{b}{T}\,\dot{h}\!\left(\frac{t}{T}\right)
\end{equation}
so that from\refeq{nrbridge} and\refeq{nrder} we have
\begin{equation*}
    \frac{\dot{x}(t)-\frac{b}{T}\,\dot{h}\!\left(\frac{t}{T}\right)}{\frac{1}{T}\,\dot{g}\!\left(\frac{t}{T}\right)}=a=
    \frac{x(t)-b\,h\!\left(\frac{t}{T}\right)}{g\!\left(\frac{t}{T}\right)}
\end{equation*}
and rearranging the terms we find the \ode
\begin{equation}\label{nreq}
    \dot{x}(t)=\frac{b\left[g\!\left(\frac{t}{T}\right)\dot{h}\!\left(\frac{t}{T}\right)-\dot{g}\!\left(\frac{t}{T}\right)h\!\left(\frac{t}{T}\right)\right]+x(t)\dot{g}\!\left(\frac{t}{T}\right)}{T\,g\!\left(\frac{t}{T}\right)}
\end{equation}
whose solutions\refeq{nrbridge} will connect every possible initial
condition $a$ to the same final value $b$ inscribed into it.

The simplest example, used for the most widespread stochastic
generalization to Brownian bridges (see~\cite{kars,revuz} for
details), adopts the following linear functions
\begin{equation}\label{lin}
    g(s)=1-s\qquad\quad h(s)=s
\end{equation}
so that the corresponding connecting trajectories
\begin{equation*}
    x(t)=a\left(1-\frac{t}{T}\right)+b\,\frac{t}{T}
\end{equation*}
are $t$-linear and satisfy the \ode
\begin{equation*}
    \dot{x}(t)=\frac{b-x(t)}{T-t}
\end{equation*}
But this by no means constitutes a unique possibility. We could for
instance take
\begin{equation*}
    g(s)=(1-s)^2\qquad\quad h(s)=s^2
\end{equation*}
and in this case we would have the parabolic trajectories
\begin{equation*}
    x(t)=a\left(1-\frac{t}{T}\right)^2+b\left(\frac{t}{T}\right)^2
\end{equation*}
solutions of the \ode
\begin{equation*}
    \dot{x}(t)=\frac{2}{T}\,\frac{b\,t-Tx(t)}{T-t}
\end{equation*}
As final examples among many we can either put
\begin{equation*}
    g(s)=\cos\frac{\pi s}{2}\qquad\quad h(s)=\sin\frac{\pi s}{2}
\end{equation*}
and get as trajectories
\begin{equation*}
    x(t)=a\,\cos\frac{\pi t}{2T}+b\,\sin\frac{\pi t}{2T}
\end{equation*}
and as \ode
\begin{equation*}
    \dot{x}(t)=\frac{\pi}{2T}\,\frac{b-x(t)\,\sin\frac{\pi t}{2T}}{\cos\frac{\pi t}{2T}}
\end{equation*}
or instead
\begin{equation*}
    g(s)=\cos^2\frac{\pi s}{2}\qquad\quad h(s)=\sin^2\frac{\pi s}{2}
\end{equation*}
and get as trajectories
\begin{equation*}
    x(t)=a\,\cos^2\frac{\pi t}{2T}+b\,\sin^2\frac{\pi t}{2T}
\end{equation*}
and as \ode
\begin{equation*}
    \dot{x}(t)=\frac{b-x(t)}{T}\,\pi\tan\frac{\pi t}{2T}
\end{equation*}

\subsubsection{Brownian bridges}

A straightforward stochastic generalization of the \ode\refeq{nreq}
discussed in the previous section is obtained just by adding a
Brownian noise $W(t)$ with constant diffusion coefficient $2D$
\begin{equation}\label{brbreq}
    dX(t)=\frac{b\left[g\!\left(\frac{t}{T}\right)\dot{h}\!\left(\frac{t}{T}\right)-\dot{g}\!\left(\frac{t}{T}\right)h\!\left(\frac{t}{T}\right)\right]
    +X(t)\dot{g}\!\left(\frac{t}{T}\right)}{T\,g\!\left(\frac{t}{T}\right)}\,dt+dW(t)
\end{equation}
and apparently results in a \sde\ with linear coefficients and
initial condition $X(0)=a,\,\Pqo$ Therefore, according to the
notations adopted in the Appendix\myref{elem}, the
equation\refeq{brbreq} is a \sde\ of the form\refeq{linear} with
\begin{eqnarray*}
  \ha_0(t) &=& \frac{b}{T}\,\frac{g\!\left(\frac{t}{T}\right)\dot{h}\!\left(\frac{t}{T}\right)-\dot{g}\!\left(\frac{t}{T}\right)h\!\left(\frac{t}{T}\right)}
  {g\!\left(\frac{t}{T}\right)}=b\,g\!\left(\frac{t}{T}\right)\frac{d}{dt}\left[\frac{h\!\left(\frac{t}{T}\right)}{g\!\left(\frac{t}{T}\right)}\right] \\
  \ha_1(t) &=& \frac{1}{T}\,\frac{\dot{g}\!\left(\frac{t}{T}\right)}{g\!\left(\frac{t}{T}\right)}
  =\frac{d}{dt}\left[\ln g\!\left(\frac{t}{T}\right)\right]\\
  \hb_0(0) &=& 1  \\
  \hb_1(t) &=& 0
\end{eqnarray*}
and hence it is easy to find that
\begin{equation*}
    \overline{Z}(t)=\int\ha_1(t)\,dt=\ln
    g\!\left(\frac{t}{T}\right)\qquad\quad e^{\overline{Z}(t)}=g\!\left(\frac{t}{T}\right)
\end{equation*}
so that its solution after some algebra becomes
\begin{equation}\label{brbr}
    X(t)=a\,g\!\left(\frac{t}{T}\right)+b\,h\!\left(\frac{t}{T}\right)+g\!\left(\frac{t}{T}\right)\int_0^t\frac{1}{g\!\left(\frac{s}{T}\right)}\,dW(s)
\end{equation}
Remark that the first two terms of this solution exactly coincide
with the non-random interpolating trajectories\refeq{nrbridge}. In
particular when we take the linear functions\refeq{lin} we get the
solution
\begin{equation}\label{rectibrbr}
    X(t)=a\left(1-\frac{t}{T}\right)+b\,\frac{t}{T}+(T-t)\int_0^t\frac{1}{T-s}\,dW(s)
\end{equation}
which is exactly the \emph{rectilinear Brownian bridge} discussed
in~\cite{kars}

Following the same line of reasoning of~\cite{kars} we could now
prove that the solution\refeq{brbr} is Gaussian with $\Pqo$
continuous paths, with expectation
\begin{equation*}
    m(t)=\EXP{X(t)}=a\,g\!\left(\frac{t}{T}\right)+b\,h\!\left(\frac{t}{T}\right)
\end{equation*}
and with covariance
\begin{equation*}
    C(s,t)=\cov{X(s)}{X(t)}=D\,g\bigg(\frac{s}{T}\bigg)g\!\left(\frac{t}{T}\right)\int_0^{s\wedge
    t}\frac{1}{g^2\!\left(\frac{u}{T}\right)}\,du
\end{equation*}
so that its laws can now be deemed completely known. In particular
for the variance we have
\begin{equation*}
    \VAR{X(t)}=C(t,t)=D\,g^2\!\left(\frac{t}{T}\right)\int_0^t\frac{1}{g^2\!\left(\frac{s}{T}\right)}\,ds
\end{equation*}
In the case of the rectilinear Brownian bridge\refeq{rectibrbr}
these formulas give
\begin{eqnarray*}
  m(t) &=& a\left(1-\frac{t}{T}\right)+b\,\frac{t}{T} \\
  C(s,t) &=& D(T-s)(T-t)\int_0^{s\wedge
  t}\frac{du}{(T-u)^2}=D\left((s\wedge t)-\frac{st}{T}\right)
\end{eqnarray*}
and its distributions coincide with that of a Wiener process
conditioned at both the endpoints with $X(0)=a$ and $X(T)=b$: in
fact, if we take
\begin{equation*}
    \phi(x,t|y)=\frac{e^{-\frac{(x-y)^2}{2Dt}}}{\sqrt{2\pi Dt}}
\end{equation*}
the finite dimensional distributions of the rectilinear Brownian
bridge coincide with the following conditional \pdf's of a Wiener
process for $0=t_0<t_1<\ldots<t_n<T$
\begin{eqnarray*}
  f(x_1,t_1;\ldots;x_n,t_n\,|\,a,0;b,T) &=& \frac{f(x_1,t_1;\ldots;x_n,t_n;b,T\,|\,a,0)}{f(b,T\,|\,a,0)} \\
   &=& \frac{f(b,T|x_n,t_n)\ldots
   f(x_2,t_2|x_1,t_1)f(x_1,t_1|a,0)}{f(b,T\,|\,a,0)}\\
   &=&\frac{\phi(b,T-t_n|x_n)\ldots\phi(x_2,t_2-t_1|x_1)\phi(x_1,t_1|a)}{\phi(b,T|a)}
\end{eqnarray*}
In fact we will call \emph{rectilinear Brownian bridge} every
stochastic process with such finite dimensional distributions. In
particular it could be shown that
\begin{equation*}
    a\left(1-\frac{t}{T}\right)+b\,\frac{t}{T}+\left(W(t)-\frac{t}{T}W(T)\right)
\end{equation*}
also is a rectilinear Brownian bridge

\subsection{Eigenfunction expansion}\label{eigenexp}

For a \fpe\refeq{fpt} with both $a(x)$ and $b(x)$ time-independent
coefficients we will here replace $a(x)$ with $\fvel(x)$ (the symbol
that we will adopt in the Section\myref{stochmech} for the forward
velocity within the framework of the stochastic mechanics), while
for short we will adopt the notation $B(x)=D\,b^2(x)$. The \pdf\
$f(x,t)$ of our continuous Markov process $X(t)$ is then a solution
of the \fpe
\begin{equation}\label{fpeq}
    \partial_tf=\partial_x^2(B f)-\partial_x(\fvel f)=
     \partial_x\bigl[\partial_x(B f)-\fvel f\bigr]
\end{equation}
defined for $x\in[a,b]$ (beware also the temporary new meaning of
the symbols $a$ and $b$) and $t\geq s$. We will further suppose that
$\fvel(x)$ has no singularities in $(a,b)$, and that $\fvel(x)$ and
$B(x)$ are both continuous and differentiable functions. The
conditions imposed on the probabilistic solutions are of course
\begin{equation}\label{cond}
    f(x,t)\geq0\qquad
     \int_a^bf(x,t)\,dx=1\qquad\quad a<x<b\qquad s\leq t
\end{equation}
while from\refeq{fpeq} the second condition\refeq{cond} also takes
the form
\begin{equation*}
    \bigl[\partial_x(Bf)-\fvel f\bigr]_a^b=0\,,\qquad s\leq t
\end{equation*}
Suitable initial conditions will be added to produce the required
evolution: for example a transition pdf $f(x,t|y,s)$ will be
selected by the initial condition
\begin{equation}\label{transcond}
   \lim_{t\to s^+}f(x,t)=f(x,s^+)=\delta(x-y)
\end{equation}
It is also possible to show by direct calculation that
\begin{equation}\label{ssol}
    \widetilde{f}(x)=Z^{-1}\,e^{-\int\frac{B'(x)\,-\fvel(x)}{B(x)}\,dx}\quad\qquad
       Z=\int_a^b e^{-\int\frac{B'(x)\,-\fvel(x)}{B(x)}\,dx}\,dx
\end{equation}
is an invariant solution of\refeq{fpeq} satisfying the
conditions\refeq{cond} (for its coherence with the Boltzmann
distribution\refeq{boltz} see Appendix\myref{boltzdistr}). Remark
that\refeq{fpeq} is not in the standard self-adjoint form, but if we
define the new function $g(x,t)$ by means of
\begin{equation*}
    f(x,t)=\sqrt{\widetilde{f}(x)}\,g(x,t)
\end{equation*}
it would be easy to show\mycite{pla,jpha} that $g(x,t)$ obeys now an
equation of the form
\begin{equation*}
    \partial_t\,g=\Lop{g}
\end{equation*}
where the operator ${\cal L}$ defined on a test function $\varphi$
as
\begin{align*}
  & \qquad\qquad\qquad \Lop{\varphi}=
     {d\over dx}\left[p(x)\,{d\varphi(x)\over
     dx}\right]-q(x)\varphi(x)\\
  & p(x) = B(x)>0 \qquad\quad
  q(x) = {\bigl[B'(x)-\fvel(x)\bigr]^2\over4B(x)}\,-\,
                   {\bigl[B'(x)-\fvel(x)\bigr]'\over2}
\end{align*}
is now self-adjoint. Then, by separating the variables by means of
$g(x,t)=\gamma(t)G(x)$ we have $\gamma(t)={\rm e}^{-\lambda t}$
while $G$ must be solution of a typical Sturm-Liouville problem
associated to the equation
\begin{equation}\label{slpr}
    \Lop{G(x)}+\lambda G(x)=0
\end{equation}
with the boundary conditions
\begin{eqnarray}
\bigl[B'(a)-\fvel(a)\bigr]G(a)+2B(a)G'(a) &=& 0\label{bcond} \\
  \bigl[B'(b)-\fvel(b)\bigr]G(b)+2B(b)G'(b) &=& 0\label{bcond1}
\end{eqnarray}
It easy to see that $\lambda=0$ is always an eigenvalue for the
problem\refeq{slpr} with\refeq{bcond} and\refeq{bcond1}, and that
the corresponding eigenfunction is $G_0(x)=\sqrt{\widetilde{f}(x)}$.
The other simple eigenvalues $\lambda_n$ will then constitute an
infinite, increasing sequence and the corresponding eigenfunction
$G_n(x)$ will have $n$ simple zeros in $(a,b)$. This also means that
$\lambda_0=0$, corresponding to the eigenfunction $G_0(x)$ which
never vanishes in $(a,b)$, is the lowest eigenvalue so that all
other eigenvalues are strictly positive. The eigenfunctions will
constitute a complete orthonormal set of functions in
$L^2\bigl([a,b]\bigr)$ so that the general solution of\refeq{fpeq}
with\refeq{cond} will have the form
\begin{equation}\label{gensol}
    f(x,t)=\sum_{n=0}^{\infty}c_n e^{-\lambda_nt}\sqrt{\widetilde{f}(x)}G_n(x)
    =\sum_{n=0}^{\infty}c_n e^{-\lambda_nt}G_0(x)G_n(x)
\end{equation}
with $c_0=1$ for normalization (remember that $\lambda_0=0$). The
coefficients $c_n$ for a particular solution selected by an initial
condition
\begin{equation*}
   f(x,s^+)=f_0(x)
\end{equation*}
are then calculated from the orthonormality relations as
\begin{equation*}
    c_n=\int_a^bf_0(x)\,{G_n(x)\over G_0(x)}\,dx
\end{equation*}
and in particular for the transition \pdf\ we have
from\refeq{transcond} that
\begin{equation*}
    c_n={G_n(x_0)\over G_0(x_0)}
\end{equation*}
Since $\lambda_0=0$ and $\lambda_n>0$ for $n\geq1$, the general
solution\refeq{gensol} has a precise time evolution: all the
exponential factors vanish with $t\to+\infty$ with the only
exception of the term $n=0$ which is constant, so that exponentially
fast we will always have
\begin{equation*}
    \lim_{t\to+\infty}f(x,t)=c_0\sqrt{\widetilde{f}(x)}\,G_0(x)=\widetilde{f}(x)
\end{equation*}
namely the general solution will relax toward the invariant solution
$\widetilde{f}(x)$

\subsection{Invariant Boltzmann distributions}\label{boltzdistr}

Let us consider the \sde\ with time-independent coefficients
\begin{equation*}
    dX(t)=a(X(t))\,dt+b(X(t))\,dW(t)
\end{equation*}
and the corresponding \fpe\ (see also Section\myref{thefpe})
\begin{eqnarray}
    \partial_tf_X(x,t)&=&-\partial_x\left[a(x)f_X(x,t)\right]+D\,\partial^2_x\left[b^2(x)f_X(x,t)\right]\nonumber\\
    &=&-\partial_x\left[a(x)f_X(x,t)-D\,\partial_x\left(b^2(x)f_X(x,t)\right)\right]\label{statfpx}
\end{eqnarray}
We know then from the Section\myref{constcoeff} that the
transformation $Y(t)=h(X(t))$ with a monotonic $h(x)$ and
\begin{equation}\label{transf}
    y=h(x)=\int\frac{dx}{b(x)}\qquad h'(x)=\frac{1}{b(x)}\qquad
    x=g(y)=h^{-1}(x)
\end{equation}
will bring us to a new \sde
\begin{equation*}
    dY(t)=\ha(Y(t))\,dt+dW(t)
\end{equation*}
and to a new \fpe
\begin{equation}
    \partial_tf_Y(y,t)=-\partial_y\left[\ha(y)f_Y(y,t)\right]+D\,\partial^2_yf_Y(y,t)
    =-\partial_y\left[\ha(y)f_Y(y,t)-D\,\partial_yf_Y(y,t)\right]\label{statfpy}
\end{equation}
where now
\begin{equation*}
    \ha(y)=\left[\frac{a(x)}{b(x)}-Db'(x)\right]_{x=g(y)}
\end{equation*}
We claimed in the Section\myref{elem} that if we introduce a
potential $\phi(y)$ according to
\begin{equation*}
    \ha (y)=-\frac{D}{kT}\,\phi'(y)
\end{equation*}
then -- provided that $\phi(y)$ is such that
$e^{-\frac{\phi(y)}{kT}}$ is an integrable function -- it is
possible to show that
\begin{equation}\label{boltz1}
   \widetilde{f}_Y(y)=Z^{-1}e^{-\frac{\phi(y)}{kT}}= Z^{-1}e^{\frac{1}{D}\int\ha(y)\,dy}\qquad\quad Z=\int_\R e^{-\frac{\phi(z)}{kT}}dz
\end{equation}
is an invariant Boltzmann distribution. It is apparent indeed that,
since it is
\begin{equation*}
    \partial_y\widetilde{f}_Y(y)=\frac{\ha(y)}{D}\,\widetilde{f}_Y(y)
\end{equation*}
then $\widetilde{f}_Y(y)$ is a stationary solution
of\refeq{statfpy}. On the other hand, when $Y(t)=h(X(t))$ according
to the transformation\refeq{transf}, the respective \pdf's are also
connected by the transformations
\begin{equation*}
    f_Y(h(x),t)=|b(x)|f_X(x,t)\qquad\quad
    f_Y(y,t)=\big[|b(x)|f_X(x,t)\big]_{x=g(y)}
\end{equation*}
and hence in particular we have
\begin{equation}\label{transtat}
    \widetilde{f}_X(x)=\frac{\widetilde{f}_Y(h(x))}{|b(x)|}
\end{equation}
that we claim to be a stationary solution for the original
\fpe\refeq{statfpx}

In order to check our last statement, we will first express
$\widetilde{f}_X(x)$ in terms of the coefficients $a(x)$ and $b(x)$
of\refeq{statfpx}: to this end let us remark that
\begin{equation*}
    \int\ha(y)\,dy=\int\left[\frac{a(x)}{b(x)}-Db'(x)\right]_{x=g(y)}dy=\left[\int\left(\frac{a(x)}{b(x)}-Db'(x)\right)\frac{dx}{b(x)}\right]_{x=g(y)}
\end{equation*}
because, taking for short
\begin{equation*}
    A(x)=\frac{a(x)}{b(x)}-Db'(x)
\end{equation*}
we have\footnote{Given a function $f(x)$, its primitives are the
functions $F(x)$ such that
\begin{equation*}
    F(x)=\int f(x)\,dx\qquad\quad F'(x)=f(x)
\end{equation*}
Let us suppose now to have an invertible transformation of variables
$y=h(x)$ with
\begin{equation*}
    y=h(x)\qquad\quad x=h^{-1}(y)=g(y)\qquad\quad g(h(x))=x\qquad\quad g'(y)=\frac{1}{h'\left(g(y)\right)}
\end{equation*}
It is then easy to show that the function
\begin{equation*}
    F_1(x)=\left[\int f(g(y))g'(y)\,dy\right]_{y=h(x)}=G(h(x))\qquad\qquad
    G(y)=\int f(g(y))g'(y)\,dy
\end{equation*}
is again a primitive of $f(x)$ because
\begin{equation*}
    F'_1(x)=G'(h(x))h'(x)=h'(x)\left[f(g(y))g'(y)\right]_{y=h(x)}=h'(x)\left[\frac{f(g(y))}{h'\left(g(y)\right)}\right]_{y=h(x)}=f(x)
\end{equation*}
As a consequence the rules for the change of variables in the
indefinite integrals are
\begin{equation*}
  \int f(x)\,dx=\left[\int f(g(y))\,g'(y)\,dy\right]_{y=h(x)}\qquad\quad\left[\int f(x)\,dx\right]_{x=g(y)}=\int f(g(y))\,g'(y)\,dy
\end{equation*}
To have the formula in the text just take $f(x)=A(x)h'(x)$
 }
for the change of
variable $y=h(x)$
\begin{eqnarray*}
    \int\ha(y)\,dy&=&\int A(g(y))\,dy=\left[\int
    A(x)h'(x)\,dx\right]_{x=g(y)}\\
    &=&\left[\int\left(\frac{a(x)}{b(x)}-Db'(x)\right)\frac{dx}{b(x)}\right]_{x=g(y)}
\end{eqnarray*}
or also the equivalent formulation
\begin{equation*}
   \left[\int\ha(y)\,dy\right]_{y=h(x)}=\int\left(\frac{a(x)}{b(x)}-Db'(x)\right)\frac{dx}{b(x)}
\end{equation*}
so that from\refeq{boltz1} and\refeq{transtat} we find
\begin{equation*}
    \widetilde{f}_X(x)=\frac{e^{\frac{1}{D}\int\left(\frac{a(x)}{b(x)}-D\,b'(x)\right)\frac{1}{b(x)}\,dx}}{Z|b(x)|}
\end{equation*}
On the other hand, since it is easy to see that
\begin{equation*}
    e^{\,\int\frac{b'(x)}{b(x)}\,dx}=e^{\,\int\frac{d}{dx}|\ln
    b(x)|\,dx}=e^{\,\ln|b(x)|}=|b(x)|
\end{equation*}
we finally have
\begin{equation}\label{boltz2}
    \widetilde{f}_X(x)=Z^{-1}e^{\frac{1}{D}\int\left(\frac{a(x)}{b(x)}-2D\,b'(x)\right)\frac{1}{b(x)}\,dx}
\end{equation}
Now it is easy to check that this $\widetilde{f}_X(x)$ is the
invariant solution of the \fpe\refeq{statfpx}: we find indeed that
\begin{eqnarray*}
    D\partial_x\left[b^2(x)\widetilde{f}_X(x)\right]&=&\widetilde{f}_X(x)\left[2Db(x)b'(x)+b^2(x)\left(\frac{a(x)}{b(x)}-2D\,b'(x)\right)\frac{1}{b(x)}\right]\\
    &=&a(x)\widetilde{f}_X(x)
\end{eqnarray*}
It is finally apparent that the \pdf\refeq{boltz2} also coincides
with the invariant \pdf\refeq{ssol} because within the notation of
the Section\myref{eigenexp} we have
\begin{equation*}
    -\frac{B'(x)-\fvel(x)}{B(x)}=\left(\frac{a(x)}{b(x)}-2D\,b'(x)\right)\frac{1}{b(x)}
\end{equation*}
By summarizing we can say that, if the \pdf\refeq{boltz1} is the
invariant solution (if it exists) of the stationary
\fpe\refeq{statfpy} for the process $Y(t)=h(X(t))$, then the
\pdf\refeq{boltz2} is the corresponding invariant solution of the
\fpe\refeq{statfpx} for the process $X(t)$

\section{An \emph{anamnesis} of stochastic mechanics}\label{stochmech}

Initially proposed as a possible interpretation of quantum
mechanics\mycite{nels,guerra} with the challenging aim of shedding
new light on its enduring mysteries, over the years the stochastic
mechanics evolved into a more general theory dealing with
conservative diffusion processes\mycite{gm}. Its tools are therefore
valuable nowadays even beyond the strict quantum precinct, and in
particular they have been employed in the broad field of the
stochastic control\mycite{control}. We will refrain however from
giving here a comprehensive review of these topics referring the
interested readers to the quoted literature, and we will rather
confine ourselves in the present appendix to recall just the few
results deemed to be instrumental for the follow-up of the present
enquiry

From a quantum wave function $\psi(x,t)$ solution of a
(one-dimensional) Schr\"odinger equation
\begin{equation}\label{se}
    i\hbar\partial_t\psi(x,t)=-\,{\hbar^2\over 2m}\,\partial_x^2\psi(x,t)+V(x,t)\psi(x,t)
\end{equation}
we can deduce the form of the forward velocity $\fvel(x,t)$ -- here
the upper arrow just means \emph{forward}, and does not denote a
\emph{vector}, while $\bvel(x,t)$ will be understood as a
\emph{backward} velocity -- appearing in both the \fpe
\begin{eqnarray}\label{fpe}
    \partial_tf(x,t)&=&-\partial_x[\,\fvel(x,t)f(x,t)]+D\,\partial_x^2f(x,t)\nonumber\\
    &=&\partial_x\left[D\partial_xf(x,t)-\fvel(x,t)f(x,t)\right]
\end{eqnarray}
for the \pdf\ $f(x,t)=|\psi(x,t)|^2$, and the associated \ito\ \sde\
\begin{equation}\label{sde}
    dX(t)=\fvel(X(t),t)dt+dW(t)
\end{equation}
for the corresponding Markov process $X(t)$ in the framework of the
stochastic mechanics: here $W(t)$ is a Wiener process with a
constant diffusion coefficient $2D=\frac{\hbar}{m}$, namely such
that $\EXP{W(t)^2}=2D\,t$. If $\psi(x,t)$ is an arbitrary solution
of\refeq{se}, it is well known indeed that with the usual {\it
Ansatz}
\begin{equation}\label{ansatz}
    \psi(x,t)=R(x,t)\,{\rm e}^{iS(x,t)/\hbar}
\end{equation}
where $R$ and $S$ are real functions, $R^2=|\psi|^2$ comes out to be
a particular solution of the \fpe\refeq{fpe} with forward velocity
field of the form
\begin{equation}\label{fvel}
    \fvel(x,t)={\partial_xS(x,t)\over m}\,+{\hbar\over 2m}\,\partial_x\left[\ln
    R^2(x,t)\right]
\end{equation}
as it is deduced by separating the real and imaginary parts
of\refeq{se}. Remark that the explicit dependence of $\fvel$ on  the
form of $R$ clearly indicates that to have a solution of\refeq{fpe}
which makes quantum sense we must pick-up just one, suitable
solution. In the stochastic mechanical framework, indeed, the system
is ruled not only by the \fpe\refeq{fpe}, but also by a second,
dynamical equation (the imaginary part)
\begin{equation}\label{hjme}
    \partial_tS(x,t)+{\left[\partial_xS(x,t)\right]^2\over 2m}+V(x,t)\,-\,{\hbar^2\over2m}\,
                    {\partial^2_xR(x,t)\over R(x,t)}=0
\end{equation}
known as Hamilton-Jacobi-Madelung equation

Let us consider now the Schr\"odinger equation\refeq{se} in the case
of a time-independent potential $V(x)$, with a Hamiltonian
\begin{equation*}
    \widehat H\psi(x)=-\,{\hbar^2\over
    2m}\,\partial_x^2\psi(x)+V(x)\psi(x)
\end{equation*}
with purely discrete spectrum and stationary, normalizable states,
and let us use the following notations for these states, and their
eigenvalues and eigenfunctions:\begin{eqnarray*}
  \psi_n(x,t) &=& \phi_n(x)\,{\rm e}^{-iE_nt/\hbar} \\
  \widehat H\phi_n(x) &=& -{\hbar^2\over
  2m}\,\phi''_n(x)+V(x)\phi_n(x)=E_n\phi_n(x)
\end{eqnarray*}
For later convenience we will also introduce the constant
\begin{equation}\label{diffus}
    D={\hbar\over 2m}
\end{equation}
so that the previous eigenvalue equation can be recast in the
following form
\begin{equation*}
    D\phi''_n(x)={V(x)-E_n\over\hbar}\,\phi_n(x)
\end{equation*}
For a stationary solution $\psi_n(x,t)$ the
\emph{Ansatz}\refeq{ansatz} will give
\begin{equation*}
    S(x,t)=-E_nt\,,\qquad R(x,t)=\phi_n(x)
\end{equation*}
so that for our stationary states the velocity fields are
\begin{equation}\label{statfvel}
   \fvel_n(x)=2D\,{\phi'_n(x)\over\phi_n(x)}
\end{equation}

\subsection{The \fpe\ for stationary states}

From\refeq{statfvel} we see that the forward velocities $\fvel_n(x)$
for stationary states are time-independent, and that they have
singularities in the zeros (nodes) of the eigenfunction. Since the
$n$-th eigenfunction of a quantum system with bound states has
exactly $n$ simple nodes $x_1,\dots,x_n$, the coefficients of the
\fpe
\begin{equation}\label{statfpe}
    \partial_tf(x,t)=-\partial_x[\,\fvel_n(x)f(x,t)]+D\,\partial_x^2f(x,t)=\partial_x\left[D\partial_xf(x,t)-\fvel_n(x)f(x,t)\right]
\end{equation}
diverge in these $n$ points and we will be obliged to solve it in
separate intervals by imposing suitable boundary conditions
connecting the different sections (see Appendix\myref{thefpe} for
further details). As a matter of fact, these singularities
effectively separate the real axis in $n+1$ sub-intervals with
impenetrable (to the probability current) walls. Hence the process
will not have an \emph{unique} invariant measure and will never
cross the boundaries fixed by the singularities of $\fvel_n(x)$: if
we start at $t_0$ in one of the intervals in which the axis is
divided we will always remain therein. As a consequence, with an
arbitrary initial distribution, we must require that the integrals
\begin{equation*}
    \int_{x_k}^{x_{k+1}}f(x,t)\,dx
\end{equation*}
be kept at a constant value for $t\geq t_0$: this values will not,
in general, be equal to $1$ (only their sum will amount to  $1$)
and, since the separate intervals can not communicate, they will be
fixed by the choice of the initial conditions. The boundary
conditions are hence imposed by the conservation of the probability
in $[x_k,x_{k+1}]$ and that entails the vanishing of the probability
current in\refeq{statfpe} at the end points of the intervals:
\begin{equation*}
    \bigl[D\partial_xf(x,t)-\fvel_n(x)f(x,t)\bigr]_{x_k,x_{k+1}}=0\,,\qquad t\geq t_0
\end{equation*}

Therefore, since every particular solution is selected by the
initial conditions, we are first interested in finding the
transition \pdf\ $f(x,t|y,s)$ which is singled out by the
condition\begin{equation}\label{degenincond}
   \lim_{t\to s^+}f(x,t)=f(x,s^+)=\delta(x-y)
\end{equation}
The evolution of any other initial condition $f_0(x)$ at $t=s$ is
subsequently ruled by the Chapman-Kolmogorov equation
\begin{equation}\label{chk}
    f(x,t)=\int_{-\infty}^{+\infty}f(x,t|y,s)f_0(y)\,dy
\end{equation}
To solve\refeq{statfpe} in every interval $[x_k,x_{k+1}]$ (both
finite or infinite), when we already know the invariant,
time-independent solution $\phi_n^2(x)$, we usually put
\begin{equation*}
    f(x,t)=\phi_n(x)g(x,t)
\end{equation*}
in order to reduce\refeq{statfpe} to the form
\begin{equation}\label{sadjform}
    \partial_tg={\cal L}_ng
\end{equation}
where ${\cal L}_n$ is now the self-adjoint operator defined on
$[x_k,x_{k+1}]$ as
\begin{align*}
  & \quad{\cal L}_ng(x)={{\rm d}\over{\rm d}x}
        \left[p(x){{\rm d}g(x)\over{\rm d}x}\right]
                     -q_n(x)g(x)  \\
  &  p(x)=D>0\quad\qquad q_n(x)={\fvel_n^{\,2}(x)\over 4D}+{\fvel_n'(x)\over 2}
\end{align*}
To solve\refeq{sadjform} it is then advisable to separate the
variables into $g(x,t)=\gamma(t)G(x)$, so that we immediately have
$\gamma(t)={\rm e}^{-\lambda t}$, while $G(x)$ must be a solution of
the Sturm-Liouville eigenvalue problem associated to the equation
\begin{equation}\label{sl}
    {\cal L}_nG(x)+\lambda G(x)=0
\end{equation}
with the boundary conditions
\begin{equation}\label{slcond}
    \bigl[2DG'(x)-\fvel_n(x)G(x)\bigr]_{x_k,\,x_{k+1}}=0
\end{equation}
The general behavior of the solutions obtained as expansions in the
system of the eigenfunctions of\refeq{sl} has already been discussed
elsewhere\mycite{pla,jpha}

\subsection{The \ito\ \sde\ for stationary states}

If on the other hand we would like to see the problem from the
standpoint of the \ito\ \sde's for some Markov process $X(t)$ we
will be confronted with the Smoluchowsky equations of the
type\refeq{sdeb1}
\begin{equation}\label{statsde}
    dX(t)=\fvel_n(X(t))dt+dW(t)
\end{equation}
which are the path-wise counterparts of the \fpe's\refeq{statfpe}.
It is well known from our discussion inthe
Appendix\myref{boltzdistr} that in this case, if
\begin{equation*}
    e^{\frac{1}{D}\int\fvel_n(x)dx}
\end{equation*}
is an integrable function, then
\begin{equation}\label{statsol}
    \frac{e^{\frac{1}{D}\int\fvel_n(x)dx}}{\int_{\RE} e^{\frac{1}{D}\int\fvel_n(x)dx}dx}
\end{equation}
is the \pdf\ of the stationary solution of\refeq{statsde}. On the
other hand we also know that there is no simple way to find the
general solutions of\refeq{statsde} without scrutinizing it in
particular cases

\subsection{\qho\ stationary states}

Let us then consider in detail the particular example of a \qho\
with the potential
\begin{equation*}
    V(x)={m\over2}\,\omega^2x^2
\end{equation*}
It is well-known that its eigenvalues are
\begin{equation*}
    E_n=\hbar\omega\left(n+{1\over2}\right)\,;\qquad n=0,1,2\dots
\end{equation*}
while, with the notation
\begin{equation}\label{diffusvar}
    \sigma^2={\hbar\over2m\omega}=\frac{D}{\omega}
\end{equation}
the eigenfuncions are
\begin{equation*}
    \phi_n(x)=
    {1\over\sqrt{\sigma\sqrt{2\pi}2^nn!}}\,{\rm e}^{-x^2/4\sigma^2}\,
                 H_n\left({x\over\sigma\sqrt{2}}\right)
\end{equation*}
where $H_n$ are the Hermite polynomials (see\mycite{grad} 8.95). The
corresponding forward velocity fields are then easily calculated
from\refeq{statfvel}
\begin{equation}\label{qhofvel}
    \fvel_n(x)=\omega\sigma\sqrt{2}\,\frac{H'_n\left({x\over\sigma\sqrt{2}}\right)}{H_n\left({x\over\sigma\sqrt{2}}\right)}-\omega x
    =\omega\sigma\sqrt{2}\left[\frac{2nH_{n-1}(z)}{H_n(z)}-z\right]_{z={x\over\sigma\sqrt{2}}}
\end{equation}
and the first examples are
\begin{figure}
\begin{center}
\includegraphics*[width=14cm]{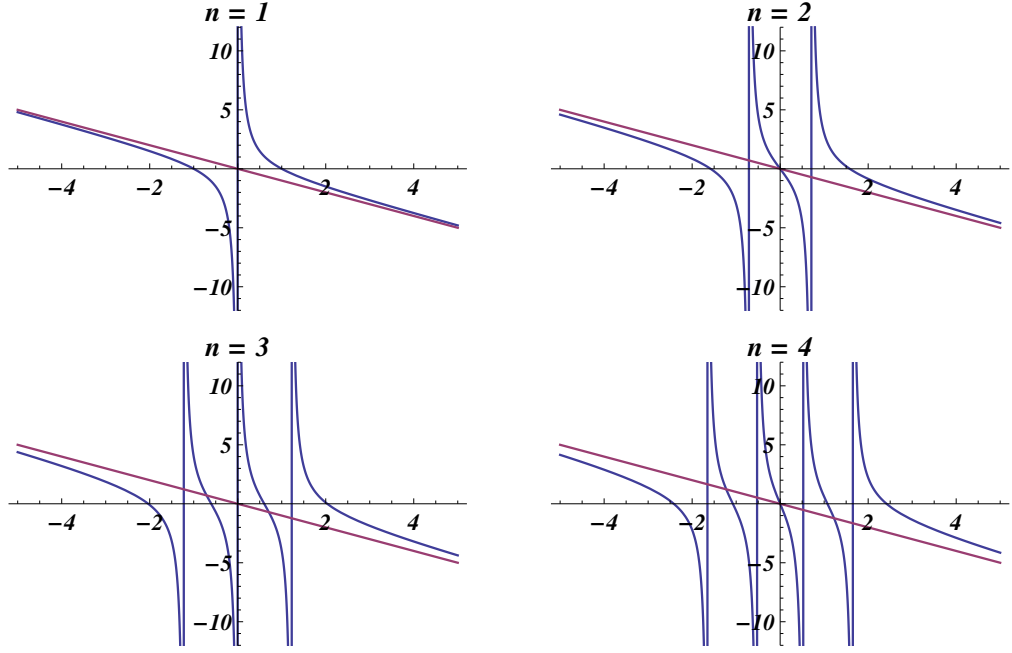}
\caption{The dimensionless forward velocities
$\frac{\fvel_n}{\omega\sigma\sqrt{2}}$ for $n=1,2,3,4$ plotted in
function of the dimensionless variable $\frac{x}{\sigma\sqrt{2}}$.
The case $n=0$ leads to the straight line}\label{hermfig}
\end{center}
\end{figure}
\begin{eqnarray*}
  \fvel_0(x) &=& -\omega x \\
  \fvel_1(x) &=& -\omega x+\omega\sigma\,{2\sigma\over x} \\
  \fvel_2(x) &=& -\omega x+\omega\sigma\,{4\sigma x\over x^2-\sigma^2} \\
  \fvel_3(x) &=& -\omega x+\omega\sigma\,\frac{6\sigma(x^2-\sigma^2)}{x(x^2-3\sigma^2)}\\
  \fvel_4(x) &=& -\omega x+\omega\sigma\,\frac{8\sigma x(x^2-3\sigma^2)}{x^4-6\sigma^2x^2+3\sigma^4}
\end{eqnarray*}
with singularities in the zeros $x_k$ of the Hermite polynomials.
These velocities are \emph{piecewise} (between two subsequent
singularities) \emph{monotonic decreasing} functions as displayed in
the Figure\myref{hermfig}

\subsection{The processes for the \qho\ stationary states}

We first recall a few general remarks about the solution methods of
the eigenvalue problem\refeq{sl} which for our forward velocities
$\fvel_n(x)$, with $\epsilon=\hbar\lambda$, can be written as
\begin{equation*}
    -\,{\hbar^2\over2m}\,G''(x)+
        \left({m\over2}\,\omega^2x^2-\hbar\omega\,{2n+1\over2}\right)G(x)
                    =\epsilon\, G(x)
\end{equation*}
in every interval $[x_k,x_{k+1}],\;k=0,1,\dots,n$ between two
subsequent singularities of $\fvel_n(x)$, with the boundary
conditions
\begin{equation*}
   \left[\phi_n(x)G'(x)-\phi'_n(x)G(x)\right]_{x_k,x_{k+1}}=0
\end{equation*}
Since $\phi_n(x)$ (but not $\phi'_n(x)$) vanishes in the $x_k$'s,
the actual boundary conditions to impose are
\begin{equation*}
    G(x_k)=G(x_{k+1})=0
\end{equation*}
where it is understood that in $x_0=-\infty$ and $x_{n+1}=+\infty$
this respectively means
\begin{equation*}
    \lim_{x\to-\infty}G(x)=0\quad\qquad\lim_{x\to+\infty}G(x)=0
\end{equation*}
In a dimensionless form, by using $z=\frac{x}{\sigma}$,
$\mu=\frac{\epsilon}{\hbar\omega}=\frac{\lambda}{\omega}$ and
$\chi(z)=G(\sigma z)$, our eigenvalue problem then becomes
\begin{equation}\label{hgeom}
  \chi''(z)-\left({z^2\over4}-{2n+1\over2}-\mu\right)\chi(z) = 0
  \qquad\quad  \chi(z_k)=\chi(z_{k+1})= 0
\end{equation}
where $z_k, z_{k+1}$ are the new dimensionless endpoints. If now
$\mu_m$ and $\chi_m(z)$ are the eigenvalues and eigenfunctions, the
general solution of the \fpe\refeq{statfpe} will be
\begin{equation*}
    f(x,t)=\sum_{m=0}^{\infty}c_m e^{-\mu_m\omega t}
                  \phi_n(x)\chi_m\left({x\over\sigma}\right)
\end{equation*}
where the coefficients $c_m$ will be fixed by the initial conditions
and by the non-negativity and normalization requirements for
$f(x,t)$ along all its evolution. We finally remember that two
linearly independent solutions of the ordinary differential
equation\refeq{hgeom} are
\begin{eqnarray*}
   \chi^{(1)}(z)&=&e^{-z^2/4}M\left(-\,{\mu+n\over2},{1\over2};{z^2\over2}\right)\\
   \chi^{(2)}(z)&=&z\,e^{-z^2/4}M\left(-\,{\mu+n-1\over2},{3\over2};{z^2\over2}\right)
\end{eqnarray*}
where $M(a,b;z)$ are the confluent hypergeometric
functions\mycite{grad}

\subsubsection{The ground state $n=0$}

When $n=0$ the \fpe\refeq{statfpe} takes the form
\begin{equation}\label{fpe0}
    \partial_tf(x,t)=\omega\sigma^2\partial_x^2f(x,t)+
    \partial_x[\omega xf(x,t)]
\end{equation}
while the corresponding \ito\ \sde\refeq{statsde} is
\begin{equation}\label{sde0}
    dX(t)=-\omega X(t)dt+dW(t)
\end{equation}
where $W(t)$ is a Wiener process with diffusion coefficient
$D=\omega\sigma^2$. In both cases we at once recognize an
\emph{\ou}\ process with transition \pdf
\begin{equation}\label{trans0}
    f(x,t|y,s)={1\over\beta(t-s)\sqrt{2\pi}}\,
    {\rm e}^{-[x-\alpha(t-s)]^2/2\beta^2(t-s)}\quad\qquad t\geq s
\end{equation}
where we used the notation
\begin{equation}\label{oucoeff}
    \alpha(t)=ye^{-\omega t}\quad\qquad
   \beta^2(t)=\sigma^2\left(1-e^{-2\omega t}\right)
\end{equation}
This solution of\refeq{fpe0}, which obeys the initial condition
$f(x,s)=\delta(x-y)$, is also the \pdf\ of the solution of the
\sde\refeq{sde0} with initial condition $X(s)=y,\;\Pqo$, namely
\begin{equation*}
   X(t)=ye^{-\omega(t-s)}+\int_s^te^{-\omega(t-t')}dW(t')
\end{equation*}
The stationary process is instead selected by the initial condition
$X(0)\sim\norm(0,\sigma^2)$ namely by the invariant initial \pdf
\begin{equation*}
    f(x,0)={1\over\sigma\sqrt{2\pi}}\,e^{-x^2/2\sigma^2}
\end{equation*}
which is also the asymptotic \pdf\ for every initial condition when
the evolution is ruled by\refeq{trans0}, so that the invariant
distribution plays also the role of the limit distribution. It is
remarkable that this invariant \pdf\ also coincides with the quantum
stationary \pdf\ $\phi_0^2(x)=|\psi_0(x,t)|^2$: in other words, the
process associated by the stochastic mechanics to the ground state
of a \qho\ is a stationary \ou\ process

\subsubsection{The the first excited state $n=1$}

In the case $n=1$ the forward velocity $\fvel_1(x)$ has a
singularity in $x=0$, the \fpe\ is
\begin{equation}\label{fpe1}
    \partial_tf(x,t)=\omega\sigma^2\partial_x^2f(x,t)+
    \partial_x\left[\left(\omega x-\frac{2\omega\sigma^2}{x}\right)f(x,t)\right]
\end{equation}
while the corresponding \ito\ \sde\refeq{statsde} is
\begin{equation}\label{sde1}
    dX(t)=\left(-\omega X(t)+\frac{2\omega\sigma^2}{X(t)}\right)dt+dW(t)
\end{equation}
The \fpe\refeq{fpe1} is then reduced to the eigenvalue
problem\refeq{hgeom} with $x_0=-\infty$, $x_1=0$ and $x_2=+\infty$,
namely (with $z=x/\sigma$)
\begin{equation}\label{hgeom1}
  \chi''(z)-\left({z^2\over4}-{3\over2}-\mu\right)\chi(z) = 0
  \qquad\quad  \chi(-\infty)=\chi(0)=\chi(+\infty)= 0
\end{equation}
This problem should be separately solved for $z\leq0$ and for
$z\geq0$. The eigenvalues turn out to be $\mu_m=2m$ with
$m=0,1,\dots$ and the complete set of eigenfunctions (for both
$z\geq0$, and $z\leq0$) is
\begin{equation*}
    \chi_m(z)=z e^{-z^2/4}M\left(-m,{3\over2};{z^2\over2}\right)=
             {(-1)^m m!\over\sqrt{2}(2m+1)!}\, e^{-z^2/4}
                  H_{2m+1}\left({z\over\sqrt2}\right)
\end{equation*}
In particular it is easy to see that the relation with the quantum
eigenfunction $\phi_1$ is
\begin{equation*}
    \phi_1(x)={\chi_0(x/\sigma)\over\sqrt{\sigma\sqrt{2}}}\qquad\quad \chi_0(z)=z e^{-z^2/4}
\end{equation*}
and that the solution of\refeq{fpe1} for an initial condition
imposed at the time $s$ is
\begin{equation*}
    f(x,t)=\sum_{m=0}^{\infty}c_m e^{-2m\omega(t-s)}
                \phi_1(x)\chi_m\left(\frac{x}{\sigma}\right)\qquad\quad
                t\geq s
\end{equation*}
where the $c_m$'s are fixed by the initial condition. For the
transition \pdf\ from\refeq{degenincond} we have
\begin{equation*}
    c_m={2\over\sqrt{\sigma\sqrt{2\pi}}}\,
{(2m+1)!!\over(2m)!!}\,{\chi_m(y/\sigma)\over \chi_0(y/\sigma)}
\end{equation*}
and by summing up the series\mycite{pla,jpha} we will have with the
notations\refeq{oucoeff}
\begin{equation}\label{trans1}
    f(x,t|y,s)=\Theta(xy)\,{x\over\alpha(t-s)}\,
{ e^{-\frac{[x-\alpha(t-s)]^2}{2\beta^2(t-s)}}-
         e^{-\frac{[x+\alpha(t-s)]^2}{2\beta^2(t-s)}}
           \over \beta(t-s)\sqrt{2\pi}}
\end{equation}
where $\Theta(z)$ is the Heaviside function. In particular we have
\begin{equation*}
    \lim_{t\to+\infty}f(x,t|y,s)=2\Theta(xy)\,
                {x^2 e^{-x^2/2\sigma^2}\over\sigma^3\sqrt{2\pi}}=
                        2\Theta(xy)\phi_1^2(x)
\end{equation*}
and for an arbitrary initial condition $f(x,s^+)=f_0(x)$ we have
\begin{eqnarray*} \lim_{t\to+\infty}f(x,t) &=&
\lim_{t\to+\infty}
                \int_{-\infty}^{+\infty}f(x,t|y,s)f_0(y)\,dy \\
   &=& 2\phi_1^2(x)\int_{-\infty}^{+\infty}\Theta(xy)f_0(y)\,dy
                =\Gamma(q;x)\phi_1^2(x)
\end{eqnarray*}
where we have defined the function
\begin{equation*}
    \Gamma(q;x)=q\Theta(x)+(2-q)\Theta(-x)\,;\qquad
                             q=2\int_0^{+\infty}f_0(y)\,dy
\end{equation*}
Remark that $q=1$ when the initial probability is equally shared on
the two (positive and negative) half-lines, and in this case we have
$\Gamma(1;x)=1$ so that the asymptotical \pdf\ exactly coincides
with the quantum stationary pdf $\phi_1^2(x)$. If on the other hand
$q\not=1$ the asymptotical \pdf\ has the same shape of $\phi_1^2(x)$
but with different weights on the two half-lines. The transition
\pdf\refeq{trans1} is however less elementary of what it looks like
at first sight. It is possible for instance to calculate
expectations, but the results are not very simple: for instance we
have
\begin{eqnarray*}
    \EXP{X(t)\,|\,X(0)=y}&=&\int_0^{\infty}{x^2\over\alpha(t)}\,
{ e^{-\frac{[x-\alpha(t)]^2}{2\beta^2(t)}}-
         e^{-\frac{[x+\alpha(t)]^2}{2\beta^2(t)}}
           \over \beta(t)\sqrt{2\pi}}\\
           &=&
           \sqrt{\frac{2}{\pi}}\beta(t)e^{-\frac{\alpha^2(t)}{2\beta^2(t)}}+\frac{\alpha^2(t)+\beta^2(t)}{\alpha(t)}
           \left[2\Phi\left(\frac{\alpha(t)}{\beta(t)}\right)-1\right]
\end{eqnarray*}
where as usual
\begin{equation*}
    \Phi(z)=\int_{-\infty}^z\frac{e^{-\frac{x^2}{2}}}{\sqrt{2\pi}}\,dx
\end{equation*}
while $\alpha(t)$ and $\beta(t)$ are defined in\refeq{oucoeff}. If
on the other hand $X(0)$ is not degenerate in $y$ but has a \pdf\
$f_0(y)$, the expectation $\EXP{X(t)}$ would be calculated as an
$y$-integral of $\EXP{X(t)\,|\,X(0)=y}\,f_0(y)\,dy$, by recalling of
course that the variable $y$ is hidden into $\alpha(t)$

The transition \pdf\refeq{trans1} anyhow completely solves the
problem of the \fpe\refeq{fpe1}, and hence also completely defines
the law of the Markov process solution of the (non linear)
\sde\refeq{sde1}: it is then puzzling to remark that apparently we
do not know any simple procedure to solve\refeq{sde1} as a \sde. For
instance we could think of changing\refeq{sde1} in some other, more
manageable form by means of a transformation $Y(t)=h(X(t))$ with the
\ito\ formula (see Appendix\myref{theito}). Given the form of our
equation, the simplest idea could seem to be to take $Y(t)=X^2(t)$.
Now\refeq{sde1} is an \ito\ \sde\ with the following coefficients
\begin{equation}\label{ab1}
    a(x)=\fvel_1(x)=-\omega x+\frac{2\omega\sigma^2}{x}\qquad\qquad b(x)=1
\end{equation}
but our transformation $y=h(x)=x^2\geq0$ (with
$h'(x)=2x,\;h''(x)=2$) is only \emph{piecewise monotonic} separately
on the two real half-lines, so that we should separate the process
in two regions according to the sign of $x=g(y)=\pm\sqrt{y}$: a
procedure a bit confusing if done by hand. This of course
corresponds to the fact that the $a(x)$ has a singularity in $x=0$
which effectively separates the two  half-lines. The \ito\ calculus
implies that $Y(t)$ will satisfy a new \sde\ with the
coefficients\refeq{itoy1tind} and\refeq{itoy2tind}, namely (with
$D=\omega\sigma^2$)
\begin{equation*}
    \ha(y)=-2\omega y+6\omega\sigma^2\qquad\qquad\hb(y)=\pm2\sqrt{y}
\end{equation*}
so that we will have one of the two equations for $Y(t)\geq0$
\begin{equation*}
    dY(t)=\left(6\omega\sigma^2-2\omega
    Y(t)\right)dt\pm2\sqrt{Y(t)}dW(t)
\end{equation*}
which in fact are not much easier to handle than the original one,
and also seem not to be on a very firm ground because of the
piecewise monotonicity of $h(x)$

On the other hand there not seems to be any hope of
transforming\refeq{sde1} in a \sde\ with \emph{linear coefficients}
by means of some other clever transformation because the
compatibility conditions\mycite{gs} are not satisfied: in fact,
from\refeq{ab1} and with
\begin{equation*}
    q(x)=\frac{a(x)}{b(x)}-Db'(x)=a(x)=-\omega x+\frac{2\omega\sigma^2}{x}
\end{equation*}
we have
\begin{equation*}
    \frac{1}{q'(x)}\frac{d}{dx}[b(x)q'(x)]=\frac{-4\omega\sigma^2}{\omega x(x^2+2\sigma^2)}
\end{equation*}
which apparently is not $x$-independent, as instead it is required
by the compatibility conditions (see\mycite{gs} p 39). In any case,
since we completely know the law of $X(t)$ through the transition
\pdf\refeq{trans1}, some effort could be produced in the hope of
gaining an insight into the process $X(t)$ solution of\refeq{sde1}
from the fact that\refeq{trans1} turns out to be some kind of
combination of symmetrically separated \ou\ solutions. The proposed
transformation $y=x^2$, however, looks rather preposterous because
$h(x)=x^2$ is not a monotonic function: we will discuss a possible
monotonic modification in the Section\myref{the1} of the present
appendix

\subsubsection{The second excited state $n=2$}

If $n=2$ the velocity $\fvel_2(x)$ has singularities in $x=\pm
\sigma$, the \fpe\ is
\begin{equation}\label{fpe2}
    \partial_tf(x,t)=\omega\sigma^2\partial_x^2f(x,t)+
    \partial_x\left[\left(\omega x-\frac{4\omega\sigma^2x}{x^2-\sigma^2}\right)f(x,t)\right]
\end{equation}
while the corresponding \ito\ \sde\refeq{statsde} is
\begin{equation}\label{sde2}
    dX(t)=\left(-\omega X(t)+\frac{4\omega\sigma^2X(t)}{X^2(t)-\sigma^2}\right)dt+dW(t)
\end{equation}
The \fpe\refeq{fpe2} is then reduced to the eigenvalue
problem\refeq{hgeom} with $x_0=-\infty$, $x_1=-\sigma$, $x_2=\sigma$
and $x_3=+\infty$, namely (with $z=x/\sigma$)
\begin{equation}\label{hgeom1}
  \chi''(z)-\left({z^2\over4}-{5\over2}-\mu\right)\chi(z) = 0
  \qquad\quad  \chi(-\infty)=\chi(\pm1)=\chi(+\infty)= 0
\end{equation}
that should be separately solved in the three intervals
$(-\infty,-1]$, $[-1,1]$ and $[1,+\infty)$. The two linearly
independent solutions are now
\begin{equation*}
    \chi^{(1)}(z)= e^{-z^2/4}M\left(-\,{\mu+2\over2},{1\over2};{z^2\over2}\right)
       \qquad \chi^{(2)}(z)=
z e^{-z^2/4}M\left(-\,{\mu+1\over2},{3\over2};{z^2\over2}\right)
\end{equation*}
and it is easy to check that $\mu=0$ is an eigenvalue for all the
three intervals with eigenfunction
\begin{equation*}
    \chi_0(z)= e^{-z^2/4}M\left(-1,{1\over2};{z^2\over2}\right)=
          e^{-z^2/4}H_2\left({z\over\sqrt{2}}\right)=
         2 e^{-z^2/4}(z^2-1)
\end{equation*}
so that the relation with the quantum eigenfunction now is
\begin{equation*}
    \phi_2(x)={\chi_0(x/\sigma)\over\sqrt{8\sigma\sqrt{2\pi}}}
\end{equation*}
As for the other eigenvalues and eigenfunction they can be obtained
only numerically: for example it can be shown that, beyond
$\mu_0=0$, the first eigenvalues in the interval $[-1,1]$ can be
calculated as the first values such that
\begin{equation*}
    M\left(-\,{\mu+1\over2},{3\over2};{1\over2}\right)=0
\end{equation*}
and are $\mu_1\sim 7.44$, $\mu_2\sim 37.06$, $\mu_3\sim 86.41$. Also
for the unbounded interval $[1,+\infty)$ (the analysis is similar
for $(-\infty,-1]$) the eigenvalues are derivable only numerically.

\subsection{The macroscopic limit}

From the \sde\refeq{statsde} one could hope to derive some
macroscopic, deterministic equation describing the behavior of some
global characteristic of the process. However it is easy to see that
by taking the expectation of\refeq{statsde} we just have
\begin{equation*}
    \frac{d\EXP{X(t)}}{dt}=\EXP{\fvel_n(X(t))}
\end{equation*}
but that we can not deduce any equation for $\EXP{X(t)}$ because the
velocities $\fvel_n(x)$ are not linear functions (with the only
exception of $\fvel_0(x)$), and hence the term $\EXP{\fvel_n(X(t))}$
can not be put in the form $\fvel_n\left(\EXP{X(t)}\right)$. We
could surmise that (at least when the form of $\fvel_n(x)$ is
explicitly given, as for the \qho) some function $h_n(x)$ can be
found such that
\begin{equation*}
     \EXP{\fvel_n(X(t))}=h_n\left(\EXP{X(t)}\right)
\end{equation*}
but this seems not to be an easy task, even for the simplest case of
$\fvel_1(x)$. The outlook seems not to be much brighter if we take
the \emph{medians} of\refeq{statsde} instead of the expectations: as
a matter of fact in this case -- because of the non linearity of the
functional $\MED{\ldots}$ -- we can not even suppose that
$\MED{dX(t)}$ is simply reduced to $d\MED{X(t)}$, and hence the hope
to obtain a differential equation for the median is dashed from the
beginning

The unique viable way to have deterministic equations
from\refeq{statsde} seems then to be to \emph{switch off} the Wiener
noise by taking a vanishing diffusion coefficient $D=0$, while still
keeping a non-zero Planck constant $\hbar$ in order to have
non-trivial values for $\sigma$ in the $\fvel_n(x)$: this means of
course that the relation\refeq{diffus} no longer holds, so that now
in some sense we are out of the framework of the stochastic
mechanics. In any case it could be instructive to see what kind of
deterministic trajectories are solutions of the dynamical systems
\begin{equation}\label{dynsys}
    \dot{x}(t)=\fvel_n(x(t))
\end{equation}
associated to the forward velocities of the stationary states of a
\qho. Since our forward velocities\refeq{qhofvel} are
time-independent functions, the \ode's\refeq{dynsys} can be solved
by separation of the variables for $t\geq0$ with initial condition
$x(0)=y$
\begin{equation}\label{dynsyssol}
  \int_{y/\sigma\sqrt{2}}^{x/\sigma\sqrt{2}}\frac{H_n(z)}{2nH_{n-1}(z)-zH_n(z)}\,dz=\omega t
\end{equation}
but since there is no general formula giving the solutions for every
$n$ we will be obliged to show them one by one: for $n=0$ the
equation is
\begin{equation*}
    \dot{x}(t)=-\omega x(t)
\end{equation*}
and from\refeq{dynsyssol} we have
\begin{equation*}
    \left[\ln z\right]_{y/\sigma\sqrt{2}}^{x/\sigma\sqrt{2}}=-\omega t
\end{equation*}
so that every initial condition will eventually go to $x(+\infty)=0$
according to
\begin{equation}\label{ode0}
    x(t)=ye^{-\omega t}
\end{equation}
For $n=1$ the equation is
\begin{equation*}
    \dot{x}(t)=-\omega x(t)+{2\omega\sigma^2\over x(t)}
\end{equation*}
the solution is
\begin{equation*}
    \left[\ln (z^2-1)\right]_{y/\sigma\sqrt{2}}^{x/\sigma\sqrt{2}}=-2\omega t
\end{equation*}
and hence
\begin{equation}\label{ode1}
    x^2(t)=2\sigma^2+(y^2-2\sigma^2)e^{-2\omega t}
\end{equation}
The trajectory $x(t)$ will then exponentially go from $y$ to
$\pm\sigma\sqrt{2}$ where the sign of the square root will be
decided according to the sign of the initial value $y$: in this
second case we have indeed two attracting points. For $n=2$ on the
other hand the equation is
\begin{equation*}
    \dot{x}(t)=-\omega x(t)+{4\omega\sigma^2 x\over x^2(t)-\sigma^2}
\end{equation*}
and we have
\begin{equation*}
    \left[\ln z(2z^2-5)^2\right]_{y/\sigma\sqrt{2}}^{x/\sigma\sqrt{2}}=-5\,\omega t
\end{equation*}
namely the trajectories are implicitly defined by
\begin{equation}\label{ode2}
    x(t)\left[x^2(t)-5\sigma^2\right]^2=y(y^2-5\sigma^2)^2e^{-5\omega t}
\end{equation}
and, while not easy to be calculated explicitly, they will have now
three asymptotic attracting points in $x=0$ and
$x=\pm\sigma\sqrt{5}$. The solutions are much less elementary for
$n\geq3$: for instance with $n=3$ the trajectories are implicitly
defined by
\begin{equation*}
    \left[\ln (\sqrt{57}+9-4z^2)^{19+\sqrt{57}}(\sqrt{57}-9+4z^2)^{19-\sqrt{57}}\right]_{y/\sigma\sqrt{2}}^{x/\sigma\sqrt{2}}=-76\,\omega t
\end{equation*}
and we are able to find explicitly only the four asymptotic
attracting points
\begin{equation*}
    x=\pm\sigma\sqrt{\frac{9+\sqrt{57}}{2}}=\pm2.8766\,\sigma\qquad\quad x=\pm\sigma\sqrt{\frac{9-\sqrt{57}}{2}}=\pm0.8515\,\sigma
\end{equation*}
Subsequent solutions would grow increasingly complicated and will
not be displayed here

\subsection{Looking deeper into the $n=1$ case}\label{the1}

By looking into the transition \pdf\refeq{trans1} for the $n=1$
eigenstate we see that it is a combination of two \ou\ transition
\pdf's\refeq{trans0} with opposite expectations $\pm\alpha$. This
combination, however, does not qualify as a proper \emph{mixture}
because of the opposite sign of the two terms. Remark in any case
that, these signs notwithstanding, the \pdf\refeq{trans1} is
assembled in such a way that it turns out to be always non-negative
over $\RE$, as it must be for a \pdf. It would be interesting then
to understand the nature of this combination because this could
possibly shed some light on the kind of combination of \ou\
processes constituting a solution of\refeq{sde1}
\begin{prop}\label{prop}
If $p(x)$ is a \pdf\ with finite expectation $\alpha=\int_{\RE}
xp(x)\,dx$, and if
\begin{equation}\label{cond2}
    p(x)\geq p(-x)\qquad\quad\forall x\geq0
\end{equation}
then $\alpha\geq0$, and when $\alpha>0$
\begin{equation}\label{f}
    f(x)=\Theta(x)\frac{x}{\alpha}\,[p(x)-p(-x)]=\left\{
                                                 \begin{array}{cc}
                                                   \frac{x}{\alpha}[p(x)-p(-x)] & \quad x\geq0 \\
                                                   0 & \quad x\leq0
                                                 \end{array}
                                               \right.
\end{equation}
is a \pdf, while
\begin{equation}\label{fbar}
    \overline{f}(x)=\frac{x}{2\alpha}[p(x)-p(-x)]
\end{equation}
is a symmetric \pdf\ in the sense that
$\overline{f}(-x)=\overline{f}(x)$. On the other hand, if $\alpha=0$
we have $p(x)=p(-x)$ (namely $p(x)$ must be symmetric) and
$f,\overline{f}$ must be defined -- if possible -- as a limit for
$\alpha\to0$
\end{prop}
 \indim
First of all we have
\begin{eqnarray*}
\alpha &=& \int_{-\infty}^{\infty} xp(x)\,dx=\int_{-\infty}^0
xp(x)\,dx+\int_0^{\infty}xp(x)\,dx \\
   &=&
   -\int_0^{\infty}xp(-x)\,dx+\int_0^{\infty}xp(x)\,dx=\int_0^{\infty}x[p(x)-p(-x)]\,dx\geq0
\end{eqnarray*}
because all the terms in the integral are non-negative. Then it is
immediate to check that the function
\begin{equation*}
    q(x)=x[p(x)-p(-x)]\qquad\quad x\in\RE
\end{equation*}
is symmetric (namely $q(x)=q(-x)$) and non-negative for every
$x\in\RE$ with $q(0)=0$. Now it is easy to see that
\begin{equation*}
    \int_0^\infty q(x)\,dx=\int_0^\infty xp(x)\,dx-\int_0^\infty
    xp(-x)\,dx=\int_{-\infty}^{\infty} xp(x)\,dx=\alpha
\end{equation*}
and hence also that
\begin{equation*}
    \int_{-\infty}^\infty q(x)\,dx=2\alpha
\end{equation*}
It is apparent then that\refeq{f} is a \pdf\ concentrated on the
positive half-line, while\refeq{fbar} is a symmetric \pdf\ defined
on $\RE$. Finally, since we have seen that
\begin{equation*}
    \alpha=\int_0^{\infty}x[p(x)-p(-x)]\,dx=\int_0^{\infty}q(x)\,dx
\end{equation*}
and $q(x)\geq0$, then $\alpha$ can not be zero unless $q(x)=0$,
namely $p(x)=p(-x)$. In this case $f$ and $\overline{f}$ can not be
defined as\refeq{f} and\refeq{fbar} and we must resort to a limit
for $\alpha\to0$ hoping that the Theorem of l'H\^opital brings it to
a finite result
 \findim

\noindent It is apparent then that the \pdf\refeq{trans1} is a
particular case of\refeq{f} where $p(x)$ is the Gaussian
law\refeq{trans0} of an \ou\ process. It is not clear, instead, what
kind of combination of \rv's -- if any -- admits either $f(x)$ or
$\overline{f}(x)$ as their \pdf's: let us suppose that $X$ is a \rv\
with \pdf\ $p(x)$. Then $p(-x)$ will play the role of the \pdf\ for
$-X$, and the condition\refeq{cond2} could be formulated as
\begin{equation*}
    \PR{X\in B}\geq\PR{X\in-B}
\end{equation*}
where $B\in \mathcal{B}(\RE_+)$ is a Borelian on the positive
half-line $\RE_+$, while we define $-B=\{x\in\RE\,|\,-x\in B\}$. Our
problem then can be formulated as follows: is either $f(x)$, or
$\overline{f}(x)$ the \pdf\ of some combination of $X$ and $-X$?
Could such a combination obey some simpler form of our
\sde\refeq{sde1}? We do not know an answer at this point, but we can
add just a final remark about the expectations of $f$ and
$\overline{f}$: because of the symmetry we immediately have
\begin{equation*}
    \int_{-\infty}^{\infty}x\overline{f}(x)\,dx=0
\end{equation*}
while on the other hand
\begin{eqnarray*}
  \int_{-\infty}^{\infty}xf(x)\,dx &=& \int_0^{\infty}\frac{x^2}{\alpha}[p(x)-p(-x)]\,dx=\int_0^{\infty}\frac{x^2}{\alpha}p(x)\,dx-\int_0^{\infty}\frac{x^2}{\alpha}p(-x)\,dx \\
   &=&\int_0^{\infty}\frac{x^2}{\alpha}p(x)\,dx-\int_{-\infty}^0\frac{x^2}{\alpha}p(x)\,dx=\int_{-\infty}^{\infty}\frac{x|x|}{\alpha}p(x)\,dx
\end{eqnarray*}
In other words the expectation of a \rv\ $Y$ with \pdf\ $f(x)$ seems
to coincide with the expectation of the \rv\ $\frac{1}{\alpha}X|X|$,
if $X$ has the \pdf\ $p(x)$ (namely the \pdf\ which defines $f(x)$
according to\refeq{f} in the Proposition\myref{prop}). Remark
however that, if $X$ has the \pdf\ $p(x)$, $f(x)$ defined
in\refeq{f} would not be the \pdf\ of $Y=X|X|$ which instead, after
a short calculation, would be
\begin{equation*}
   f_Y(y)= \frac{1}{2\sqrt{|y|}}p\left(\frac{|y|}{y}\sqrt{|y|}\right)
\end{equation*}

This discussion, however, in some sense suggests a role for the
transformation $Y(t)=h(X(t))$
\begin{equation*}
    y=h(x)=x|x|=\left\{
                 \begin{array}{ll}
                   x^2 &\qquad x\geq0 \\
                   -x^2 &\qquad x\leq 0
                 \end{array}
               \right.
\end{equation*}
which, at variance with $x^2$, is now a monotonic transformation
with
\begin{equation*}
    h'(x)=2|x|\qquad h''(x)=2\frac{|x|}{x}\qquad x=g(y)=\left\{
                                                          \begin{array}{ll}
                                                            \sqrt{y} & \qquad y\geq0 \\
                                                            -\sqrt{-y} & \qquad y\leq0
                                                          \end{array}
                                                        \right.
\end{equation*}
Now, if $X(t)$ is a solution of the \ito\ \sde\refeq{sde1} with the
coefficients\refeq{ab1}, the \ito\ calculus implies that
$Y(t)=h(X(t))$ will satisfy a new \sde\ with the
coefficients\refeq{itoy1tind} and\refeq{itoy2tind}, namely (with
$D=\omega\sigma^2$)
\begin{equation*}
    \ha(y)=-2\omega y+6\omega\sigma^2\frac{|y|}{y}\qquad\qquad\hb(y)=2\sqrt{|y|}
\end{equation*}
so that we will have
\begin{equation*}
    dY(t)=\left(6\omega\sigma^2\frac{|Y(t)|}{Y(t)}-2\omega
    Y(t)\right)dt+2\sqrt{|Y(t)|}dW(t)
\end{equation*}

\section{Quantiles and medians: a reminder}\label{quantmed}

\subsection{Definitions}\label{def}

The law of a \rv\ (\emph{random variable}) $X$ whatsoever is
characterized by a \cdf\ (\emph{cumulative distribution function})
$F(x)=\PR{X\leq x}$ which is a monotonic, non decreasing function of
$x$ confined between $0$ and $1$, and right-continuous wherever it
jumps. It can also show flat spots where no probability is present.
The \qf\ (\emph{quantile function}) is then usually defined as
\begin{equation}\label{qf}
    Q(p)=\inf\{x\in\RE :p\leq F(x)\}\qquad\quad0\leq
p\leq1
\end{equation}
This results in a well defined, \emph{one-valued} function with
$Q(0)=-\infty$, while $Q(1)=+\infty$ when $F(x)$ only asymptotically
reaches the value $1$. In the case of continuous laws (no jumps),
however, the definition\refeq{qf} can be reduced to
\begin{equation}\label{qf1}
    Q(p)=\inf\{x\in\RE :p= F(x)\}
\end{equation}
and when $F(x)$ is also strictly increasing (no flat spots) we
finally have
\begin{figure}
\begin{center}
\includegraphics*[width=14cm]{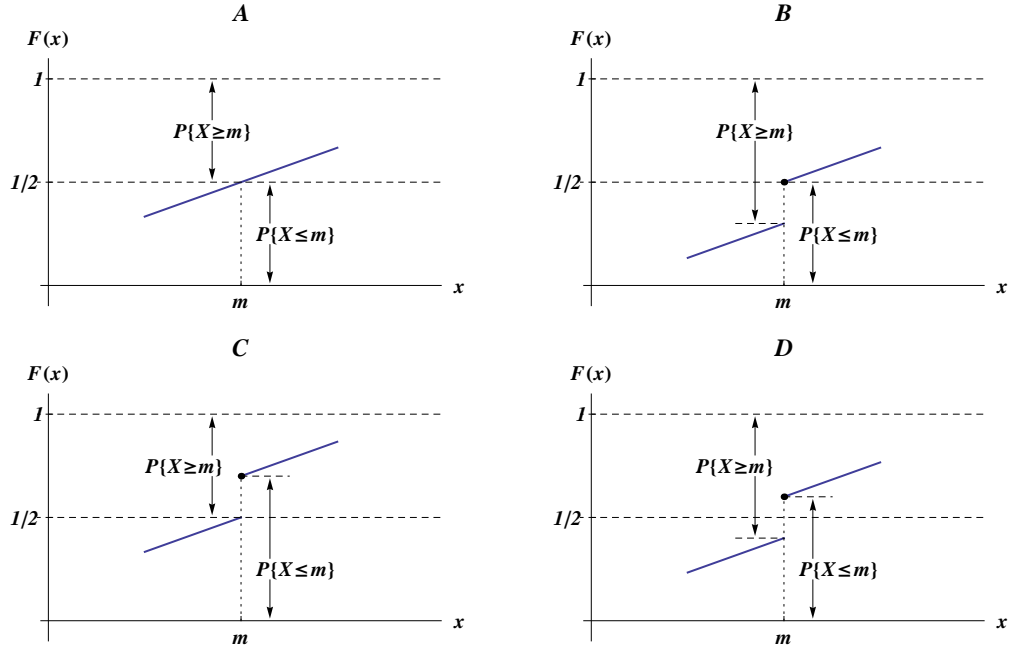}
\caption{Medians according to the Definition\myref{def2} for \cdf's
$F(x)$ with possible jumps, but without flat spots: the value of $m$
is unique and coincides with the $m_X$ of the
Definition\myref{def1}}\label{fig1}
\end{center}
\end{figure}
\begin{equation}\label{qf2}
    Q(p)=F^{-1}(p)
\end{equation}
It is apparent now that $Q(p)$ jumps wherever $F(x)$ has flat spots,
while it has flat spots wherever $F(x)$ jumps. A closer inspection
of the definition\refeq{qf} shows moreover that in its
discontinuities $Q(p)$ is left-continuous. In fact every
non-decreasing, left-continuous function is a possible \qf. An
extensive presentation of both the properties and the statistical
applications of\refeq{qf2} can be found in\mycite{gil}. In the
framework of this notation all the \emph{quantiles} -- and in
particular the \emph{median} -- have a non ambiguous definition
\begin{defn}\label{def1}
The \textbf{median} $\MED{X}$ of a \emph{\rv}\ $X$ is the quantile
of order $\,^1/_2$, namely
\begin{equation}\label{med1}
   m_X=\MED{X}=Q\left(\,\!^1/_2\right)
\end{equation}
\end{defn}
\noindent There is however another, more general definition which
allows for multi-valued medians in the sense that they can also
coincide with a full, closed interval of numbers (the \emph{median
segment}), this definition being of interest mostly when we deal
with the limits of sums of independent random
variables\mycite{loeve,stroock}
\begin{defn}\label{def2}
The \textbf{median} $\MED{X}$ of a \emph{\rv}\ $X$ is any number $m$
such that
\begin{equation}\label{med3}
   \PR{X\leq m}\geq\,^1/_2 \qquad\quad\hbox{and}\qquad\quad\PR{X\geq m}\geq\,^1/_2
\end{equation}
\end{defn}
\noindent It is possible to see in fact that when the median defined
as in\refeq{med3} corresponds to an interval of numbers, then -- due
to the presence of the $\inf\{\ldots\}$ -- the median $m_X$ defined
as in\refeq{med1} coincides with the left endpoint of the said
interval. On the other hand, when the Definition\myref{def2} gives
rise to a unique value, the two definitions apparently coincide.
This shows that in any case $m_X$ of Definition\myref{def1} always
take one (the smallest) of the possible values $m$ of the
Definition\myref{def2}. Examples of applications of the
Definition\myref{def2} are displayed in the Figures\myref{fig1}
and\myref{fig2}, and in particular it is easy to see that, according
to the two definitions, for a Bernoulli \rv\ $X$ taking values $1,0$
with probabilities respectively $p$ and $1-p$ (for $0\leq p\leq1$)
we have
\begin{align*}
   & \hbox{\qquad\qquad\emph{Definition\myref{def1}}} \qquad\qquad\qquad\qquad\qquad\quad
    \hbox{\emph{Definition\myref{def2}}} \\
   & \MED{X}=\left\{
             \begin{array}{ll}
               0 & \quad 0\leq p\leq\frac{1}{2} \\
               1 & \quad \frac{1}{2}<p\leq1
             \end{array}
           \right.
    \qquad\qquad
    \MED{X}=\left\{
                    \begin{array}{cl}
                      0 & \quad 0\leq p<\frac{1}{2} \\
                      $[0,1]$ & \quad p=\frac{1}{2}\\
                      1 & \quad \frac{1}{2}<p\leq1
                    \end{array}
                  \right.
\end{align*}
We finally remember that sometimes in statistics, when the median
segment does not degenerate in a single point, the median is not its
left endpoint but rather some other intermediate point, for example
its middle point: we will not however elaborate more on these
additional possibilities here
\begin{figure}
\begin{center}
\includegraphics*[width=14cm]{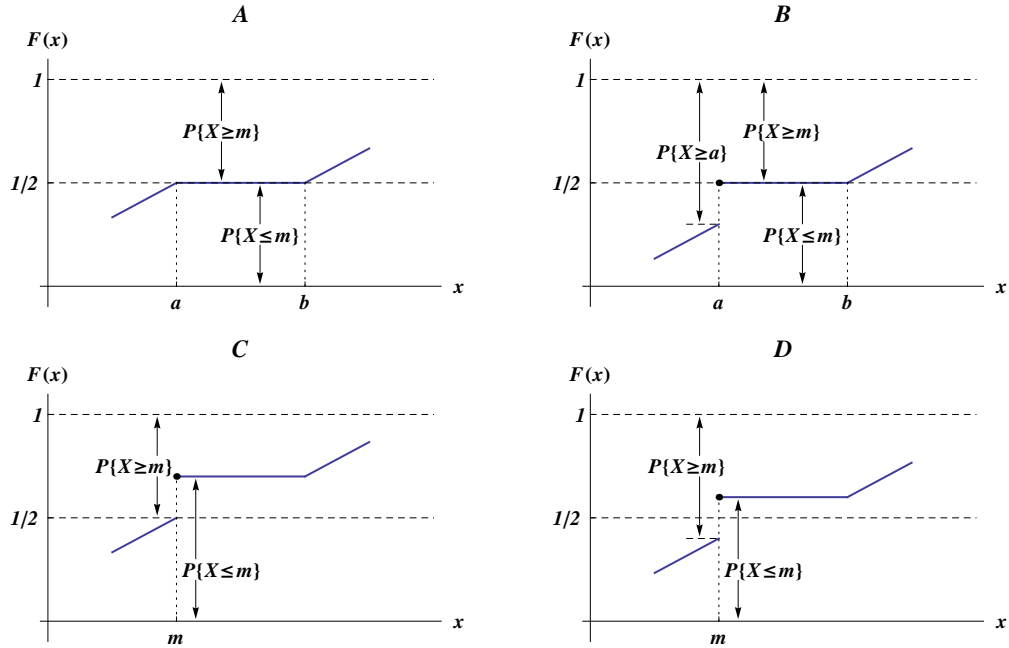}
\caption{Medians according to the Definition\myref{def2} for \cdf's
$F(x)$ with flat spots and possible jumps: the values of $m$ fill a
non degenerate interval $[a,b]$ whenever the flat spot falls at
level $\,^1/_2$ as in $A$ and $B$. In these cases the median $m_X$
of the Definition\myref{def1} coincides with the left endpoint $a$
of the said interval}\label{fig2}
\end{center}
\end{figure}

\subsection{Properties}\label{prop}

\begin{prop}\label{prop1}
Every \emph{\rv}\ $X$ admits a median according to the
Definition\myref{def2}, and if $\,T(x)$ is a monotonic function then
it is
\begin{equation*}
    \MED{T(X)}=T\left(\MED{X}\right)
\end{equation*}
In particular, with $\lambda,\eta\in\RE$, we always have
\begin{equation*}
    \MED{-X}=-\MED{X}\qquad\qquad\quad\MED{\lambda+\eta X}=\lambda+\eta\MED{X}
\end{equation*}
\end{prop}
 \indim
If $m$ is a median for $X$, then from Definition\myref{def2} both
the following inequalities must hold
\begin{equation*}
    \PR{X\leq m}\geq\,^1/_2 \qquad\quad\hbox{and}\qquad\quad\PR{X\geq m}\geq\,^1/_2
\end{equation*}
But with $T(x)$ monotonic these are also equivalent to the pair of
inequalities
\begin{equation*}
    \PR{T(X)\leq T(m)}\geq\,^1/_2 \qquad\quad\hbox{and}\qquad\quad\PR{T(X)\geq T(m)}\geq\,^1/_2
\end{equation*}
so that, always from Definition\myref{def2}, $T(m)$ is a median for
$T(X)$
 \findim

\noindent On the other hand it is neither easy to compute the
medians of a sum in terms of the medians of the summands, nor to
relate the medians of an integrable random variable to its mean
value. Nonetheless many relevant results can be deduced and we will
select here a few among them
\begin{prop}\label{prop2}\label{prop2}
Given a \emph{\rv}\ $X$ with $\EXP{|X|^p}<+\infty$ for some $p\geq1$
(namely endowed at least with the expectation, and possibly also
with higher order moments), for every possible value of $\MED{X}$
according to the Definition\myref{def2}, and for every $a\in\RE$ we
always have
\begin{equation}\label{ineq1}
    \big|\MED{X}-a\big|\leq\big(2\,\EXP{|X-a|^p}\big)^{\frac{1}{p}}
\end{equation}
In particular for $p=2$ and $a=\EXP{X}$ it is straightforward to see
that
\begin{equation*}
    \big|\MED{X}-\EXP{X}\big|\leq\sqrt{2\VAR{X}}
\end{equation*}
\end{prop}
 \indim
From the Chebyshev inequality for every $p\geq1,\;a\in\RE$ and
$\epsilon>0$ we first of all have
\begin{equation*}
    \PR{|X-a|\geq\epsilon}=\PR{|X-a|^p\geq\epsilon^p}\leq\frac{\EXP{|X-a|^p}}{\epsilon^p}
\end{equation*}
so that by taking $\epsilon^p=2\EXP{|X-a|^p}$ we get
\begin{equation*}
    \PR{|X-a|\geq\big(2\EXP{|X-a|^p}\big)^\frac{1}{p}}\leq\frac{1}{2}
\end{equation*}
and since
\begin{eqnarray*}
    \lefteqn{\PR{|X-a|\geq\left(2\EXP{|X-a|^p}\right)^\frac{1}{p}}}\qquad\qquad \\
    &=&\PR{X\geq a+\left(2\EXP{|X-a|^p}\right)^\frac{1}{p}}+\PR{X\leq a-\left(2\EXP{|X-a|^p}\right)^\frac{1}{p}}
\end{eqnarray*}
both the following two inequalities hold simultaneously
\begin{equation*}
    \PR{X\geq a+\left(2\EXP{|X-a|^p}\right)^\frac{1}{p}}\leq\frac{1}{2}\qquad\qquad\PR{X\leq a-\left(2\EXP{|X-a|^p}\right)^\frac{1}{p}}\leq\frac{1}{2}
\end{equation*}
Now it is apparent that
$a-\left(2\EXP{|X-a|^p}\right)^\frac{1}{p}\leq
a+\left(2\EXP{|X-a|^p}\right)^\frac{1}{p}$, so that from the
Definition\myref{def2} we deduce that $\MED{X}$ must fall somewhere
in between these two numbers, namely
\begin{equation*}
    a-\left(2\EXP{|X-a|^p}\right)^\frac{1}{p}\leq\MED{X}\leq
a+\left(2\EXP{|X-a|^p}\right)^\frac{1}{p}
\end{equation*}
and hence\refeq{ineq1} is completely proved
 \findim

\noindent Remark however that there are many, perfectly legitimate,
laws which possess a median, but do not have a well defined
expectation (such as, for example, the Cauchy law), and hence are
not in the framework of the Proposition\myref{prop2} hypotheses: as
a consequence the scope of this result is less wide than it looks
like at first sight.

It is well known that for \rv's with $\EXP{|X|^2}<+\infty$ the
expectation $\EXP{X}$ can be characterized as the value of the
variable $a\in\RE$ that minimizes the mean \emph{square} error
$\EXP{|X-a|^2}$. A similar result holds for the medians of \rv's $X$
with $\EXP{|X|}<+\infty$, but in terms of the mean \emph{absolute}
error
\begin{prop}
Given a \emph{\rv}\ $X$ with $\EXP{|X|}<+\infty$, a number $m$ is a
possible value of the median $\MED{X}$ according to the
Definition\myref{def2} if and only if
\begin{equation*}
    \EXP{|X-m|}=\min_{a\in\RE}\EXP{|X-a|}
\end{equation*}
\end{prop}
 \indim
See\mycite{stroock} p.\ 43
 \findim

\noindent A \rv\ $X$ and its law are said to be
\emph{\textbf{$\mu$-symmetric}} if $\mu-X$ has the same distribution
as $X-\mu$ for some parameter $\mu$. This means that for
$0\,$-symmetric (or \emph{symmetric} tout-court) \rv's we have
\begin{equation*}
    \PR{X\leq x}=\PR{-X\leq x}\qquad\quad\forall\,x\in\RE
\end{equation*}
namely (by changing for convenience the explicit sign of $x$)
\begin{equation*}
    \PR{X\leq -x}=\PR{-X\leq -x}=\PR{X\geq x}=1-\PR{X<x}\qquad\quad\forall\,x\in\RE
\end{equation*}
or, in terms of the \cdf\ $F(x)$,
\begin{equation}\label{symm}
  F(-x)+F(x^-)=1\qquad\quad\forall\,x\in\RE
\end{equation}
%Remark that an arbitrary \rv\ $X$ can always be \emph{symmetrized}
%by defining a new \rv\ $\widetilde{X}$
%\begin{equation}\label{symm3}
%    \widetilde{X}=X-X'
%\end{equation}
%where $X$ and $X'$ are \iid\ (\emph{independent and identically
%distributed}).
In particular for $x=0$ the equation\refeq{symm} implies that
\begin{equation*}
    F(0)+F(0^-)=1
\end{equation*}
and hence either $F(0)=\,\!^1/_2$ (when $F(x)$ is continuous in
$x=0$), or $F(0^-)$ and $F(0)$ are symmetrically located around the
central value $\,\!^1/_2$ (when $F(x)$ jumps in $x=0$). Obviously,
when the expectations exist, we have $\EXP{X}=0$ for every symmetric
\rv, and also the most natural choice for the median of a symmetric
random variable seems to be 0, but in this case some qualification
is in order
\begin{prop}
Given a symmetric \emph{\rv}\ $X$, either $\MED{X}=0$ if the median
is unique, or the non-degenerate median segment is $[-a,a]$ for some
suitable $a>0$
\end{prop}
 \indim
When $F(x)$ is continuous in $x=0$ (without being constant in a
neighborhood of it) then we have seen that $F(0)=\,\!^1/_2$ and
hence the median has the unique value $\MED{X}=0$. The same result
holds when $F(x)$ jumps in $x=0$ because in this case we know that
$F(0-)<\,\!^1/_2$ and $F(0)>\,\!^1/_2$. On the other hand when
$F(x)$ in continuous and constantly takes the value $\,\!^1/_2$ in a
neighborhood of $x=0$, then we have a non-degenerate median segment
which must be symmetric around $x=0$. If indeed $m>0$ belongs to
this segment, also $-m$ must be a median value because
\begin{equation*}
    F(-m)=1-F(m^-)=1-F(m)=1-\,\!^1/_2=\,\!^1/_2
\end{equation*}
Then the median segment apparently is $[-a,a]$ for some suitable
$a>0$
 \findim

\noindent When a \cdf\ $F(x)$ is continuous without flat spots on
all its support the median is always uniquely defined, the two
Definitions\myref{def1} and\myref{def2} give rise to the same value,
and we can simply take $Q\left(\,\!^1/_2\right)$ as the median. In
this case the \qf\ $Q(p)$ has many other properties that can be
found for example in\mycite{gil}. In particular it is possible to
describe -- by means of suitable parameters -- entire \emph{types}
of \qf's starting from some \emph{basic form} $S(p)$: if for
instance $S(p)$ is the \qf\ of some \rv\ $X$, then it is easy to see
that all the functions
\begin{equation*}
    Q(p)=\lambda+\eta S(p)
\end{equation*}
with $\lambda\in\RE$ and $\eta>0$, would be good \qf's for the \rv's
the type $\lambda+\eta X$ spanned by $X$. If in fact $F(x)$ and
$G(x)$ respectively are the \cdf's of $Q(p)$ and $S(p)$, it is easy
to see that from the previous equation we also have
\begin{equation*}
    F(x)=G\left(\frac{x-\lambda}{\eta}\right)
\end{equation*}
namely $F(x)$ belongs to the type spanned by $G(z)$: if indeed $X$
is distributed as $G(x)$, then $\lambda+\eta\,X$ is distributed as
$F(x)$. This property however can be put in a more general form as
follows
\begin{prop}
If $Q_X(p)$ is the \qf\ of the continuous \rv\ $X$, and $Z=T(X)$ we
have
\begin{equation*}
    Q_Z(p)=\left\{
             \begin{array}{ll}
               T\left(Q_X(p)\right) &\quad \hbox{if $T(x)$ is continuous, monotonic increasing} \\
               T\left(Q_X(1-p)\right) &\quad \hbox{if $T(x)$ is continuous, monotonic decreasing}
             \end{array}
           \right.
\end{equation*}
\end{prop}
 \indim
When $T(x)$ is increasing we have
\begin{equation*}
  p=F_Z(z) = \PR{Z\leq z}=\PR{T(X)\leq z}=\PR{X\leq T^{-1}(z)}=F_X\left(T^{-1}(z)\right)
\end{equation*}
and hence $ T^{-1}(z)=Q_X(p)$, namely $z=T\left(Q_X(p)\right)$, so
that finally $Q_Z(p)=z=T\left(Q_X(p)\right)$. If on the other hand
$T(x)$ is decreasing, being our \rv's continuous we find that
\begin{equation*}
  p=F_Z(z) = \PR{Z\leq z}=\PR{T(X)\leq z}=\PR{X\geq T^{-1}(z)}=1-F_X\left(T^{-1}(z)\right)
\end{equation*}
and hence, as before, we get $Q_Z(p)=T\left(Q_X(1-p)\right)$
 \findim

\noindent Remark finally that, when $T(x)$ is simply monotonic, from
both the previous results we get
\begin{equation*}
    Q_Z\left(\,\!^1/_2\right)=T\big(Q_X\left(\,\!^1/_2\right)\big)
\end{equation*}
so that
\begin{equation*}
    \MED{T(X)}=\MED{Z}=m_Z=T(m_X)=T\left(\MED{X}\right)
\end{equation*}
namely we again obtain the result about the medians stated in the
Proposition\myref{prop1}

\subsection{Expectations and medians for symmetric laws}\label{medians}

To avoid possible ambiguities in this section we will confine
ourselves to \emph{continuous} (namely \emph{absolutely continuous})
laws equipped with a non vanishing \pdf. Generally speaking
expectations and medians do not coincide, except in particular
cases: for instance, the \pdf\ $f(x)$ of a $\mu$-symmetric \rv\ $X$
is an even function around the parameter $\mu$
\begin{equation*}
    f(\mu+x)=f(\mu-x)
\end{equation*}
then (if the expectation exists) it is easy to see that
$\EXP{X}=\MED{X}=\mu$. We know on the other hand from the
Proposition\myref{prop1} that the medians show a particular property
not shared with the expectations: when $y=T(x)$ is a monotonic
function defined on the set of values of the \rv\ $X$ and we define
$Y=T(X)$, we get
\begin{equation}\label{medianfunct}
    \MED{T(X)}=T\left(\MED{X}\right)
\end{equation}
It is apparent then that if $X$ is $\mu$-symmetric we also have
\begin{equation}\label{meanmed}
    \MED{T(X)}=T\left(\EXP{X}\right)
\end{equation}
and in particular this is true for the \textbf{Gaussian \rv's
$X\sim\norm(\mu,\sigma^2)$} that are famously $\mu$-symmetric and
will be briefly discussed herein

If, to begin with, $X$ is a \emph{dimensionless} Gaussian \rv, then
$Y= e^X\sim\lnorm(\mu,\sigma^2)$ is a dimensionless log-normal, and
hence we have
\begin{align}
  & \qquad\qquad\quad \EXP{X}=\MED{X}=\mu \qquad\qquad\VAR{X}=\sigma^2 \label{norm}\\
  &   \EXP{Y}=e^{\mu+\frac{\sigma^2}{2}} \qquad\qquad\MED{Y}=e^\mu\qquad\qquad
    \VAR{Y}=e^{2\mu+\sigma^2}\left(e^{\sigma^2}-1\right)\label{lnorm}
\end{align}
In particular, in agreement with\refeq{meanmed}, we have the
following relations
\begin{equation}\label{med}
   \EXP{\ln Y}=\EXP{X}=\ln\MED{Y}\qquad\qquad\MED{e^X}=\MED{Y}=e^{\EXP{X}}
\end{equation}
that are instrumental in the discussion of the present paper, and
then will deserve a few additional remarks

First of all the equations\refeq{med} apparently hold     for
dimensionless Gaussian \rv's $X$ and for their log-normal
counterparts $Y=e^X$, but a dimensionally complete formulation is
always possible: let us suppose that our Gaussian
$X\sim\norm(\mu,\sigma^2)$ is no longer dimensionless, but it is for
instance a length. In this case, by taking advantage of the standard
deviation $\sigma$ that is a length too, we remark first that
$\,\!^X/_\sigma\sim\norm(\,\!^\mu/_\sigma,1)$ is now dimensionless,
and then that $Z=e^{X/\sigma}\sim\lnorm(\,\!^\mu/_\sigma,1)$ is a
dimensionless log-normal with
    \begin{equation*}
        \EXP{Z}=e^{\frac{\mu}{\sigma}+\frac{1}{2}}\qquad\quad\MED{Z}=e^{\frac{\mu}{\sigma}}\qquad\quad\VAR{Z}=e^{\frac{2\mu}{\sigma}+1}\left(e-1\right)
    \end{equation*}
Going then to the dimensional variable $Y=\sigma Z=\sigma
e^{X/\sigma}$, we at once have
    \begin{equation*}
        \EXP{Y}=\sigma e^{\frac{\mu}{\sigma}+\frac{1}{2}}\qquad\qquad\VAR{Y}=\sigma^2 e^{\frac{2\mu}{\sigma}+1}\left(e-1\right)
    \end{equation*}
while from the Proposition\myref{prop1} for the median we have
    \begin{equation*}
       \MED{Y}=\sigma\MED{Z}=\sigma e^{\frac{\mu}{\sigma}}
    \end{equation*}
In particular  the relations\refeq{med} are accordingly changed into
    \begin{equation}\label{med2}
       \EXP{\sigma\ln \frac{Y}{\sigma}}=\EXP{X}=\sigma\ln\frac{\MED{Y}}{\sigma}\qquad\qquad\MED{\sigma e^{\frac{X}{\sigma}}}=\MED{Y}=\sigma e^{\frac{\EXP{X}}{\sigma}}
    \end{equation}
Remark that these results should be suitably adjusted when dealing
with a Gaussian \emph{process} $X(t)$ because its variance could
possibly be time-dependent: in this case it would be appropriate to
find some other constant parameter to play the role of $\sigma$

We must remark  moreover that, if face the problem of writing down
some deterministic evolution either as expectation, or as median of
a process, it is possible to consider several alternatives by taking
for instance into account also the variances. For example we found
that the expectation of $Y=\sigma e^{X/\sigma}$
    \begin{equation}\label{mean}
        \EXP{Y}=\sigma e^{\frac{\mu}{\sigma}+\frac{1}{2}}
    \end{equation}
is a function of the expectation $\mu$ and of the variance
$\sigma^2$ of the original normal \rv\ $X$: as a consequence, since
for the processes discussed in the present paper, both $\EXP{X(t)}$
and $\VAR{X(t)}$ are explicitly known, for our required evolutions
we can always resort to the expectations rather than to the medians

Of course the relations\refeq{med}, \refeq{med2} and\refeq{mean}
perfectly agree with\refeq{medianfunct} drawn from the
Proposition\myref{prop1} but for the fact that here, being $X$ a
Gaussian \rv, we also take advantage of the relation
$\EXP{X}=\MED{X}$ which only holds for $\mu$-symmetric distributions
as discussed at the beginning of this section: for example with $X$
arbitrarily distributed and $Y=\sigma e^{X/\sigma}$, we always have
\begin{equation}\label{medianexp}
    \MED{Y}=\MED{\sigma e^{X/\sigma}}=\sigma e^{\frac{\MED{X}}{\sigma}}
\end{equation}
while\refeq{med2} only holds for a Gaussian $X$

\end{appendix}

\end{document}